\documentclass[10pt]{amsart}
\usepackage{amsmath,amscd}
\usepackage{amssymb,xypic,graphicx}
\usepackage{amsthm,verbatim}
\usepackage{xcolor}
\usepackage[all,color]{xy}
\usepackage{hyperref}
\usepackage{tikz}
\usepackage{tikz-cd}
\usetikzlibrary{matrix,arrows,decorations.pathmorphing}
\newcommand*{\rom}[1]{\expandafter\@slowromancap\romannumeral #1@}

\newcommand{\fkg}{{\mathfrak{{g}}}}

\newcommand{\KK}{{\mathcal{K}}}

\newcommand{\BB}{{\mathcal{B}}}

\newcommand{\CC}{{\mathcal{C}}}
\newcommand{\Mat}{\operatorname{Mat}}

\newcommand{\im}{\operatorname{im}}

\numberwithin{equation}{subsection}

\newtheorem{thm}{Theorem}[section]
\newtheorem{prop}[thm]{Proposition}
\newtheorem{lem}[thm]{Lemma}

{  \theoremstyle{definition}
\newtheorem{defi}[thm]{Definition}

\newtheorem{rem}[thm]{Remark}

}

\newcommand{\Pf}{\noindent {\it Proof}}
\newcommand{\id}{\operatorname{id}}

\newcommand{\ra}{\rightarrow}
\renewcommand{\AA}{{\mathcal A}}

\newcommand{\UU}{{\mathcal U}}

\newcommand{\Om}{\Omega}

\newcommand{\Hom}{\operatorname{Hom}}

\renewcommand{\a}{\alpha}

\newcommand{\R}{{\Bbb R}}
\newcommand{\Z}{{\Bbb Z}}

\newcommand{\ot}{\otimes}

\newcommand{\ed}{\qed\vspace{3mm}}

\title{Homotopy L-infinity spaces and Kuranishi manifolds, I: categorical structures }
\author{Junwu Tu}
\begin{document}
\begin{abstract}
Motivated by the definition of homotopy $L_\infty$ spaces, we develop a new theory of Kuranishi manifolds, closely related to Joyce's recent theory. We prove that Kuranishi manifolds form a $2$-category with invertible $2$-morphisms, and that certain fiber product property holds in this $2$-category. In a subsequent paper, we construct the virtual fundamental cycle of a compact oriented Kuranishi manifold, and prove some of its basic properties.

Manifest from this new formulation is the fact that $[0,1]$-type homotopy $L_\infty$ spaces are naturally Kuranishi manifolds. The former structured spaces naturally appear as derived enhancements of Maurer-Cartan moduli spaces from Chern-Simons type gauge theory. In this way, Kuranishi manifolds theory can be applied to study path integrals in such type of gauge theories.
\end{abstract}

\maketitle

\section{Introduction}\label{intro-sec}

\subsection{Main results and literature} Global Kuranishi theory was pioneered by Fukaya-Oh-Ohta-Ono's seminal works~\cite{FO},~\cite{FOOO}. There are now quite a few similar, yet different versions of global Kuranishi theory~\cite{Joyce},~\cite{Joyce2},~\cite{MW},~\cite{MW2},~\cite{Pardon},~\cite{Yang}. We refer to Joyce's papers~\cite{Joyce},~\cite{Joyce2} for a thorough discussion of the different definitions. In this paper, we propose yet another new definition of Kuranishi manifolds (with trivial isotropy group).  We prove that they form a $2$-category $\mathfrak{Kur}$ with invertible $2$-morphisms, and that it holds certain fiber product property. The proof uses a $2$-localization construction of $2$-categories. This new Kuranishi manifold theory is closely related to Joyce's version where similar results were obtained. For this reason, the author considers the current theory no more than a simplified version of Joyce's theory. Although, the $2$-category $\mathfrak{Kur}$ presented in this paper is not homotopy equivalent to Joyce's version. We refer to Remark~\ref{rem:joyce} for a detailed discussion of the differences. It is particularly interesting how we arrived at this definition: by writing down an existing definition in the theory of $L_\infty$ spaces~\cite{Cos},~\cite{Tu}.  More precisely, we consider $[0,1]$-type $L_\infty$ spaces, those that have tangent complexes concentrated at degrees $0$ and $1$. From this $L_\infty$ point of view, Joyce's definition is extremely natural. The relationship of the three classes of structured spaces are illustrated in the following diagram~\footnote{The horizontal implication was pointed out by Joyce to the author.}.
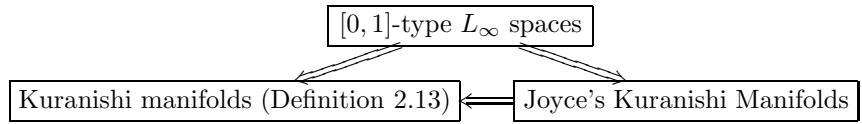
\[\begin{xy}
(0,0)*=<100pt,15pt>[F]{\mbox{$[0,1]$-type $L_\infty$ spaces}}*\frm{-}="A";
(-30,-10)*=<170pt,15pt>[F]{\mbox{Kuranishi manifolds (Definition~\ref{def:kuranishi})}}="B";
(30,-10)*=<130pt,15pt>[F]{\mbox{Joyce's Kuranishi Manifolds}}="C";
{\ar@{=>} "A"; "B"};
{\ar@{=>} "A"; "C"};
{\ar@{=>} "C"; "B"};
\end{xy}\]

\subsection{A new look at Kuranishi charts} A small neighborhood $U$ of a point in a moduli space of interest can often be described via a homeomorphism
\[\psi: U \ra s^{-1}(0)\]
onto the zero locus of a map
\[ s: V \ra \R^m\]
where $V\subset\R^n$ is an open neighborhood of the origin, and $s(0)=0$. Denote by $\underline{\R^m}:=C^\infty_V\ot_\R \R^m$
the trivial vector bundle on $V$ generated by $\R^m$. We may consider $s$  as a section  of $\underline{\R^m}$. The quadruple $(V,\underline{\R^m}, s,\psi)$ is called a {\sl Kuranishi chart} of $U$. 

Let $(V,\underline{\R^m},s,\psi)$ be a Kuranishi chart of a topological space $U$.  The infinite prolongation 
\[ \widetilde{s} \in \Gamma(V, J(\underline{\R^m}))\]
of $s$ gives a section in the infinite jet bundle of $\underline{\R^m}$. Since $\underline{\R^m}$ is by definition trivialized, we have an isomorphism $J(\underline{\R^m})\cong J(\underline{\R}) \ot_\R \R^m$. Furthermore, since $V\subset \R^n$ is an open subset of a vector space, the linear structure induces a torsion-free (and flat) connection on its tangent bunle $T_V$, which further induces a canonical trivialization $ J(\underline{\R})\cong \hat{S}_{C^\infty_V} \Om_V$. Putting the two isomorphisms together we get
\[ J(\underline{\R^m})\cong (\hat{S}_{C^\infty_V}\Om_V) \ot_\R \R^m \cong \prod_{j}^\infty\Hom_{C^\infty_V} ( S^j T_V, \underline{\R^m} ).\]
We denote by $\prod_j^\infty \widetilde{s}_j$ the image of $\widetilde{s}$ under this isomorphism. 

\begin{lem}
The morphisms $\widetilde{s}_j$ with $j\geq 0$ defines a curved $L_\infty$ algebra structure on the graded bundle
\[ \fkg:= T_V[-1]\oplus \underline{\R^m}[-2].\]
\end{lem}

\Pf. Simply by degree reason, any collection of morphisms from $S^j T_V \ra \underline{\R^m}$ defines a $L_\infty$ algebra.\ed

Let $(C^*\fkg, Q)$ be the Chevalley-Eilenberg algebra of this curved $L_\infty$ algebra bundle. To characterize the extra property that the morphisms $\widetilde{s}_j, (j\geq 0)$ are determined by only the first map $\widetilde{s}_0=s$, 
the author introduced a flat connection $D$ on $(C^*\fkg,Q)$ in~\cite{Tu}. The triple $(V,\fkg,D)$ is called a $L_\infty$ space. 

\medskip
\begin{lem}~\label{intro-space}
Let $V\subset \R^n$, and set $\fkg:= T_V[-1]\oplus \underline{\R^m}[-2]$ as before. Then to give a $L_\infty$ space $(V,\fkg,D)$ with $D$ extending the flat structures on $\fkg$, is equivalent to the data of a section $s\in \underline{\R^m}$.
\end{lem}

\Pf. Starting from the data of a section $s\in \underline{\R^m}$, we have seen how to construct the $L_\infty$ algebra structure on the bundle $\fkg$. For the definition of the connection $D$, see~\cite{Tu}. In the reverse direction, the section $s$ is simply the curvature term of the $L_\infty$ algebra $\fkg$.\ed

Let $U$ be a topological space. An $L_\infty$ enhancement of $U$ is a quadruple $(V,\fkg, D, \psi)$ where $(V,\fkg,D)$ is a $L_\infty$ space of the form in Lemma~\ref{intro-space}, and $\psi: U \ra s^{-1}(0)\subset V$ is a homeomorphism onto the zero locus of $s$ (the curvature term of $\fkg$). By Lemma~\ref{intro-space}, to give a $L_\infty$ enhancement of $U$ is equivalent to give a Kuranishi chart of $U$. This $L_\infty$ point of view of a section of a bundle, lies at the heart of the paper: it suggests a homotopy theoretic approach to define Kuranishi structures.

\subsection{Morphisms and homotopies between morphisms} The first obstacle in obtaining a global Kuranishi theory is the delicate question of how to define coordinate changes between Kuranishi charts. 

So what if we think about this question in terms of $L_\infty$ spaces? The notion of homomorphisms between two $L_\infty$ spaces, as well as the notion of homotopies between homomorphisms
have already been worked out in~\cite{Tu}. The work is to translate them into the languages of Kuranishi theory. Surprisingly, this $L_\infty$ approach immediately brings us to the frontier of the current research on Kuranishi manifolds (closely related to Joyce's work~\cite{Joyce}).

Recall in~\cite{Tu} a morphism $$(f,f^\sharp): (V_\alpha, \fkg_\alpha, D_\alpha) \ra (V_\beta,\fkg_\beta,D_\beta)$$ consists of a smooth map
\[ f: V_\alpha \ra V_\beta,\]
together with a $L_\infty$ morphism
\[ f^\sharp: \fkg_\alpha \ra f^*\fkg_\beta\]
that is compatible with the flat connections $D_\alpha$ and $f^*D_\beta$. Let $(V_\alpha, \fkg_\alpha, D_\alpha,\psi_\alpha)$ and $(V_\beta,\fkg_\beta,D_\beta,\psi_\beta)$ be two $L_\infty$ enhancements of a topological space $U$. A morphism of enhancements is simply a morphism $(f,f^\sharp)$ of $L_\infty$ spaces such that 
\[  f \circ \psi_\alpha =\psi_\beta.\]

\medskip
\begin{lem}~\label{intro-mor}
Let the two $L_\infty$ spaces $(V_\alpha, \fkg_\alpha, D_\alpha)$,  $(V_\beta,\fkg_\beta,D_\beta)$ be as in Lemma~\ref{intro-space}. Then the data of a morphism $(f,f^\sharp): (V_\alpha, \fkg_\alpha, D_\alpha) \ra (V_\beta,\fkg_\beta,D_\beta)$ between $L_\infty$ spaces is equivalent to the data of a map $f: V_\alpha \ra V_\beta$, and a bundle map $\hat{f}: \underline{\R^{m_\alpha}} \ra f^*\underline{\R^{m_\beta}}$ such that 
\[ f^*s_\beta=\hat{f} s_\alpha,\]
where $s_\alpha$ and $s_\beta$ are the two curvature terms of $\fkg_\alpha$ and $\fkg_\beta$.
\end{lem}

Fixing a topological space $U$, the above lemma implies that the category of $L_\infty$ enhancements of $U$ of the form $(V,\fkg,D,\psi)$ with $D$ extends the canonical flat structures on $\fkg$, is equivalent to the category of Kuranishi charts of $U$. The question of defining a notion of coordinate changes between Kuranishi charts boils down to find a suitable notion of equivalences in either one of these two equivalent categories. 

\medskip
The main advantage of the $L_\infty$ interpretation is that it suggests a natural solution to this question: homotopy equivalences between $L_\infty$ enhancements. 
Indeed, given two morphisms
\[ (f_0,f_0^\sharp), (f_1,f_1^\sharp): (V_\alpha, \fkg_\alpha, D_\alpha,\psi_\alpha)\ra (V_\beta,\fkg_\beta,D_\beta,\psi_\beta)\]
between $L_\infty$ enhancements, the notion of a homotopy between them was introduced in~\cite[Section 4.4]{Tu}. The following lemma is a geometric reinterpretation of such a homotopy. The proof is omitted, as we only use the lemma to motivate Definition~\ref{def:homotopy} in the next section.
\begin{comment}
Recall as defined in~\cite{Tu}, two morphisms
\[ (f_0,f_0^\sharp), (f_1,f_1^\sharp): (V_\alpha, \fkg_\alpha, D_\alpha,\psi_\alpha)\ra (V_\beta,\fkg_\beta,D_\beta,\psi_\beta)\]
between $L_\infty$ enhancements are called homotopic if there exist
\begin{itemize}
\item [(1.)]  A map
\[ F: V_\alpha\times \Delta^1 \ra V_\beta,\]
such that $F|_{V_\alpha\times 0}=f_0$, $F|_{V_\alpha\times 1}=f_1$ and that the following diagram is commutative.

\medskip
\item[(2.)] A morphism of $L_\infty$ algebra bundles
\[ F^\sharp+G^\sharp dt: \pi_1^*\fkg_\alpha \ra F^*\fkg_\beta\oplus F^*\fkg_\beta\cdot dt,\]
which is compatible with the two connections $D_\alpha$ and $D_\beta$. Furthermore, we have $$i_0^*F^\sharp=f_0^\sharp, \mbox{\;\; and \;\;} i_1^*F^\sharp=f_1^\sharp$$
where $i_0, i_1: V_\alpha \ra V_\alpha\times \Delta^1$ are inclusions at $t=0, 1$ respectively.
\end{itemize}

\end{comment}
\medskip
\begin{lem}~\label{intro-homotopy}
Let $(f_0,f_0^\sharp)$ and $(f_1,f_1^\sharp)$ be two morphisms between $L_\infty$ enhancements, and let $(f_0,\hat{f}_0)$, $(f_1,\hat{f}_1)$ be the associated morphism between Kuranishi charts. Then a homotopy in the sense of~\cite[Section 4.4]{Tu} is equivalent to the following data:
\begin{itemize}
\item[(a.)] A family of map
\[ F(t): V_\alpha \ra V_\beta,\; t\in [0,1]\]
such that $F(0)=f_0$, $F(1)=f_1$.
\begin{comment} and that
\[ F(t)\psi_\alpha = \psi_\beta,\; \forall t\in [0,1].\]\end{comment}
\item[(b.)] A family of bundle map $\widehat{F}(t): \underline{\R^{m_\alpha}} \ra F(t)^*\underline{\R^{m_\beta}}$ such that $\hat{F}(0)=\hat{f}_0$, $\widehat{F}(t)=\hat{f}_1$, and satisfy the equation
\[ F(t)^*s_\beta=\hat{F}(t) s_\alpha, \; \forall t\in [0,1].\]
\item[(c.)] Two families of bundle maps
\begin{align*}
\Lambda(t):& \underline{\R^{m_\alpha}} \ra F(t)^* \underline{\R^{n_\beta}},\\
\Xi(t):& \wedge^2\underline{\R^{m_\alpha}} \ra F(t)^*\underline{\R^{m_\beta}}.
\end{align*}
such that
\begin{align*}
\frac{dF(t)}{dt}|_x &= \Lambda(t)(s_\alpha(x)), \;\; \forall x\in V_\alpha\\
\frac{d\hat{F}(t)}{dt}|_{(x,\xi)} &= ds_\beta|_{F(t)(x)}\big(\Lambda(t)(x,\xi)\big)+\Xi(t)(s_\alpha(x)\wedge \xi),\;\; \forall (x,\xi)\in \underline{\R^{m_\alpha}}.
\end{align*}
Here in the second equation, $\frac{d\hat{F}(t)}{dt}|_{(x,\xi)}$ stands for the derivative valued in the fiber direction, and $ds_\beta$ is computed using the trivialization by considering $s_\beta$ as a map $V_\beta \ra \R^{m_\beta}$.
\end{itemize}
\end{lem}

The actual Definition~\ref{def:homotopy} we shall use is not given by Conditions $(a.)$,$(b.)$,$(c.)$, as that turned out to be unnecessarily complicated. We shall work out some corollaries from these conditions with which we found most convenient to develop a global Kuranishi theory. Finally, we note that the Conditions above imply Joyce's equivalence relation among Kuranishi morphisms. This shows that $[0,1]$-type $L_\infty$ spaces are automatically Kuranishi manifolds in Joyce's sense.

\subsection{Contents of the paper} In Section~\ref{sec:kur}, we define the notion of a Kuranishi manifold. In Section~\ref{sec:cat}, we define morphisms between Kuranishi manifolds and $2$-morphisms between morphisms. Then we  prove the main result (Theorem~\ref{thm:main}) that Kuranishi manifolds form a $2$-category $\mathfrak{Kur}$ with invertible $2$-morphisms. Section~\ref{sec:fiber} deals with certain $2$-fiber product property in $\mathfrak{Kur}$ (Theorem~\ref{thm:fiber}). The Appendix~\ref{app} contains technical proofs. 

\subsection{Notations and Conventions} The following are used throughout the paper.

\medskip
\noindent {\bf (A.)} Over a topological space $V_\beta$, we use the notation $\underline{\R^m}$ for the trivial bundle of rank $m$. If $f: V_\alpha\ra V_\beta$ is a continuous map, the pull-back $f^*\underline{\R^{m}}$ is still a trivial bundle, over $V_\alpha$ instead of $V_\beta$. We shall abuse the notation $\underline{\R^{m}}$ to stand for the trivial bundle of rank $m$ over possibly different topological spaces. We shall write down the underlying space whenever confusion may occur.

\medskip
\noindent {\bf (B.)} In a $2$-category $\CC$, we often form the horizontal composition $\lambda\circ_1 \eta: fh\ra gk$ between $2$-morphisms $\lambda: f\ra g$ and $\eta: h\ra k$. In the special case when $f=g$, and $\lambda=\id_f$, we use the notation $f\circ_1\eta$ instead of $\id_f\circ_1\eta$. Similarly, we also use $\lambda\circ_1 h$ instead of $\lambda\circ_1 \id_h$.

\medskip
\noindent {\bf (C.)} Throughout the paper, we continue to adopt Joyce's convention to not write down explicitly the domains of Kuranishi morphisms, since the actual domain is clear from the indices used. For example, in Definition~\ref{def:kuranishi}, the notation
\[ f_{ij}: [V_i,\underline{\R^{m_i}}, s_i, \psi_i] \ra [V_j,\underline{\R^{m_j}}, s_j, \psi_j]\]
means a germ of Kuranishi morphism $[f_{ij}, \hat{f}_{ij}]: [V_{\underline{i}j}, \underline{\R^{m_i}}, s_i, \psi_i] \ra [V_{i\underline{j}}, \underline{\R^{m_j}}, s_j, \psi_j]$ of the restricted Kuranishi charts, with $V_{\underline{i}j}\subset V_i$, $V_{i\underline{j}}\subset V_j$ such that $\psi_i^{-1}(V_{\underline{i}j})=U_{ij}$, and $\psi_j^{-1}(V_{i\underline{j}})=U_{ij}$. Similarly, the notation
\[[\Lambda_{ijk}]: [f_{ik}, \hat{f}_{ik}]\cong [f_{jk},\hat{f}_{jk}]\circ [f_{ij},\hat{f}_{ij}]\]
means we find open subsets 
\[ V_{\underline{i}jk}\subset V_{\underline{i}j}\cap V_{\underline{i}k}, \;\; V_{i\underline{j}k} \subset V_{i\underline{j}}\cap V_{\underline{j}k}, \;\; V_{ij\underline{k}}\subset V_{i\underline{k}}\cap V_{j\underline{k}},\]
with $\psi_i^{-1}(V_{\underline{i}jk})=U_{ijk}$, $\psi_j^{-1}(V_{i\underline{j}k})=U_{ijk}$, $\psi_k^{-1}(V_{ij\underline{k}})=U_{ijk}$, so that the restrictions of $(f_{ij},\hat{f}_{ij})$, $(f_{jk},\hat{f}_{jk})$, $(f_{ik},\hat{f}_{ik})$, and the composition $(f_{jk},\hat{f}_{jk})\circ (f_{ij},\hat{f}_{ij})$ are well-defined. Then $\Lambda_{ijk}$ is defined over $V_{\underline{i}jk}$.

As Joyce puts it, the idea is we always take the maximal possible domain where all compositions are well-defined.

\subsection{Acknowledgment}  I am grateful to Professor Joyce for pointing out an important mistake in an earlier draft of the current paper, and for his encouragement to work out a complete categorical framework for Kuranishi manifolds.

\section{Definition of Kuranishi manifolds}~\label{sec:kur}

In this section, we define the notion of Kuranishi manifolds. 

\subsection{Kuranishi charts}

Let $X$ be a Hausdorff, second countable topological space. Let $d\in \Z$ be an integer.

\begin{defi}
A Kuranishi chart of $X$ of dimension $d$ is a quadruple $(V,\underline{\R^m},s,\psi)$ where $V\subset \R^n$ is an open subset, $s$ a smooth section of the trivial bundle $\underline{\R^m}$, and $$\psi: U\ra s^{-1}(0)$$ a homeomorphism from an open subset $U\subset X$ onto the zero locus of $s$. We also require that $d=n-m$. We call $U$ the footprint of the chart~\footnote{This terminology was suggested in the article~\cite{MW}.}.
\end{defi}

\begin{defi}
Let $(V_\alpha,\underline{\R^{m_\alpha}} ,s_\alpha, \psi_\alpha)$ and  $(V_\beta,\underline{\R^{m_\beta}},s_\beta,\psi_\beta)$ be two Kuranishi charts with the same footprint $U\subset X$. A Kuranishi morphism 
is given by a pair $(f,\hat{f})$ where $f: V_\alpha\ra V_\beta$ is a smooth map such that $f\circ \psi_\alpha=\psi_\beta$, and $\hat{f}: \underline{\R^{m_\alpha}} \ra f^*\underline{\R^{m_\beta}}$ a bundle map such that $\hat{f}\circ s_\alpha= f^*s_\beta$.
\end{defi}

To have the right notion of a ``coordinate change" between Kuranishi charts lies at the heart of the problem of obtaining a global theory of Kuranishi manifolds. There are various different definitions in the literature, which we refer to Joyce's recent work~\cite{Joyce} for an excellent survey. As mentioned in the introduction, we shall propose a new definition motivated by Lemma~\ref{intro-homotopy}.

Indeed, consider $(f,\hat{f}): (V_\alpha,\underline{\R^{m_\alpha}} ,s_\alpha, \psi_\alpha) \ra (V_\beta,\underline{\R^{m_\beta}},s_\beta,\psi_\beta)$  a morphism between Kuranishi charts, with the same footprint $U\subset X$. Differentiating the identity $\hat{f} s_\alpha= f^* s_\beta$, we get a morphism $(df,\hat{f})$ between $2$-term complexes
\begin{equation}~\label{diag:tangent}\begin{CD}
0@>>> \underline{\R^{n_\alpha}} @>ds_\alpha>> \underline{\R^{m_\alpha}} @>>> 0\\
@. @V df VV         @V\hat{f} VV @. \\
0@>>> \underline{\R^{n_\beta}} @> ds_\beta>> \underline{\R^{m_\beta}}@>>> 0
\end{CD}\end{equation}
We consider this as a diagram defined over $U$ by pulling back via $\psi_\alpha$.

\medskip
\begin{prop}~\label{prop:homology}
Let $(f_0,\hat{f}_0), (f_1,\hat{f}_1): (V_\alpha,\underline{\R^{m_\alpha}} ,s_\alpha, \psi_\alpha) \ra (V_\beta,\underline{\R^{m_\beta}},s_\beta,\psi_\beta)$ be two morphisms between Kuranishi charts, with the same footprint $U$. Assume that they satisfy the conditions $(a.)$,$(b.)$,$(c.)$ in Lemma~\ref{intro-homotopy}. Then the induced maps $(df_0,\hat{f}_0)$, $(df_1,\hat{f}_1)$ between the $2$-term complexes are homotopic.
\end{prop}

\Pf. We use the equations in Condition $(c.)$ of Lemma~\ref{intro-homotopy} to prove the existence of the following commutative diagram.
\[\begin{tikzcd}
0 \arrow{r} &\underline{\R^{n_\alpha}} \arrow{d}[swap]{df_1-df_0}\arrow{r}{ds_\alpha} &\underline{\R^{m_\alpha}} \arrow[dotted]{ld}[description]{\Lambda}\arrow{d}{\hat{f}_1-\hat{f}_0}\arrow{r} &0 \\
0 \arrow{r} &\underline{\R^{n_\beta}} \arrow{r}{ds_\beta} &\underline{\R^{m_\beta}} \arrow{r} &0
\end{tikzcd}\]
Indeed, integrating the identity $\frac{dF(t)}{dt}|_x = \Lambda(t)(s_\alpha(x))$ over $[0,1]$ implies that
\begin{equation}~\label{eq:difference} f_1(x)-f_0(x)=\big( \int_0^1\Lambda(t)dt\big) (s_\alpha(x)).\end{equation}
Denote by $\Lambda:= \int_0^1\Lambda(t)dt$, which a morphism $\underline{\R^{m_\alpha}}\ra \underline{\R^{n_\beta}}$. Differentiating the above equation in the $x$-direction gives
\[ df_1-df_0= d\Lambda (s_\alpha(x))+\Lambda ds_\alpha.\]
Now restricting this identity to $U$ via $\psi_\alpha$ implies that
\[ df_1-df_0=\Lambda ds_\alpha,\]
since $s_\alpha(x)=0$ if $x$ is in the image of $\psi_\alpha$.

Equation~\ref{eq:difference} above also implies that $f_0(x)=f_1(x)$ if $x\in \im \psi_\alpha$. In fact by the same proof, for all $t\in [0,1]$ and $x\in \im\psi_\alpha$, the image $F(t)(x)$ is independent of $t$. Thus, we may integrate the identity $\frac{d\hat{F}(t)}{dt}|_{(x,\xi)} = ds_\beta|_{F(t)(x)}\big(\Lambda(t)(x,\xi)\big)+\Xi(t)(s_\alpha(x)\wedge \xi)$ to obtain
\[ \hat{f}_1-\hat{f}_0=ds_\beta|_{f_0(x)} \Lambda.\]
Note that the second term $\Xi(t)(s_\alpha(x)\otimes\xi)=0$ after restriction to $U$ via $\psi_\alpha$. The proposition is proved.\ed

Again, consider two morphisms 
\[(f_0,\hat{f}_0), (f_1,\hat{f}_1): (V_\alpha,\underline{\R^{m_\alpha}} ,s_\alpha, \psi_\alpha) \ra (V_\beta,\underline{\R^{m_\beta}},s_\beta,\psi_\beta)\]
between Kuranishi charts of a topological space $U$. Before giving the definition of a homotopy between Kuranishi morphisms, we define the space of {\sl quotient bundle maps}~\footnote{The reason to use this quotient space is to have an associative composition of homotopies (see the next subsection).} as follows. A quotient bundle map is given by an equivalence of bundle maps  
$\Lambda\in \Hom_{V_\alpha}(\underline{\R^{m_\alpha}}, \underline{\R^{n_\beta}})$ such that 
\[ \Lambda\big(s_\alpha(x)\big)=f_1(x)-f_0(x),\]
where two such bundle maps $\Lambda_1$, $\Lambda_2$ are equivalent if their restrictions to $U$ are the same. Formally, the space of quotient bundle maps is given by
\[ \KK\Hom(\underline{\R^{m_\alpha}}, \underline{\R^{n_\beta}}):=\left\{ \Lambda\in \Hom_{V_\alpha}(\underline{\R^{m_\alpha}}, \underline{\R^{n_\beta}})|\;\;\Lambda\big(s_\alpha(x)\big)=f_1(x)-f_0(x).\right\}/ R,\]
where $R$ is the equivalence relation defined by
\[ R:=\left\{ \big(\Lambda_1,\Lambda_2\big) \mid \Lambda_1|_U=\Lambda_2|_U\right\}.\]
Motivated by Proposition~\ref{prop:homology}, we make the following 

\medskip
\begin{defi}~\label{def:homotopy}
Two morphisms 
\[(f_0,\hat{f}_0), (f_1,\hat{f}_1): (V_\alpha,\underline{\R^{m_\alpha}} ,s_\alpha, \psi_\alpha) \ra (V_\beta,\underline{\R^{m_\beta}},s_\beta,\psi_\beta)\]
between Kuranishi charts of a topological space $X$, with the same footprint $U\subset X$, are called homotopic if there exists a quotient
bundle map
\[ \Lambda\in \KK\Hom(\underline{\R^{m_\alpha}} , \underline{\R^{n_\beta}})\]
such that the following diagram defined over $U$  is commutative.
\[\begin{tikzcd}
0 \arrow{r} &\underline{\R^{n_\alpha}} \arrow{d}[swap]{df_1-df_0}\arrow{r}{ds_\alpha} &\underline{\R^{m_\alpha}} \arrow{ld}[swap]{\Lambda}\arrow{d}{\hat{f}_1-\hat{f}_0}\arrow{r} &0 \\
0 \arrow{r} &\underline{\R^{n_\beta}} \arrow{r}{ds_\beta} &\underline{\R^{m_\beta}} \arrow{r} &0
\end{tikzcd}\]
Denote this homotopy relation by $\Lambda:(f_0,\hat{f}_0)\cong (f_1,\hat{f}_1)$.  
\end{defi}

\begin{rem}
As shown in the proof of Proposition~\ref{prop:homology}, differentiating the identity $\Lambda\big(s_\alpha(x)\big)=f_1(x)-f_0(x)$ automatically implies the commutativity of the upper triangle. So the data of a homotopy boils down to the previous identity and the commutativity of the low triangle.
\end{rem}

\medskip
\begin{rem}~\label{rem:joyce}
The above definition should be compared with Joyce's Definition $2.3$ in~\cite{Joyce}, which leads to the main differences between Joyce's Kuranishi manifold theory and the current one. More precisely, the differences are
\begin{itemize}
\item[(a.)] We use the quotient space $\KK\Hom(\underline{\R^{m_\alpha}} , \underline{\R^{n_\beta}})$ as opposed to $\Hom(\underline{\R^{m_\alpha}} , \underline{\R^{n_\beta}})$ itself. This is technically useful to eliminate different choices of homotopies.
\item[(b.)] We work with local charts while Joyce's definition uses more general global charts. This flexibility leads to Joyce's big-$O$ notation, since one needs to choose a connection to differentiate.
\item[(c.)] Apart from the big-$O$ notation, the two definitions look rather similar. However, they are not equivalent, due to the fact that we only impose part $(2.)$, or the quation $\hat{f}_1-\hat{f}_0=ds_\beta \circ \Lambda$ over the footprint $U$, rather than on an actual chart. Indeed, consider the following example. 
\begin{align*}
& U=\left\{ 0\right\} \hookrightarrow \R \stackrel{s_\alpha}{\ra} \underline{\R^2}; \;\; s_\alpha(x)=(x^2,0)\\
& U=\left\{ 0\right\} \hookrightarrow \R \stackrel{s_\beta}{\ra} \underline{\R}; \;\; s_\beta(x)=x^3\\
& f_0(x)=x, \;\; \hat{f}_0(x)(v_1,v_2)=xv_1+xv_2\\
& f_1(x)=x, \;\; \hat{f}_1(x)(v_1,v_2)=xv_1+x^2v_2.
\end{align*}
One can check that the two Kuranishi morphisms $(f_0,\hat{f}_0), (f_1,\hat{f}_1)$ are homotopic in the sense of Definition~\ref{def:homotopy}, but are not equivalent in Joyce's definition. The reason being that the term $$\hat{f}_1(x)-\hat{f}_0(x)= (0,x-x^2)$$ vanishes if restricted to $U$, which implies that any morphism $\Lambda$ of the form $(0,a(x)): \underline{\R^2}\ra \underline{\R}$ would be a homotopy in our sense. However, in Joyce's definition, one needs to write the function $x-x^2$ into the form $ds_\beta  \Lambda+ O(x^2)$ in an actual open neighborhood of $0$ in $\R$.  But, since $s_\beta=x^3$, we have $ds_\beta=3x^2$, which implies that $ds_\beta  \Lambda\in O(x^2)$.
\end{itemize}
On the other hand, we still expect that a Kuranishi manifold in our sense (Definition~\ref{def:kuranishi}) should be equivalent to Joyce's definition. What could happen is that even though the equivalence relation used on morphisms between Kuranishi charts is different, it still produces the same notion of an isomorphism between Kuranishi charts.
\end{rem}

\begin{prop}
The relation $\cong$ is an equivalence relation.
\end{prop}

\Pf. Any morphism $(f,\hat{f})$ is homotopic to itself via $\Lambda=0$. If $(f_0,\hat{f}_0)\cong (f_1,\hat{f}_1)$ via $\Lambda$, then $(f_1,\hat{f}_1)\cong (f_0,\hat{f}_0)$ via $-\Lambda$. Finally if we have $(f_0,\hat{f}_0)\cong (f_1,\hat{f}_1)$ via $\Lambda_{01}$, and $(f_1,\hat{f}_1)\cong (f_2,\hat{f}_2)$ via $\Lambda_{12}$, then one can check that 
\[ (f_0,\hat{f}_0) \cong (f_2,\hat{f}_2)\]
via the homotopy $\Lambda_{01}+\Lambda_{12}$. \ed

\subsection{The $2$-category of Kuranishi charts}~\label{subsec:hori} Let $$(f_0,\hat{f}_0)\cong (f_1,\hat{f}_1):(V_\alpha,\underline{\R^{m_\alpha}} ,s_\alpha, \psi_\alpha) \ra (V_\beta,\underline{\R^{m_\beta}},s_\beta,\psi_\beta)$$ be homotopic morphisms between Kuranishi charts, via the homotopy data $\Lambda_{\alpha\beta}: \underline{\R^{m_\alpha}} \ra \underline{\R^{n_\beta}}$. Let $$(g_0,\hat{g}_0)\cong (g_1,\hat{g}_1): (V_\beta,\underline{\R^{m_\beta}},s_\beta,\psi_\beta)\ra(V_\gamma,\underline{\R^{m_\gamma}} ,s_\gamma, \psi_\gamma)$$ be another such pair of morphisms, with homotopy $\Lambda_{\beta\gamma}: \underline{\R^{m_\beta}} \ra \underline{\R^{n_\gamma}}$. 
We would like to construct a third homotopy
\[ \Lambda_{\alpha\gamma}: (g_0f_0,\hat{g}_0\hat{f}_0) \cong (g_1f_1,\hat{g}_1\hat{f}_1).\]
First recall the following form of Taylor's theorem.
\medskip
\begin{lem}~\label{lem:taylor}
Let $h: V\ra \R^l$ be a smooth map where $V\subset \R^k$ is an open and convex domain.
Then there exists a $l\times k$ matrix valued smooth map
\[ \delta h: V\times V \ra \Mat_{l\times k}\]
such that $ h(x)-h(y)=\delta h(x,y) (x-y)$, and $\delta h(x,x)= dh|_x$. Furthermore, if $h_1$ and $h_2$ are composible maps, then we have
\begin{equation}~\label{eq:taylor} 
\delta h_2\big( h_1(x), h_1(y)\big) \circ \delta h_1(x,y)(x-y) = \delta (h_2\circ h_1)(x,y)(x-y).
\end{equation}
\end{lem}

\Pf. By fundamental theorem of calculus, we have
\[ h(x)-h(y)=\int_0^1\frac{d}{dt} h\big(tx+(1-t)y\big)dt.\]
This is well-defined by the convexity of the domain $V$. Differentiating with chain rule yields
\[ h(x)-h(y)=\int_0^1 \sum_{j=1}^k \partial_j h \big(tx+(1-t)y\big) dt\cdot (x_j-y_j).\]
Thus the $l\times k$ matrix is given by
\[\delta h(x,y):= [ \int_0^1 \partial_1 h \big(tx+(1-t)y\big) dt \cdots \int_0^1  \partial_k h \big(tx+(1-t)y\big) dt].\]
Clearly, if $x=y$, this is just $dh|_x$. For the last identity, we have
\begin{align*}
\delta h_2\big( h_1(x), h_1(y)\big) \circ \delta h_1(x,y)(x-y) &= \delta h_2\big( h_1(x), h_1(y)\big) \big(h_1(x)-h_1(y)\big)\\
&= (h_2h_1)(x)-(h_2h_1)(y)\\
&=\delta(h_2h_1)(x,y)(x-y).
\end{align*}
\ed

Back to the construction of a homotopy $\Lambda_{\alpha\gamma}:(g_0f_0,\hat{g}_0\hat{f}_0) \cong (g_1f_1,\hat{g}_1\hat{f}_1)$, we may define the required homotopy $\Lambda_{\alpha\gamma}: \underline{\R^{m_\alpha}} \ra \underline{\R^{n_\gamma}}$ over $V_\alpha$, by the assignment
$$x\mapsto \Lambda_{\beta\gamma}|_{f_1(x)}\circ \hat{f}_1|_x+\delta g_0( f_1(x),f_0(x))\circ\Lambda_{\alpha\beta}|_x$$ 
at the point $x\in V_\alpha$. For simplicity, we denote this homotopy by
\[ \Lambda_{\alpha\gamma}:= \Lambda_{\beta\gamma} \hat{f}_1+\delta g_0 \Lambda_{\alpha\beta}.\]

\begin{lem}
The element $\Lambda_{\alpha\gamma}\in \KK\Hom(\underline{\R^{m_\alpha}}, \underline{\R^{n_\gamma}})$ defined above is a homotopy between $(g_0f_0,\hat{g}_0\hat{f}_0)$ and $(g_1f_1,\hat{g}_1\hat{f}_1)$.
\end{lem}

\Pf. we have
\begin{align*}
g_1f_1(x)-g_0f_0(x)&= g_1f_1(x)-g_0f_1(x)+g_0f_1(x)-g_0f_0(x)\\
&=\Lambda_{\beta\gamma}\big(s_\beta(f_1(x))\big)+\delta g_0 (f_1(x),f_0(x)) (f_1(x)-f_0(x))\\
&=\Lambda_{\beta\gamma}\hat{f}_1(s_\alpha(x))+\delta g_0( f_1(x),f_0(x))\Lambda_{\alpha\beta}(s_\alpha(x)).\\
&=\Lambda_{\alpha\gamma}(s_\alpha(x)).
\end{align*}
Furthermore, if $x=\psi_\alpha(u)$, then we have $f_1(x)=f_0(x)=\psi_\beta(u)$. Thus after restricting to $U$, by Lemma~\ref{lem:taylor}, we have
\[ \Lambda_{\alpha\gamma}|_U= \Lambda_{\beta\gamma} \circ \hat{f}_1+dg_0\circ \Lambda_{\alpha\beta}.\]
This implies that
\begin{align*}
\Lambda_{\alpha\gamma}|_U \circ ds_\alpha &= \Lambda_{\beta\gamma} \hat{f}_1 ds_\alpha+dg_0 \Lambda_{\alpha\beta} ds_\alpha\\
&= \Lambda_{\beta\gamma} ds_\beta df_1+ dg_0(df_1-df_0)\\
&=(dg_1-dg_0)df_1+dg_0(df_1-df_0)\\
&=dg_1df_1-dg_0df_0.
\end{align*}
Similarly, we also have
\begin{align*} 
ds_\gamma \circ \Lambda_{\alpha\gamma}|_U &= ds_\gamma\Lambda_{\beta\gamma} \hat{f}_1 +ds_\gamma dg_0 \Lambda_{\alpha\beta} \\
&=(\hat{g}_1-\hat{g}_0)\hat{f}_1+\hat{g}_0 ds_\beta \Lambda_{\alpha\beta}\\
&=(\hat{g}_1-\hat{g}_0)\hat{f}_1+\hat{g}_0 (\hat{f}_1-\hat{f}_0)\\
&= \hat{g}_1\hat{f}_1-\hat{g}_0\hat{f}_0.
\end{align*}
This proves the required property for $\Lambda_{\alpha\gamma}$. \ed

The homotopy $\Lambda_{\alpha\gamma}$ can viewed as the horizontal composition of $\Lambda_{\alpha\beta}$ and $\Lambda_{\beta\gamma}$. Denote this by $\Lambda_{\beta\gamma}*\Lambda_{\alpha\beta}:=\Lambda_{\alpha\gamma}$. 

\begin{prop}
The $*$ composition is associative.
\end{prop}

\Pf. Given three composible homotopies
\begin{align*}
\Lambda_{\alpha\beta}&:(f_0,\hat{f}_0)\cong (f_1,\hat{f}_1):(V_\alpha,\underline{\R^{m_\alpha}} ,s_\alpha, \psi_\alpha) \ra (V_\beta,\underline{\R^{m_\beta}},s_\beta,\psi_\beta)\\
\Lambda_{\beta\gamma}&:(g_0,\hat{g}_0)\cong (g_1,\hat{g}_1):(V_\beta,\underline{\R^{m_\beta}} ,s_\beta, \psi_\beta) \ra (V_\gamma,\underline{\R^{m_\gamma}},s_\gamma,\psi_\gamma)\\
\Lambda_{\gamma\eta}&:(h_0,\hat{h}_0)\cong (h_1,\hat{h}_1):(V_\gamma,\underline{\R^{m_\gamma}} ,s_\gamma, \psi_\gamma) \ra (V_\eta,\underline{\R^{m_\eta}},s_\eta,\psi_\eta),
\end{align*}
we want to verify that $(\Lambda_{\gamma\eta}*\Lambda_{\beta\gamma})*\Lambda_{\alpha\beta}=\Lambda_{\gamma\eta}*(\Lambda_{\beta\gamma}*\Lambda_{\alpha\beta})$, as an element of the space $\KK\Hom(\underline{\R^{m_\alpha}}, \underline{\R^{n_\eta}})$. Since both sides are homotopies between $(h_0g_0f_0,\hat{h}_0\hat{g}_0\hat{f}_0)$ and  $(h_1g_1f_1,\hat{h}_1\hat{g}_1\hat{f}_1)$, we have
\[ (\Lambda_{\gamma\eta}*\Lambda_{\beta\gamma})*\Lambda_{\alpha\beta} (s_\alpha(x))=\Lambda_{\gamma\eta}*(\Lambda_{\beta\gamma}*\Lambda_{\alpha\beta}) (s_\alpha(x)) = h_1g_1f_1(x)-h_0g_0f_0(x).\]
It remains to verify the two agrees after restriction onto $U$. We have
\begin{align*}
&(\Lambda_{\gamma\eta}*\Lambda_{\beta\gamma})*\Lambda_{\alpha\beta}|_{U}\\
=& (\Lambda_{\gamma\eta}\hat{g}_1+dh_0\Lambda_{\beta\gamma})*\Lambda_{\alpha\beta}|_{U}\\
=& \Lambda_{\gamma\eta}\hat{g}_1\hat{f}_1+dh_0\Lambda_{\beta\gamma}\hat{f}_1+dh_0dg_0\Lambda_{\alpha\beta}.
\end{align*}
\begin{align*}
&\Lambda_{\gamma\eta}*(\Lambda_{\beta\gamma}*\Lambda_{\alpha\beta})|_U\\
=&\Lambda_{\gamma\eta}*(\Lambda_{\beta\gamma}\hat{f}_1+dg_0\Lambda_{\alpha\beta})|_U\\
=&dh_0\Lambda_{\beta\gamma}\hat{f}_1+dh_0dg_0\Lambda_{\alpha\beta}+\Lambda_{\gamma\eta}\hat{g}_1\hat{f}_1.
\end{align*}
This completes the proof.\ed

\begin{rem}
To have the above proposition is the main reason we work with the space of quotient bundle maps $\KK\Hom(\underline{\R^{m_\alpha}}, \underline{\R^{n_\gamma}})$, rather than ordinary bundle maps
\end{rem}

We have already shown that the horizontal composition is associative. To confirm that the above structures define indeed a $2$-category, there is one more property to verify, the interchange law. 

\begin{prop}
We have the identity
\[ (\Lambda_{\beta\gamma}*\Lambda_{\alpha\beta})+(\Lambda_{\beta\gamma}'*\Lambda_{\alpha\beta}')=(\Lambda_{\beta\gamma}+\Lambda_{\beta\gamma}')*(\Lambda_{\alpha\beta}+\Lambda_{\alpha\beta}').\]
\end{prop}

\Pf. The two sides, being a homotopy between $(g_0f_0,\hat{g}_0\hat{f}_0)$ and $(g_2f_2,\hat{g}_2\hat{f}_2)$, assume the same value when evaluated at $s_\alpha(x)$. Thus it remains to prove the identity after restriction to $U$. Taking the difference of the two, and restricting to $U$ yields
\begin{align*}
&(\Lambda_{\beta\gamma}*\Lambda_{\alpha\beta})+(\Lambda_{\beta\gamma}'*\Lambda_{\alpha\beta}')-(\Lambda_{\beta\gamma}+\Lambda_{\beta\gamma}')*(\Lambda_{\alpha\beta}+\Lambda_{\alpha\beta}')\\
=&(\Lambda_{\beta\gamma}\hat{f}_1+dg_0\Lambda_{\alpha\beta})+(\Lambda_{\beta\gamma}'\hat{f}_2+dg_1\Lambda_{\alpha\beta}')-
(\Lambda_{\beta\gamma}+\Lambda_{\beta\gamma}')\hat{f}_2-dg_0(\Lambda_{\alpha\beta}+\Lambda_{\alpha\beta}')\\
=&(dg_1-dg_0)\Lambda_{\alpha\beta}'-\Lambda_{\beta\gamma}(\hat{f}_2-\hat{f}_1)\\
=&\Lambda_{\beta\gamma}ds_{\beta}\Lambda_{\alpha\beta}'-\Lambda_{\beta\gamma}ds_{\beta}\Lambda_{\alpha\beta}'\\
=&0.
\end{align*}
This completes the proof.\ed

To summerize, let $U$ be a topological space. We can define a $2$-category of Kuranishi charts of $U$:
\begin{itemize}
\item objects are Kuranishi charts $(V,\underline{\R^m},s,\psi)$ with footprint $U$,
\item $1$-morphisms are $(f,\hat{f}):(V_\alpha,\underline{\R^{m_\alpha}} ,s_\alpha, \psi_\alpha) \ra (V_\beta,\underline{\R^{m_\beta}},s_\beta,\psi_\beta)$, composition of $1$-morphisms is given by $(g,\hat{g})(f,\hat{f})=(gf,\hat{g}\hat{f})$.
\item $2$-morphisms between two $1$-morphisms $(f_0,\hat{f}_0)$ and $(f_1,\hat{f}_1)$ is a homotopy $\Lambda\in \KK\Hom(\underline{\R^{m_\alpha}},\underline{\R^{n_\beta}})$ over $V_\alpha$,
\item Vertical composition between $2$-morphisms
\[ \Lambda: (f_0,\hat{f}_0)\ra (f_1,\hat{f}_1), \;\;\; \Lambda': (f_1,\hat{f}_1)\ra (f_2,\hat{f}_2)\]
is defined by addition $\Lambda+\Lambda'$ in $\KK\Hom(\underline{\R^{m_\alpha}},\underline{\R^{n_\beta}})$,
\item Horizontal composition between $2$-morphisms
\begin{align*}
\Lambda_{\alpha\beta}&:(f_0,\hat{f}_0)\cong (f_1,\hat{f}_1):(V_\alpha,\underline{\R^{m_\alpha}} ,s_\alpha, \psi_\alpha) \ra (V_\beta,\underline{\R^{m_\beta}},s_\beta,\psi_\beta)\\
\Lambda_{\beta\gamma}&:(g_0,\hat{g}_0)\cong (g_1,\hat{g}_1):(V_\beta,\underline{\R^{m_\beta}} ,s_\beta, \psi_\beta) \ra (V_\gamma,\underline{\R^{m_\gamma}},s_\gamma,\psi_\gamma)
\end{align*}
is defined by $\Lambda_{\beta\gamma}*\Lambda_{\alpha\beta}$.
\end{itemize}

\medskip
\subsection{Germs}~\label{subsec:germ}  Let $(V, \underline{\R^m}, s, \psi)$ be a Kuranishi chart of a topological space with footprint $U$.  Let $U'\subset U$ be an open subset. Since $U$ is homeomorphic to $s^{-1}(0)$ endowed with the subspace topology of $V$, there exists an open subset $V'\subset V$ such that $\psi^{-1}(V')=U'$. It is clear that the quadruple $(V',\underline{\R^m}, s|_{V'}, \psi|_{U'})$ is a Kuranishi chart with footprint $U'$. We shall use the abbreviation $(V',\underline{\R^m}, s, \psi)$ for this quadruple. We call the natural morphism
\[ (V',\underline{\R^m}, s, \psi) \ra (V,\underline{\R^m}, s, \psi)\]
an {\em open inclusion}.

In the case when $U'=U$, we may still have $V'\neq V$. To solve this non-uniqueness, it is convenient to work with germs of charts. In literature, working with germs is tricky as pointed by Joyce~\cite{Joyce}. The notion of germs used here is different in the sense that only open inclusions are allowed. Namely, we say two Kuranishi charts $(V_1, \underline{\R^{m_1}}, s_1, \psi_1)$ and $(V_2,\underline{\R^{m_2}},s_2,\psi_2)$ with the same footprint $U$ equivalent if there exists a common open subset $V$ in $V_1$ and $V_2$ such that the two restricted charts $(V,\underline{\R^{m_1}}, s_1, \psi_1)$ and $(V,\underline{\R^{m_2}},s_2,\psi_2)$ are the same Kuranishi chart with footprint $U$. It is straightforward to check this gives an equivalence relation. A germ of Kuranishi chart is given by an equivalence class of Kuranishi charts. We denote by $[V,\underline{\R^m},s,\psi]$ the equivalence class of $(V,\underline{\R^m},s,\psi)$. If $U'\subset U$ is an open subset of the footprint $U$ and $[V,\underline{\R^m},s,\psi]$ a germ of chart on $U$, there exists a unique germ of Kuranishi chart $[V',\underline{\R^m},s,\psi]$ given by the class of the restriction of a representative $(V,\underline{\R^m},s,\psi)$ onto an open subset $V'$ such that $\psi^{-1}(V')=U'$.

Let $(f_1,\hat{f}_1), (f_2,\hat{f}_2): (V_\alpha,\underline{\R^{m_\alpha}} ,s_\alpha, \psi_\alpha) \ra (V_\beta,\underline{\R^{m_\beta}},s_\beta,\psi_\beta)$ be two morphisms between Kuranishi charts of $U$. They are called equivalent if there exist open subsets $V_\alpha'\subset V_\alpha$ and $V_\beta'\subset V_\beta$ such that $\psi_\alpha^{-1}(V_\alpha')=\psi_\beta^{-1}(V_\beta')=U$, and that the following two restrictions
\begin{align*}
(f_1,\hat{f}_1)&: (V_\alpha',\underline{\R^{m_\alpha}} ,s_\alpha, \psi_\alpha) \ra (V_\beta',\underline{\R^{m_\beta}},s_\beta,\psi_\beta),\\
(f_2,\hat{f}_2)&: (V_\alpha',\underline{\R^{m_\alpha}} ,s_\alpha, \psi_\alpha) \ra (V_\beta',\underline{\R^{m_\beta}},s_\beta,\psi_\beta)
\end{align*}
are well-defined, and are the same Kuranishi morphism. A germ of Kuranishi morphism is an equivalence class of Kuranishi morphisms. We define Kuranishi morphisms
\[ \Hom\big([V_\alpha,\underline{\R^{m_\alpha}} ,s_\alpha, \psi_\alpha], [V_\beta,\underline{\R^{m_\beta}},s_\beta,\psi_\beta]\big)\]
between germs of Kuranishi charts of $U$, to be germs of Kuranishi morphisms between any representing charts $(V_\alpha,\underline{\R^{m_\alpha}} ,s_\alpha, \psi_\alpha)$ and $(V_\beta,\underline{\R^{m_\beta}},s_\beta,\psi_\beta)$.

Similarly, we define germs of homotopies between germs of Kuranishi morphisms. Let 
\[ [f_0,\hat{f}_0], [f_1,\hat{f}_1]: [V_\alpha,\underline{\R^{m_\alpha}} ,s_\alpha, \psi_\alpha]\ra  [V_\beta,\underline{\R^{m_\beta}},s_\beta,\psi_\beta]\]
be two germs of Kuranishi morphisms. Two homotopies $\Lambda_1,\Lambda_2\in \KK\Hom(\underline{\R^{m_\alpha}},\underline{\R^{m_\beta}})$ defined over $V_\alpha$ are called equivalent if there exists an open subset $V\subset V_\alpha$ such that the restrictions are equal $\Lambda_1|_V=\Lambda_2|_V$. We define a germ of homotopy by an equivalence class of homotopies.

\subsection{Definition of Kuranishi manifolds} In order to apply Lemma~\ref{lem:taylor}, we shall use  \emph{hypercoverings} instead of \v Cech coverings in the following definition of Kuranishi atlases, since convex open subsets form a basis of the Euclidean topology. Roughly, a hypercovering allows refinements on every level of multiple intersections of open sets appearing in previous levels. This causes quite heavy notations to work with. For this reason, we choose to stick with the \v Cech notation while keeping in mind that we are dealing with possibly refinements of intersections. For example, if $X$ is covered by open subsets $\left\{ U_i\right\}$, we shall use the notation $U_{ij}$ to mean an open subset appearing in the refinement of intersection $U_i\cap U_j$. In other words, we have
\[ U_i\cap U_j = \cup\; U_{ij}\]
where we suppressed the indexing set of this refinement. Similarly, for $U_{ij}$, $U_{jk}$, and $U_{ik}$ appearing in the refinements of double intersections, we use $U_{ijk}$ to represent an open set appearing in the refinement of the intersection $U_{ij}\cap U_{jk}\cap U_{ik}$. Formal, we have
\[ U_{ij}\cap U_{jk}\cap U_{ik}= \cup\; U_{ijk}\]
where again we omitted the indexing set for clarity. This notation is appropriate since all proofs in the \v Cech setting carry through word by word in the hypercovering setting. For details of hypercoverings used here, see~\cite{Tu}.

\medskip

\begin{defi}~\label{def:kuranishi}
Let $X$ be a Hausdorff, second countable topological space, and $d\in \Z$.  A Kuranishi atlas on $X$ of virtual dimension $d$ is given by a collection of germs of Kuranishi charts $\Big\{[V_{i},\underline{\R^{m_i}},s_i,\psi_i]\Big\}_{i\in I}$, together with transition morphisms $\Big\{[f_{ij},\hat{f}_{ij}]\Big\}$ and homotopy data $\Big\{[\Lambda_{ijk}]\Big\}$, where
\begin{itemize}
\item[(1.)] $[V_{i},\underline{\R^{m_i}},s_i,\psi_i]$ is a Kuranishi chart with footprint $U_i\subset X$ with $V_i\subset\R^{n_i}$ an open domain, and $n_i-m_i=d$, for all $i\in I$. We require that $\cup_{i\in I} U_i=X$.
\medskip
\item[(2.)] For each $U_{ij}\neq \emptyset$ appearing in the refinement of the open set $U_i\cap U_j$, we have transition maps
\[ f_{ij}: [V_i,\underline{\R^{m_i}}, s_i, \psi_i] \ra [V_j,\underline{\R^{m_j}}, s_j, \psi_j].\]
We require that $[f_{ii}, \hat{f}_{ii}]=[\id]$ for all $i\in I$.
\medskip
\item[(3.)] For each $U_{ijk}\neq \emptyset$ appearing in the refinement of $U_{ij}\cap U_{jk}\cap U_{ik}$, there exists a homotopy in the sense of Definition~\ref{def:homotopy}:
\[ [\Lambda_{ijk}]: [f_{ik}, \hat{f}_{ik}]\cong [f_{jk},\hat{f}_{jk}]\circ [f_{ij},\hat{f}_{ij}].\]
We require that $[\Lambda_{iij}]=0$ and $[\Lambda_{ijj}]=0$.
\medskip
\item[(4.)] For each open subset $U_{ijkl}\neq\emptyset$ appearing in the refinement of $U_{ijk}\cap U_{ijl} \cap U_{ikl}\cap U_{jkl}$, we require the equation
\begin{equation}\label{eq:cocycle2}
\Lambda_{ikl}-\Lambda_{jkl}\circ \hat{f}_{ij} -\Lambda_{ijl}+df_{kl}\circ\Lambda_{ijk}=0
\end{equation}
to hold over the topological space $U_{ijkl}$.
\end{itemize}
To avoid rather long notations, even though $\Big\{[f_{ij},\hat{f}_{ij}]\Big\}_{i,j\in I}$ and $\Big\{[\Lambda_{ijk}]\Big\}_{i,j,k\in I}$ are both part of the data, we shall simply use $\Big\{[V_{i},\underline{\R^{m_i}},s_i,\psi_i]\Big\}_{i\in I}$ to denote a Kuranishi atlas. The topological space $X$ together with a Kuranishi atlas is called a Kuranishi manifold. We shall use calligraphic letters such as $\mathfrak{X}, \mathfrak{Y}$ to denote Kuranishi manifolds.
\end{defi}

\begin{rem}
Equation~\ref{eq:cocycle2} comes from the following diagram:
\[\begin{xy} 
(0,0)*{[V_{i},\underline{\R^{m_i}},s_i,\psi_i]}="A"; 
(40,-30)*{[V_{j},\underline{\R^{m_j}},s_j,\psi_j]}="B";
(80,0)*{[V_{l},\underline{\R^{m_l}},s_l,\psi_l]}="C";
(10,50)*{[V_{k},\underline{\R^{m_k}},s_k,\psi_k]}="D";
{\ar@{.>}_{[f_{il},\hat{f}_{il}]} "A"; "C"};
{\ar@{->}_{[f_{ij},\hat{f}_{ij}]} "A"; "B"};
{\ar@{->}_{[f_{jl},\hat{f}_{jl}]} "B"; "C"};
{\ar@{->}^{[f_{ik},\hat{f}_{ik}]} "A"; "D"};
{\ar@{->}_{[f_{jk},\hat{f}_{jk}]} "B"; "D"};
{\ar@{->}^{[f_{kl},\hat{f}_{kl}]} "D"; "C"};
\end{xy}\]
It asserts certain cancellation of the four homotopies on the four facets of the above tetrahedron.
\end{rem}

\medskip
\noindent The following are some examples of Kuranishi manifolds.

\begin{itemize}
\item[(A.)] If $\mathfrak{X}$ is a Kuranishi manifold with all the obstruction bundles $\underline{\R^{m_i}}=0$, which enforces $\Lambda_{ijk}=0$, then Definition~\ref{def:kuranishi} reduces to the definition of ordinary smooth manifolds.
\item[(B.)] By the discussion in the introduction, there exists a truncation 
\[ \;\;\;\;\;\;\;\;\;\;\mbox{$[0,1]$-type homotopy $L_\infty$ spaces} \Rightarrow \mbox{Kuranishi manifolds}.\]
In particular, as shown in~\cite{Tu}, moduli spaces of simple coherent sheaves on a Calabi-Yau $3$-fold are examples of Kuranishi manifolds.
\item[(C.)] Consider the fiber product $Z=X\times_M Y$ in the following diagram where $X$, $Y$, and $M$ are all smooth manifolds.
\[\begin{CD}
Z @>>> Y\\
@VVV   @V g VV\\
X @> h >> M.\end{CD}\]
Assuming $g$ and $h$ are both smooth maps, we shall prove in the next section that there exists a natural Kuranishi structure on $Z$. In particular, taking $Y$ to be a point, this shows that the preimage $h^{-1}(c)$ of any point $c\in M$ under any smooth map carries a natural Kuranishi structure. In fact, both $X$ and $Y$ can be taken to be Kuranishi manifolds $\mathfrak{X}$, $\mathfrak{Y}$, while still keeping $M$ an ordinary manifold.
\end{itemize}

\section{The $2$-category of Kuranishi manifolds}~\label{sec:cat}

In this section, we prove that Kuranishi manifolds form a $2$-category whose $2$-morphisms are all invertible.

\subsection{The $2$-category ${\sf pre-}\mathfrak{Kur}$} To define morphisms between Kuranishi manifolds, we first need a generalization of Definition~\ref{def:homotopy} for Kuranishi charts on different topological spaces. Let $U_\alpha$ and $U_\beta$ be two topological spaces. Let $(V_\alpha,\underline{\R^{m_\alpha}} ,s_\alpha, \psi_\alpha)$ and $(V_\beta,\underline{\R^{m_\beta}},s_\beta,\psi_\beta)$ be Kuranishi charts with footprints $U_\alpha$ and $U_\beta$ respectively. A morphism between them is a pair $(f,\hat{f})$ where
$f: V_\alpha\ra V_\beta$ a smooth map such that $f(\psi_\alpha(U_\alpha))\subset \psi_\beta(U_\beta)$, and $\hat{f}: \underline{\R^{m_\alpha}}\ra f^*\underline{\R^{m_\beta}}$ is a bundle map such that
\[ f^*s_\beta= \hat{f}\circ s_\alpha.\]
Since the map $f$ satisfies $f(\psi_\alpha(U_\alpha))\subset \psi_\beta(U_\beta)$, it induces a well defined map
\[ f|_{U_\alpha}: U_\alpha \ra U_\beta.\]
We abused the notation $f|_{U_\alpha}$ as $U_\alpha$ is not really a subspace of $V_\alpha$, but this should not cause any confusion. Two morphisms 
\[(f_0,\hat{f}_0), (f_1,\hat{f}_1): (V_\alpha,\underline{\R^{m_\alpha}} ,s_\alpha, \psi_\alpha) \ra (V_\beta,\underline{\R^{m_\beta}},s_\beta,\psi_\beta)\]
between Kuranishi charts of topological space $U_\alpha$ and $U_\beta$ are called homotopic if there exists a quotient bundle map
\[ \Lambda\in \KK\Hom( \underline{\R^{m_\alpha}},\underline{\R^{n_\beta}})\]
over $V_\alpha$, satisfying conditions $(1.)$, $(2.)$ in Definition~\ref{def:kuranishi}. In particular, condition $(1.)$ that 
\[ f_1(x)-f_0(x)=\Lambda(s_\alpha(x)),\]
and the fact that $s_\alpha(x)=0$ for $x\in \im \psi_\alpha$, imply that
\begin{equation}~\label{eq:space}
f_0|_{U_\alpha}=f_1|_{U_\alpha}.
\end{equation}
Again we denote this relation by $(f_0,\hat{f}_0)\cong (f_1,\hat{f}_1)$, generalizing Definition~\ref{def:homotopy}. The corresponding definition of germs of morphisms and homotopies is done in the same way as in Subsection~\ref{subsec:germ}.
\medskip
\begin{defi}~\label{def:homo}
Let $\mathfrak{X}$, $\mathfrak{Y}$ be two Kuranishi manifolds, with Kuranishi atlases given by
\begin{align*}\AA &=\big\{ [V_i,\underline{\R^{m_i}}, s_i, \psi_i], [f_{ij},\hat{f}_{ij}], [\Lambda_{ijk}]\big\}_{i,j,k\in I}, \\
\BB &= \big\{ [V_p,\underline{\R^{m_p}}, s_p, \psi_p], [f_{pq},\hat{f}_{pq}], [\Lambda_{pqr}] \big\}_{p,q,r\in P}\end{align*}
A {\sl strict Kuranishi morphism} $h: \mathfrak{X}\ra\mathfrak{Y}$ is given by the following data:
\begin{itemize}
\item[(1.)] A map of indices $\tau: I \ra P$.
\item[(2.)] For each $i\in I$, a morphism
\[ [h_{i}, \hat{h}_{i}]: [V_{i}, \underline{\R^{m_{i}}}, s_{i}, \psi_{i}] \ra [V_{\tau(i)}, \underline{\R^{m_{\tau(i)}}}, s_{\tau(i)}, \psi_{\tau(i)}]\]
between germs of Kuranishi charts.
\item[(3.)] For each pair of indices $(i,j)\in I\times I$ such that $U_{ij}\neq \emptyset$, a germ of quotient bundle map $\Delta_{ij}\in\KK\Hom(\underline{\R^{m_i}}, \underline{\R^{n_{\tau(j)}}})$, giving a homotopy
\[ [\Delta_{ij}]: [h_{j},\hat{h}_{j}]\circ [f_{ij},\hat{f}_{ij}] \cong [f_{\tau(i)\tau(j)},\hat{f}_{\tau(i)\tau(j)}]\circ [h_{i},\hat{h}_{i}]. \]
We require that $\Delta_{ii}=[0]$. Diagrammatically, the following commutes up to homotopy $\Delta_{ij}$.
\[\begin{CD}
[V_{i}, \underline{\R^{m_{i}}}, s_{i}, \psi_{i}]@>[f_{ij},\hat{f}_{ij}]>> [V_{j}, \underline{\R^{m_{j}}}, s_{j}, \psi_{j}]\\
@V [h_{i},\hat{h}_{i}] VV @VV [h_{j}, \hat{h}_{j}] V\\
[V_{\tau(i)}, \underline{\R^{m_{\tau(i)}}}, s_{\tau(i)}, \psi_{\tau(i)}] @> [f_{\tau(i)\tau(j)},\hat{f}_{\tau(i)\tau(j)}]>> [W_{\tau(j)}, \underline{\R^{m_{\tau(j)}}}, s_{\tau(j)}, \psi_{\tau(j)}]
\end{CD}\]
\item[(4.)] For a triple of indices $(i,j,k)\in I\times I\times I$ such that $U_{ijk}\neq \emptyset$, we require the equation
\begin{equation}~\label{eq:2-hom} 
\begin{split}
\Delta_{ik}-dh_{k}\circ\Lambda_{ijk}+&\Lambda_{\tau(i)\tau(j)\tau(k)}\circ \hat{h}_{i}-df_{\tau(j)\tau(k)}\circ\Delta_{ij}-\Delta_{jk}\circ \hat{f}_{ij}=0
\end{split}
\end{equation}
to hold over $U_{ijk}$.
\end{itemize}
Note that, as in the definition of Kuranishi manifolds, we consider the homotopies $[\Delta_{ij}]$ as part of the data. 
\end{defi}

\begin{rem}
We explain the appearance of Equation~\ref{eq:2-hom}. Consider the following diagram
\[\begin{xy} 
(0,0)*{V_{\tau(i)}}="A"; (70,0)*{V_{\tau(k)}}="C"; (35,20)*{V_{\tau(j)}}="B";
(0,50)*{V_i}="D"; (70,50)*{V_{k}}="F"; (35,70)*{V_{j}}="E";
{\ar@{->}^{[f_{\tau(i)\tau(k)},\hat{f}_{\tau(i)\tau(k)}]} "A"; "C"};
{\ar@{->}|-{[h_{i},\hat{h}_{i}]} "D"; "A"};
{\ar@{.>}|-{[f_{\tau(i)\tau(j)},\hat{f}_{\tau(i)\tau(j)}]} "A"; "B"};
{\ar@{.>}|-{[f_{\tau(j)\tau(k)},\hat{f}_{\tau(j)\tau(k)}]} "B"; "C"};
{\ar@{->}^{[f_{ij},\hat{f}_{ij}]} "D"; "E"};
{\ar@{->}^{[f_{jk},\hat{f}_{jk}]} "E"; "F"};
{\ar@{->}^{[f_{ik},\hat{f}_{ik}]} "D"; "F"};
{\ar@{.>}|-{[h_{j},\hat{h}_{j}]} "E"; "B"};
{\ar@{->}|-{[h_{k},\hat{h}_{k}]} "F"; "C"};
\end{xy}\]
\medskip
The two horizontal triangles are commutative up to the homotopies given by $\Lambda_{ijk}$ and $\Lambda_{\tau(i)\tau(j)\tau(k)}$ prescribed in the Kuranishi structures of $\mathfrak{X}$ and $\mathfrak{Y}$, while the three vertical squares are commutative up to the homotopies $\Delta_{ij}$, $\Delta_{jk}$, and $\Delta_{ik}$. Equation~\ref{eq:2-hom} asserts certain cancellations of these five homotopies.
\end{rem}

\begin{lem}
Let $h=\Big\{\tau, [h_{i},\hat{h}_{i}], [\Delta_{ij}]\Big\}: (X,\AA)\ra (Y,\BB)$ be a strict Kuranishi morphism as in Definition~\ref{def:homo}. Then it induces a continuous map, denoted by
\[ \underline{h}: X \ra Y\]
on the underlying topological spaces.
\end{lem}

\Pf. Since $[h_i,\hat{h}_i]$'s are morphisms of germs of Kuranishi charts, it induces a continuous map
\[ \underline{h}_i: U_i \ra U_{\tau(i)}.\]
The existence of the homotopy $\Delta_{ij}$ implies that the locally defined maps $\underline{h}_i$'s are compatible on intersections by Equation~\ref{eq:space}, which yields a global map $\underline{h}$. The continuity follows from the identity
\[ \underline{h}^{-1}(U)= \bigcup_{p\in P} \underline{h}^{-1}(U\cap U_p),\]
and local continuity.\ed

Next, we define compositions of strict Kuranishi morphisms. Let
\[ h=\Big\{\tau, [h_{i},\hat{h}_{i}], [\Delta_{ij}]\Big\}: (X,\AA)\ra (Y,\BB), \;\; g=\Big\{\sigma,[g_p,\hat{g}_p],[\Delta_{pq}]\Big\}: (Y,\BB)\ra (Z,\CC)\]
be two strict Kuranishi morphisms. We define the composition morphism 
\[  gh:=\Big\{\sigma\tau, [g_{\tau(i)}h_i, \hat{g}_{\tau(i)}\hat{h}_i], [\Delta_{\tau(i)\tau(j)}*[h_i,\hat{h}_i]+[g_{\tau(j)},\hat{g}_{\tau(j)}]*\Delta_{ij}]\Big\}: (X,\AA)\ra (Z,\CC).\]
Here $*$ is the horizontal composition defined in Subsection~\ref{subsec:hori}.

\medskip
\begin{prop}~\label{prop:pre-comp}
The composition $gh$ defined above forms a strict Kuranishi morphism from $\mathfrak{X}$ to $\mathfrak{Z}$.  Furthermore, this composition product of strict Kuranishi morphisms is associative.
\end{prop}

\Pf. The only non-trivial part in verifying that $gh$ is a strict Kuranishi morphism, is to show that Equation~\ref{eq:2-hom} holds. Indeed since both $h$ and $g$ are strict Kuranishi morphisms, we have
\begin{align*}
\Delta_{ik}-dh_{k}\circ\Lambda_{ijk}+\Lambda_{\tau(i)\tau(j)\tau(k)}\circ \hat{h}_{i}-df_{\tau(j)\tau(k)}\circ\Delta_{ij}-\Delta_{jk}\circ \hat{f}_{ij} &=0\\
\Delta_{\tau(i)\tau(k)}-dg_{\tau(k)}\circ\Lambda_{\tau(i)\tau(j)\tau(k)}+\Lambda_{\sigma\tau(i)\sigma\tau(j)\sigma\tau(k)}\circ \hat{g}_{\tau(i)}&\\
-df_{\sigma\tau(j)\sigma\tau(k)}\circ\Delta_{\tau(i)\tau(j)}-\Delta_{\tau(j)\tau(k)}\circ \hat{f}_{\tau(i)\tau(j)}&=0
\end{align*}
Post-composing the first equation with $dg_{\tau(k)}$ and pre-composing the second equation with $\hat{h}_i$, we get
\begin{align*}
&dg_{\tau(k)}\Big(\Delta_{ik}-dh_{k}\circ\Lambda_{ijk}+\Lambda_{\tau(i)\tau(j)\tau(k)}\circ \hat{h}_{i}-df_{\tau(j)\tau(k)}\circ\Delta_{ij}-\Delta_{jk}\circ \hat{f}_{ij}\Big)+\\
&\Big(\Delta_{\tau(i)\tau(k)}-dg_{\tau(k)}\circ\Lambda_{\tau(i)\tau(j)\tau(k)}+\Lambda_{\sigma\tau(i)\sigma\tau(j)\sigma\tau(k)}\circ \hat{g}_{\tau(i)}\\
&-df_{\sigma\tau(j)\sigma\tau(k)}\circ\Delta_{\tau(i)\tau(j)}-\Delta_{\tau(j)\tau(k)}\circ \hat{f}_{\tau(i)\tau(j)}\Big)\hat{h}_i=0
\end{align*}
Reorganizing the terms and using the equations
\begin{align*}
dg_{\tau(k)}df_{\tau(j)\tau(k)}&=df_{\sigma\tau(j)\sigma\tau(k)}dg_{\tau(j)}-\Delta_{\tau(j)\tau(k)}ds_{\tau(j)}\\
\hat{f}_{\tau(i)\tau(j)}\hat{h}_i&=\hat{h}_j\hat{f}_{ij}+ds_{\tau(j)}\Delta_{ij},
\end{align*}
we obtain
\begin{align*}
&\Big(\Delta_{\tau(i)\tau(k)}\hat{h}_i+dg_{\tau(k)}\Delta_{ik}\Big)-\Big(dg_{\tau(k)}dh_k\Lambda_{ijk}\Big)+\Big(\Lambda_{\sigma\tau(i)\sigma\tau(j)\sigma\tau(k)}\circ \hat{g}_{\tau(i)}\hat{h}_i\Big)\\
&-df_{\sigma\tau(j)\sigma\tau(k)}\Big(\Delta_{\tau(i)\tau(j)}\hat{h}_i+dg_{\tau(j)}\Delta_{ij}\Big)-\Big(\Delta_{\tau(j)\tau(k)}\hat{h}_j+dg_{\tau(k)}\Delta_{jk}\Big)\hat{f}_{ij}=0,
\end{align*}
which is precisely Equation~\ref{eq:2-hom} for the morphism $gh$. 

The fact that this composition is associative follows immediately from the associativity of the $*$-composition.\ed

We denote by ${\sf pre-}\mathfrak{Kur}$ the category of Kuranishi manifolds with strict Kuranishi morphisms. This can be enhanced to a $2$-category with the following definition of $2$-morphisms. 

\medskip
\begin{defi}~\label{def:equihom}
Let $$h^1=\Big\{\tau_1, [h^1_{i},\hat{h}^1_{i}], [\Delta^1_{ij}]\Big\}, h^2=\Big\{\tau_2, [h^2_{i},\hat{h}^2_{i}], [\Delta^2_{ij}]\Big\}: (X,\AA) \ra (Y,\BB)$$ be two strict morphisms between Kuranishi manifolds. Assume that the
underlying continuous maps from $X$ to $Y$ are equal: $\underline{h^1}=\underline{h^2}$. A $2$-morphism from $h^1$ to $h^2$ is given by the following data:
\begin{itemize}
\item[(a.)] For each $i\in I$, a germ of quotient bundle map $[\Upsilon_i]\in \KK\Hom(\underline{\R^{m_i}}, \underline{\R^{m_{\tau_2(i)}}})$, giving a homotopy $[\Upsilon_i]:[h^2_i, \hat{h}^2_i]\cong [f_{\tau_1(i)\tau_2(i)}, \hat{f}_{\tau_1(i)\tau_2(i)}] \circ [h^1_i,\hat{h}^1_i]$. Diagrammatically, the following is commutative, up to homotopy $[\Upsilon_i]$.
\[\begin{xy} 
(0,0)*{[V_i,\underline{\R^{m_i}},s_i,\psi_i]}="A"; (80,0)*{[V_{\tau_2(i)},\underline{\R^{m_{\tau_2(i)}}},s_{\tau_2(i)},\psi_{\tau_2(i)}]}="C"; (40,20)*{[V_{\tau_1(i)},\underline{\R^{m_{\tau_1(i)}}},s_{\tau_1(i)},\psi_{\tau_1(i)}]}="B";
{\ar@{->}^{[h^2_i,\hat{h}^2_i]} "A"; "C"};
{\ar@{->}|-{[h^1_i,\hat{h}^1_i]} "A"; "B"};
{\ar@{->}|-{[f_{\tau_1(i)\tau_2(i)},\hat{f}_{\tau_1(i)\tau_2(i)}]} "B"; "C"};
\end{xy}\]
\item[(b.)] For a pair of indices $(i,j)\in I\times I$, we require the equation
\begin{equation}~\label{eq:2-equiv} 
\begin{split}
\Delta^2_{ij}+&df_{\tau_2(i)\tau_2(j)} \circ \Upsilon_i - \Upsilon_j\circ \hat{f}_{ij} - df_{\tau_1(j)\tau_2(j)} \circ \Delta^1_{ij} \\ &-[\Lambda_{\tau_1(i)\tau_2(i)\tau_2(j)}-\Lambda_{\tau_1(i)\tau_1(j)\tau_2(j)}]\circ \hat{h}^1_i=0,
\end{split}
\end{equation}
to hold over the intersection $U_{ij}$.
\end{itemize}
\end{defi}

\begin{rem}
Equation~\ref{eq:2-equiv} is illustrated in the following diagram.
\[\begin{xy} 
(0,0)*{V_i}="A"; (70,0)*{V_{\tau_2(i)}}="C"; (35,20)*{V_{\tau_1(i)}}="B";
(0,60)*{V_j}="D"; (70,60)*{V_{\tau_2(j)}}="F"; (35,80)*{V_{\tau_1(j)}}="E";
{\ar@{->}^{[h^2_i,\hat{h}^2_i]} "A"; "C"};
{\ar@{->}|-{[f_{ij},\hat{f}_{ij}]} "A"; "D"};
{\ar@{.>}|-{[h^1_i,\hat{h}^1_i]} "A"; "B"};
{\ar@{.>}|-{[f_{\tau_1(i)\tau_2(i)},\hat{f}_{\tau_1(i)\tau_2(i)}]} "B"; "C"};
{\ar@{->}^{[h^1_j,\hat{h}^1_j]} "D"; "E"};
{\ar@{->}^{[f_{\tau_1(j)\tau_2(j)},\hat{f}_{\tau_1(j)\tau_2(j)}]} "E"; "F"};
{\ar@{->}^{[h^2_j,\hat{h}^2_j]} "D"; "F"};
{\ar@{.>}|-{[f_{\tau_1(i)\tau_1(j)},\hat{f}_{\tau_1(i)\tau_1(j)}]} "B"; "E"};
{\ar@{->}|-{[f_{\tau_2(i)\tau_2(j)},\hat{f}_{\tau_2(i)\tau_2(j)}]} "C"; "F"};
\end{xy}\]
\end{rem}

Let $\Upsilon^{12}: h^1 \ra h^2$ and $\Upsilon^{23}: h^2 \ra h^3$ be two $2$-morphisms. We define the vertical composition $\Upsilon^{23}\circ_0 \Upsilon^{12}$ as follows. For each $i\in I$, consider the following diagram.
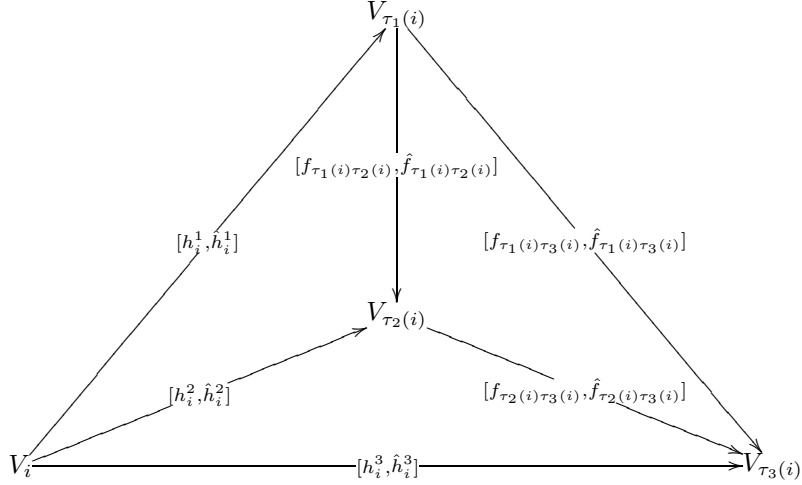
\[\begin{xy} 
(0,0)*{V_{i}}="A"; 
(50,60)*{V_{\tau_1(i)}}="B";
(100,0)*{V_{\tau_3(i)}}="C";
(50,20)*{V_{\tau_2(i)}}="D";
{\ar@{->}|-{[h^3_i,\hat{h}^3_i]} "A"; "C"};
{\ar@{->}|-{[h^1_i,\hat{h}^1_i]} "A"; "B"};
{\ar@{->}|-{[f_{\tau_1(i)\tau_3(i)},\hat{f}_{\tau_1(i)\tau_3(i)}]} "B"; "C"};
{\ar@{->}|-{[h^2_i,\hat{h}^2_i]} "A"; "D"};
{\ar@{->}|-{[f_{\tau_1(i)\tau_2(i)},\hat{f}_{\tau_1(i)\tau_2(i)}]} "B"; "D"};
{\ar@{->}|-{[f_{\tau_2(i)\tau_3(i)},\hat{f}_{\tau_2(i)\tau_3(i)}]} "D"; "C"};
\end{xy}\]
We define $[\Upsilon^{23}\circ_0\Upsilon^{12}]_i$ to be the sum of the three homotopies in the above picture:
\begin{equation}~\label{eq:vertical} 
[\Upsilon^{23}\circ_0\Upsilon^{12}]_i:= [\Upsilon^{23}]_i-\Lambda_{\tau_1(i)\tau_2(i)\tau_3(i)}*[h^1_i,\hat{h}^1_i]+[f_{\tau_2(i)\tau_3(i)},\hat{f}_{\tau_2(i)\tau_3(i)}]*[\Upsilon^{12}]_i.
\end{equation}
The $*$ composition on the right hand side is defined as in Subsection~\ref{subsec:hori}.

\medskip
\begin{lem}
The collection of homotopies $\left\{ [\Upsilon^{23}\circ_0 \Upsilon^{12}]_i \right\}$ defined above satisfies Equation~\ref{eq:2-equiv}. Thus $[\Upsilon^{23}\circ_0 \Upsilon^{12}]$ is a $2$-morphism from $h^1$ to $h^3$.
\end{lem}

\Pf. The proof is illustrated in the following diagram.
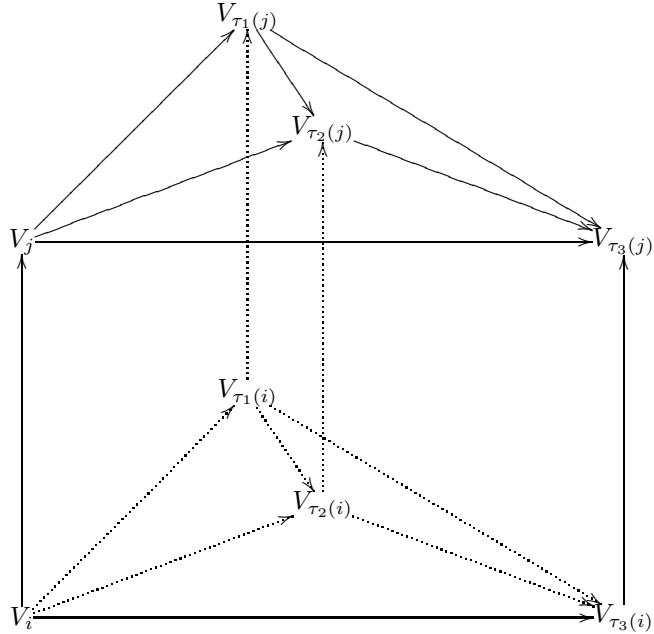
\[\begin{xy} 
(0,0)*{V_{i}}="A"; 
(30,30)*{V_{\tau_1(i)}}="B";
(80,0)*{V_{\tau_3(i)}}="C";
(40,15)*{V_{\tau_2(i)}}="D";
(0,50)*{V_{j}}="a";
(30,80)*{V_{\tau_1(j)}}="b";
(80,50)*{V_{\tau_3(j)}}="c";
(40,65)*{V_{\tau_2(j)}}="d";
{\ar@{->} "A"; "C"};
{\ar@{.>} "A"; "B"};
{\ar@{.>} "B"; "C"};
{\ar@{.>} "A"; "D"};
{\ar@{.>} "B"; "D"};
{\ar@{.>} "D"; "C"};
{\ar@{->} "a"; "c"};
{\ar@{->} "a"; "b"};
{\ar@{->} "b"; "c"};
{\ar@{->} "a"; "d"};
{\ar@{->} "b"; "d"};
{\ar@{->} "d"; "c"};
{\ar@{->} "A";"a"};
{\ar@{->} "C"; "c"};
{\ar@{.>} "B"; "b"};
{\ar@{.>} "D";"d"};
\end{xy}\]
Cancellations of the homotopies on the front and left back prisms follow from Equation~\ref{eq:2-equiv} for $\Upsilon^{12}$ and $\Upsilon^{23}$. Cancellation of the homotopies on the right back prism follows from Equation~\ref{eq:cocycle2}.\ed

\medskip
\begin{prop}
Every $2$-morphism $\Upsilon^{12}: h^1\ra h^2: \mathfrak{X}\ra \mathfrak{Y}$ is invertible with respect to the vertical composition.
\end{prop}

\Pf. We define the inverse $2$-morphism $\Upsilon^{21}: h^2 \ra h^1$ by 
\[ [\Upsilon^{21}]_i:= [\Lambda_{\tau_1(i)\tau_2(i)\tau_1(i)}]*[h^1_i,\hat{h}^1_i]-[f_{\tau_2(i)\tau_1(i)},\hat{f}_{\tau_2(i)\tau_1(i)}]*[\Upsilon^{12}]_i.\]
It is straight-forward to check the identities $\Upsilon^{12}\circ_0\Upsilon^{21}=0$ and $\Upsilon^{21}\circ_0\Upsilon^{12}=0$.\ed
 
Next, we define the horizontal compositions. For this let 
\begin{align*}
\Upsilon & : h^1\ra h^2: \mathfrak{X}\ra \mathfrak{Y},\\
\Gamma &: g^1\ra g^2: \mathfrak{Y}\ra \mathfrak{Z},
\end{align*}
be two horizontally composable $2$-morphisms. We use the notations 
\begin{align*}
h^1&=\Big\{\tau_1, [h^1_{i},\hat{h}^1_{i}], [\Delta^{h^1}_{ij}]\Big\}, h^2=\Big\{\tau_2, [h^2_{i},\hat{h}^2_{i}], [\Delta^{h^2}_{ij}]\Big\},\\
g^1&=\Big\{\sigma_1, [g^1_{i},\hat{g}^1_{i}], [\Delta^{g^1}_{ij}]\Big\}, g^2=\Big\{\sigma_2, [g^2_{i},\hat{g}^2_{i}], [\Delta^{g^2}_{ij}]\Big\}
\end{align*} 
for the structure maps of the morphisms. For an index $i\in I$, consider the following diagram
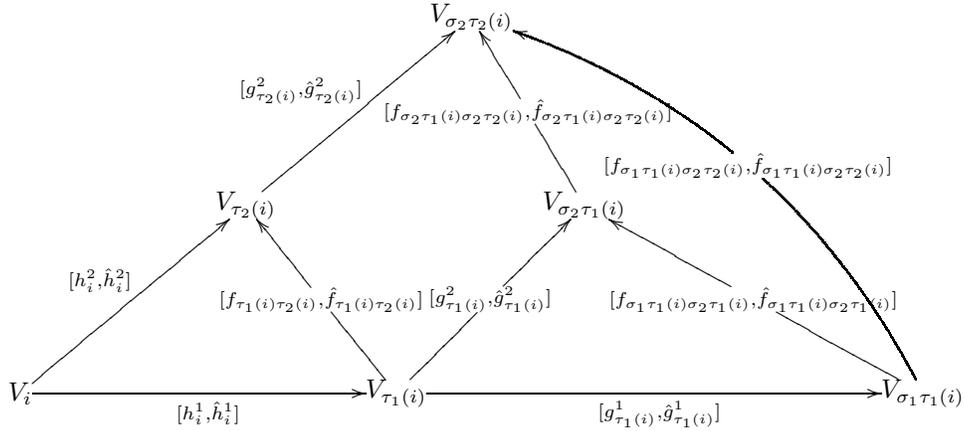
\[\begin{xy} 
(0,0)*{V_{i}}="A"; 
(50,0)*{V_{\tau_1(i)}}="B";
(30,25)*{V_{\tau_2(i)}}="C";
(120,0)*{V_{\sigma_1\tau_1(i)}}="D";
(75,25)*{V_{\sigma_2\tau_1(i)}}="E";
(60,50)*{V_{\sigma_2\tau_2(i)}}="F";
{\ar@{->}_{[h^1_i,\hat{h}^1_i]} "A"; "B"};
{\ar@{->}_{[g^1_{\tau_1(i)},\hat{g}^1_{\tau_1(i)}]} "B"; "D"};
{\ar@{->}^{[h^2_i,\hat{h}^2_i]} "A"; "C"};
{\ar@{->}|-{[f_{\tau_1(i)\tau_2(i)},\hat{f}_{\tau_1(i)\tau_2(i)}]} "B"; "C"};
{\ar@{->}|-{[g^2_{\tau_1(i)},\hat{g}^2_{\tau_1(i)}]} "B"; "E"};
{\ar@{->}^{[g^2_{\tau_2(i)},\hat{g}^2_{\tau_2(i)}]} "C"; "F"};
{\ar@{->}|-{[f_{\sigma_2\tau_1(i)\sigma_2\tau_2(i)},\hat{f}_{\sigma_2\tau_1(i)\sigma_2\tau_2(i)}]} "E"; "F"};
{\ar@{->}|-{[f_{\sigma_1\tau_1(i)\sigma_2\tau_1(i)},\hat{f}_{\sigma_1\tau_1(i)\sigma_2\tau_1(i)}]} "D"; "E"};
{\ar@{->}@/_2pc/ |-{[f_{\sigma_1\tau_1(i)\sigma_2\tau_2(i)},\hat{f}_{\sigma_1\tau_1(i)\sigma_2\tau_2(i)}]} "D";"F"};
\end{xy}\]
We define the the horizontal composition $\Gamma\circ_1 \Upsilon: g^1\circ h^1 \ra g^2\circ h^2$ by the sum of the three homotopies in the above diagram:
\begin{align}~\label{eq:horizontal} \begin{split}
[\Gamma\circ_1 \Upsilon]_i& :=  [g^2_{\tau_2(i)},\hat{g}^2_{\tau_2(i)}]*[\Upsilon_i]+[f_{\sigma_2\tau_1(i)\sigma_2\tau_2(i)},\hat{f}_{\sigma_2\tau_1(i)\sigma_2\tau_2(i)}]*[\Gamma_{\tau_1(i)}]*[h^1_i,\hat{h}^1_i]\\ &+[\Delta^{g^2}_{\tau_1(i)\tau_2(i)}]*[h^1_i,\hat{h}^1_i]-\Lambda_{\sigma_1\tau_1(i)\sigma_2\tau_1(i)\sigma_2\tau_2(i)}*[g^1_{\tau_1(i)},\hat{g}^1_{\tau_1(i)}]*[h^1_i,\hat{h}^1_i].
\end{split}\end{align}

\medskip
\begin{lem}
The collection of homotopies $\left\{ [\Gamma\circ_1 \Upsilon]_i \right\}$ defined above satisfies Equation~\ref{eq:2-equiv}. Thus it gives a $2$-morphism from $g^1\circ h^1$ to $g^2\circ h^2$.
\end{lem}

\Pf. Again, we illustrate the proof with a diagram:
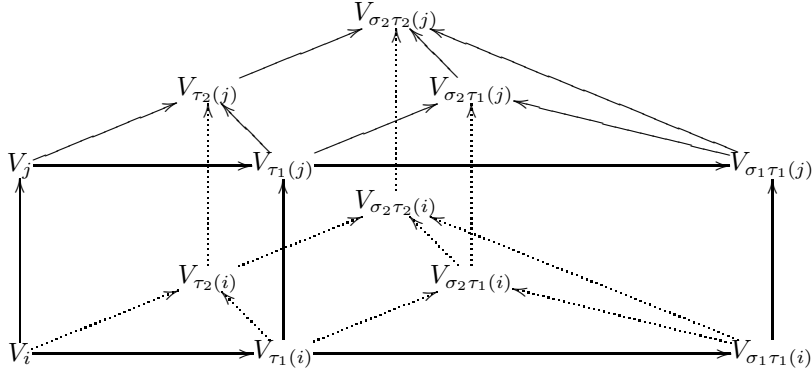
\[\begin{xy} 
(0,0)*{V_{i}}="A"; 
(35,0)*{V_{\tau_1(i)}}="B";
(25,10)*{V_{\tau_2(i)}}="C";
(100,0)*{V_{\sigma_1\tau_1(i)}}="D";
(60,10)*{V_{\sigma_2\tau_1(i)}}="E";
(50,20)*{V_{\sigma_2\tau_2(i)}}="F";
{\ar@{->} "A"; "B"};
{\ar@{->} "B"; "D"};
{\ar@{.>} "A"; "C"};
{\ar@{.>} "B"; "C"};
{\ar@{.>} "B"; "E"};
{\ar@{.>} "C"; "F"};
{\ar@{.>} "E"; "F"};
{\ar@{.>} "D"; "E"};
(0,25)*{V_{j}}="a"; 
(35,25)*{V_{\tau_1(j)}}="b";
(25,35)*{V_{\tau_2(j)}}="c";
(100,25)*{V_{\sigma_1\tau_1(j)}}="d";
(60,35)*{V_{\sigma_2\tau_1(j)}}="e";
(50,45)*{V_{\sigma_2\tau_2(j)}}="f";
{\ar@{->} "a"; "b"};
{\ar@{->} "b"; "d"};
{\ar@{->} "a"; "c"};
{\ar@{->} "b"; "c"};
{\ar@{->} "b"; "e"};
{\ar@{->} "c"; "f"};
{\ar@{->} "e"; "f"};
{\ar@{->} "d"; "e"};
{\ar@{->} "A"; "a"};
{\ar@{->} "B"; "b"};
{\ar@{->} "D"; "d"};
{\ar@{.>} "C"; "c"};
{\ar@{.>} "E"; "e"};
{\ar@{.>} "F"; "f"};
{\ar@{.>} "D";"F"};
{\ar@{->} "d";"f"};
\end{xy}\]
Cancellations of the homotopies on the two front prisms follow from Equation~\ref{eq:2-equiv} applied to $\Gamma$ and $\Upsilon$. The cancellation of homotopies on the back cube follows from Equation~\ref{eq:2-hom} applied to the morphism $g^2$. And on the back prism, we use Equation~\ref{eq:cocycle2}.\ed

With the above definition of horizontal and vertical compositions,  we have the following result. Its proof 
is straightforward, but lengthy. We put the proof in an appendix.

\medskip
\begin{thm}~\label{thm:pre}
There exists a $2$-category ${\sf pre-}\mathfrak{Kur}$ whose objects are Kuranishi manifolds (Definition~\ref{def:kuranishi}), $1$-morphisms are strict Kuranishi morphisms (Definition~\ref{def:homo}, and invertible $2$-morphisms as in Definition~\ref{def:equihom}.
\end{thm}

\subsection{Refinements} Let $\mathfrak{X}=\big( X, \left\{[V_i,\underline{\R^{m_i}},s_i,\psi_i]\right\}_{i\in I}$ be a Kuranishi manifold with the underlying topological space $X$. Let $\mathfrak{X}'=\big( X, \left\{[V_{i'},\underline{\R^{m_{i'}}},s_{i'},\psi_{i'}]\right\}_{i'\in I'}$ be another Kuranishi manifold with the same underlying topological space. 
\begin{defi}
A strict Kuranishi morphism morphism $r: \mathfrak{X}' \ra \mathfrak{X}$ is called a refinement if
\begin{itemize}
\item[(a.)] the map of indices $\iota: I'\ra I$ is surjective, and $$U_i=\bigcup_{i',\iota(i')=i} U_{i'}$$.
\item[(b.)] the morphism $r_{i'}: [V_{i'},\underline{\R^{m_{i'}}},s_{i'},\psi_{i'}]\ra [V_{\iota(i)},\underline{\R^{m_{\iota(i)}}},s_{\iota(i)},\psi_{\iota(i)}]$ is an open inclusion, for any $i'\in I$.
\item[(c.)] the homotopies $\Delta^r_{i'j'}=0$ for any indices $(i',j')\in I'\times I'$.
\end{itemize}
Note that it follows that the induced map $\underline{r}: X\ra X$ is equal to the identity map.
\end{defi}

Let $h=\Big\{\tau, [h_{i},\hat{h}_{i}], [\Delta_{ij}]\Big\}:\mathfrak{X}\ra \mathfrak{Y}$ be a strict Kuranishi morphism. Let $t=\Big\{\epsilon,[t_{p'},\hat{t}_{p'}, [\Delta_{p'q'}]\Big\}: \mathfrak{Y}'\ra \mathfrak{Y}$ be a refinement. We construct a strictly commutative diagram
\[\begin{CD}
\mathfrak{X}' @>r>> \mathfrak{X}\\
@Vh' VV  @VVh V\\
\mathfrak{Y}' @>t>> \mathfrak{Y}
\end{CD}\]
with $r: \mathfrak{X}'\ra \mathfrak{X}$ a refinement. Set the index set of $\mathfrak{X}'$ by
\[ I':= I\times_P P'=\left\{(i,p')|\tau(i)=\epsilon(p')\right\}.\]
For each index $(i,p')\in I'$, the Kuranishi chart is given by the restriction of $[V_{i},\underline{\R^{m_{i}}},s_{i},\psi_{i}]$ onto the open subset $\underline{h}^{-1}\big( U_{p'}\big)$. Note that part $(a.)$ in the above definition ensures that this is a covering of $X$. The morphism $h'_{(i,p')}$ is the restriction of $h_i$ onto $\underline{h}^{-1}\big( U_{p'}\big)$. Since all maps are defined by restrictions, the diagram is strictly commutative. We shall refer to this diagram as the {\sl canonical } pull-back of the strict Kuranishi morphism $h$ along the refinement $t$. Also we write $\mathfrak{X}'$ as
\[ \mathfrak{X}\times_{\mathfrak{Y}}\mathfrak{Y}'.\]
This should be understood only as a notation, as opposed to fiber product in categorical sense. The proof of the following Lemma is by elementary set theory of the index sets, and hence is omitted.

\medskip
\begin{lem}~\label{lem:refine}
The canonical pull-back construction satisfies the following properties:
\begin{itemize}
\item[(1.)] It is unital: $\mathfrak{X}\times_{\mathfrak{Y}}\mathfrak{Y}=\mathfrak{X}$.
\[\begin{CD}
\mathfrak{X}\times_{\mathfrak{Y}}\mathfrak{Y}=\mathfrak{X} @>\id>> \mathfrak{X}\\
@Vh VV  @VVh V\\
\mathfrak{Y} @>\id>> \mathfrak{Y}
\end{CD}\]
\item[(2.)] It is associative: $\big(\mathfrak{X}\times_{\mathfrak{Y}}\mathfrak{Y}'\big)\times_{\mathfrak{Z}}\mathfrak{Z}'=
\mathfrak{X}\times_{\mathfrak{Y}}\big(\mathfrak{Y}'\times_{\mathfrak{Z}}\mathfrak{Z}'\big)$.
\[\begin{xy}
(0,40)*{\big(\mathfrak{X}\times_{\mathfrak{Y}}\mathfrak{Y}'\big)\times_{\mathfrak{Z}}\mathfrak{Z}'}="A";
(30,40)*{\mathfrak{X}\times_{\mathfrak{Y}}\mathfrak{Y}'}="B";
(45,40)*{\mathfrak{X}}="C";
(0,0)*{\mathfrak{Z}'}="D";
(30,0)*{\mathfrak{Z}}="E";
(30,20)*{\mathfrak{Y}'}="F";
(45,20)*{\mathfrak{Y}}="G";
(65,40)*{\mathfrak{X}\times_{\mathfrak{Y}}\big(\mathfrak{Y}'\times_{\mathfrak{Z}}\mathfrak{Z}'\big)}="a";
(65,20)*{\mathfrak{Y}'\times_{\mathfrak{Z}}\mathfrak{Z}'}="b";
(110,40)*{\mathfrak{X}}="c";
(65,0)*{\mathfrak{Z}'}="d";
(95,0)*{\mathfrak{Z}}="e";
(95,20)*{\mathfrak{Y}'}="f";
(110,20)*{\mathfrak{Y}}="g";
(50,20)*{\cong};
{\ar@{->} "A";"B"};
{\ar@{->} "B";"C"};
{\ar@{->} "B";"F"};
{\ar@{->} "C";"G"};
{\ar@{->} "F";"G"};
{\ar@{->} "A"; "D"};
{\ar@{->} "D";"E"};
{\ar@{->} "F"; "E"};
{\ar@{->} "a";"c"};
{\ar@{->} "a";"b"};
{\ar@{->} "b";"d"};
{\ar@{->} "c";"g"};
{\ar@{->} "f";"g"};
{\ar@{->} "b"; "d"};
{\ar@{->} "d";"e"};
{\ar@{->} "f"; "e"};
{\ar@{->} "b"; "f"};
\end{xy}\]
\item[(3.)] It is symmetric: $\mathfrak{X}'\times_\mathfrak{X}\mathfrak{X}''=\mathfrak{X}''\times_\mathfrak{X}\mathfrak{X}'$.
\[\begin{CD}
\mathfrak{X}'\times_\mathfrak{X}\mathfrak{X}'' @>t'>> \mathfrak{X}' @. \;\;\; @. \mathfrak{X}''\times_\mathfrak{X}\mathfrak{X}' @>r'>> \mathfrak{X}''\\
@Vr' VV @VVr V @. @V t' VV @VVt V\\
\mathfrak{X}'' @>t>> \mathfrak{X} @. \;\;\; @.\mathfrak{X}' @>r>> \mathfrak{X}.
\end{CD}\]
Here, both $r$ and $t$ are refinements.
\end{itemize}
\end{lem}

\subsection{Localizing at refinements} In this subsection, we perform a $2$-categorical localization construction of ${\sf pre-}\mathfrak{Kur}$ at refinements. This gives us the desired $2$-category of Kuranishi manifolds.
\medskip
\begin{defi}~\label{defi:hom}
A Kuranishi morphism from $\mathfrak{X}$ to $\mathfrak{Y}$ is given by a roof diagram
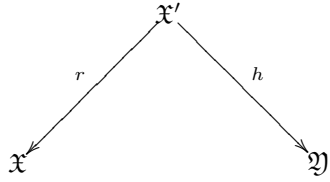
\[\begin{xy}
(0,0)*{\mathfrak{X}}="A";
(40,0)*{\mathfrak{Y}}="B";
(20,20)*{\mathfrak{X}'}="C";
{\ar@{->}_{r} "C"; "A"};
{\ar@{->}^{h} "C"; "B"};
\end{xy}\]
with $r$ a refinement, and $h$ a strict Kuranishi morphism as in Definition~\ref{def:homo}. Denote by $ hr^{-1}: \mathfrak{X}\ra \mathfrak{Y}$ this roof diagram.
\end{defi}

Note that in the ordinary localization construction of a category, one defines a morphism to be an equivalence class of roof diagrams. To perform a $2$-categorical localization, we work with honest roof diagrams, and introduce $2$-morphisms to encode the classical equivalence relation between roof diagrams. Before doing that, we first prove the existence of an associative composition among roof diagrams.

The composition between two roof diagrams $hr^{-1}: \mathfrak{X}\ra \mathfrak{Y}$ and $gt^{-1}:\mathfrak{Y}\ra\mathfrak{Z}$ is given by the following diagram
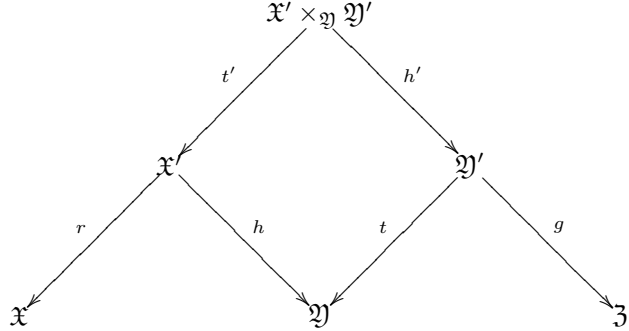
\[\begin{xy}
(0,0)*{\mathfrak{X}}="A";
(40,0)*{\mathfrak{Y}}="B";
(20,20)*{\mathfrak{X}'}="C";
(80,0)*{\mathfrak{Z}}="D";
(60,20)*{\mathfrak{Y}'}="E";
(40,40)*{\mathfrak{X}'\times_{\mathfrak{Y}}\mathfrak{Y}'}="F";
{\ar@{->}_{r} "C"; "A"};
{\ar@{->}^{h} "C"; "B"};
{\ar@{->}_{t} "E"; "B"};
{\ar@{->}^{g} "E"; "D"};
{\ar@{->}_{t'} "F"; "C"};
{\ar@{->}^{h'} "F"; "E"};
\end{xy}\]
where the top square is given by the canonical pull-back of $h$ along $t$.

\medskip
\begin{prop}
The composition of roof diagrams defined above is unital and associative.
\end{prop}

\Pf. The unit morphism of an object $\mathfrak{X}$ is given by $\id\id^{-1}: \mathfrak{X}\ra\mathfrak{X}$. The Proposition follows immediately from Lemma~\ref{lem:refine} and the associativity of $1$-morphisms in ${\sf pre-}\mathfrak{Kur}$ proved in Proposition~\ref{prop:pre-comp}.\ed

Let $\mathfrak{R}$ denote the set of all refinements. We have obtained, by localizing at $\mathfrak{R}$, an ordinary category whose objects are the same as ${\sf pre-}\mathfrak{Kur}$ and whose morphisms are roof diagrams. We denote this category by $\mathfrak{Kur}$. Since ${\sf pre-}\mathfrak{Kur}$ is a $2$-category, it is natural to ask for a $2$-category localization. In the remaining part of the subsection, we show that indeed $\mathfrak{Kur}$ admits a natural $2$-category structure. 

\begin{defi}\label{defi:2-hom}
A homotopy $\chi^{12}: (h^1)(r^1)^{-1} \ra (h^2)(r^2)^{-1}: \mathfrak{X}\ra \mathfrak{Y}$ between two roof diagrams is given by the following data (which is also illustrated in the diagram below):
\begin{itemize}
\item[--] Refinements $r_\chi^{1}: \mathfrak{X}_\chi \ra \mathfrak{X}^1$ and $r^{2}_\chi: \mathfrak{X}_\chi \ra \mathfrak{X}^2$ such that $r^{1}\circ r^{1}_\chi=r^{2}\circ r^{2}_\chi$.
\item[--] A $2$-morphism $\chi^{12}:    h^1r^{1}_\chi\ra h^2r^{2}_\chi : \mathfrak{X}_\chi \ra \mathfrak{Y}$ in the sense of Definition~\ref{def:equihom}.
\end{itemize}
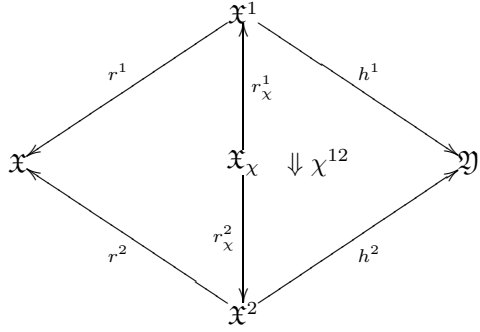
\[\begin{xy}
(0,0)*{\mathfrak{X}}="A";
(30,0)*{\mathfrak{X}_\chi}="B";
(60,0)*{\mathfrak{Y}}="C";
(40,0)*{\Downarrow\chi^{12}};
(30,20)*{\mathfrak{X}^1}="D";
(30,-20)*{\mathfrak{X}^2}="E";
{\ar@{->}_{r^1} "D"; "A"};
{\ar@{->}^{r^2} "E"; "A"};
{\ar@{->}^{h^1} "D"; "C"};
{\ar@{->}_{h^2} "E"; "C"};
{\ar@{->}_{r^{1}_\chi} "B"; "D"};
{\ar@{->}_{r^{2}_\chi} "B"; "E"};
\end{xy}\]
\end{defi}

\medskip
\noindent {\bf Notation.} Let $\chi^{12}$ be a homotopy as in the above definition. Let $r:\mathfrak{X}_\chi'\ra \mathfrak{X}_\chi$ be another refinement. We get an induced homotopy illustrated in the following diagram. 
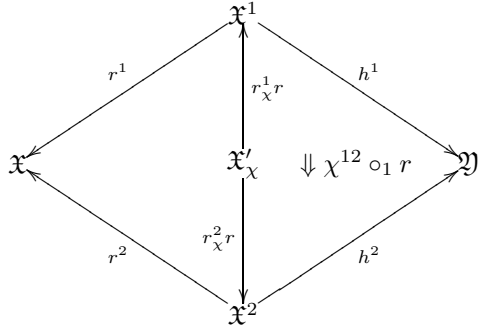
\[\begin{xy}
(0,0)*{\mathfrak{X}}="A";
(30,0)*{\mathfrak{X}_\chi'}="B";
(60,0)*{\mathfrak{Y}}="C";
(45,0)*{\Downarrow\chi^{12}\circ_1 r};
(30,20)*{\mathfrak{X}^1}="D";
(30,-20)*{\mathfrak{X}^2}="E";
{\ar@{->}_{r^1} "D"; "A"};
{\ar@{->}^{r^2} "E"; "A"};
{\ar@{->}^{h^1} "D"; "C"};
{\ar@{->}_{h^2} "E"; "C"};
{\ar@{->}_{r_\chi^{1}r} "B"; "D"};
{\ar@{->}_{r^{2}_\chi r} "B"; "E"};
\end{xy}\]
Here $\circ_1$ is the horizontal composition in ${\sf pre-}\mathfrak{Kur}$. We refer to this new homotopy $\chi^{12}\circ_1 r$ as the restriction of $\chi^{12}$ onto $\mathfrak{X}_\chi'$. We shall simply write $\chi^{12}|_{\mathfrak{X}_\chi'}$ for this homotopy, if the refinement $r:\mathfrak{X}_\chi'\ra \mathfrak{X}_\chi$ is clear from the context.

\medskip
Two homotopies $\chi^{12}$, $\lambda^{12}$ are called equivalent if there exists a common refinement $\mathfrak{X}'_{\chi\lambda}$ of $\mathfrak{X}_\chi$ and $\mathfrak{X}_\lambda$ that fit into the following diagram.
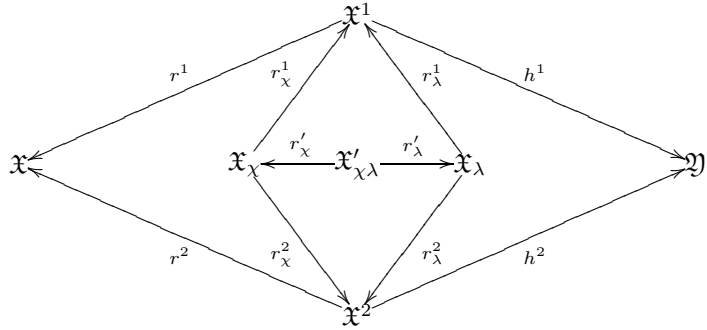
\[\begin{xy}
(0,0)*{\mathfrak{X}}="A";
(30,0)*{\mathfrak{X}_\chi}="B";
(45,0)*{\mathfrak{X}_{\chi\lambda}'}="F";
(60,0)*{\mathfrak{X}_\lambda}="G";
(90,0)*{\mathfrak{Y}}="C";
(45,20)*{\mathfrak{X}^1}="D";
(45,-20)*{\mathfrak{X}^2}="E";
{\ar@{->}_{r^1} "D"; "A"};
{\ar@{->}^{r^2} "E"; "A"};
{\ar@{->}^{h^1} "D"; "C"};
{\ar@{->}_{h^2} "E"; "C"};
{\ar@{->}^{r_\chi^{1}} "B"; "D"};
{\ar@{->}_{r^{2}_\chi} "B"; "E"};
{\ar@{->}_{r_\lambda^{1}} "G"; "D"};
{\ar@{->}^{r_\lambda^{2}} "G"; "E"};
{\ar@{->}^{r_\lambda'} "F"; "G"};
{\ar@{->}_{r_\chi'} "F"; "B"};
\end{xy}\]
In this diagram, compositions of refinements are strictly commutative. And we require that the two restricted homotopies are equal:
\[ \chi^{12}\circ_1 r_\chi'=\lambda^{12}\circ_1 r_\lambda'.\]
\begin{lem}
The equivalences between homotopies defined above form an equivalence relation on the set of homotopies from $(h^1)(r^1)^{-1}$ to $(h^2)(r^2)^{-1}$.
\end{lem}

\Pf. Only the transitivity is non-trivial. Its proof is illustrated in the following diagram. 
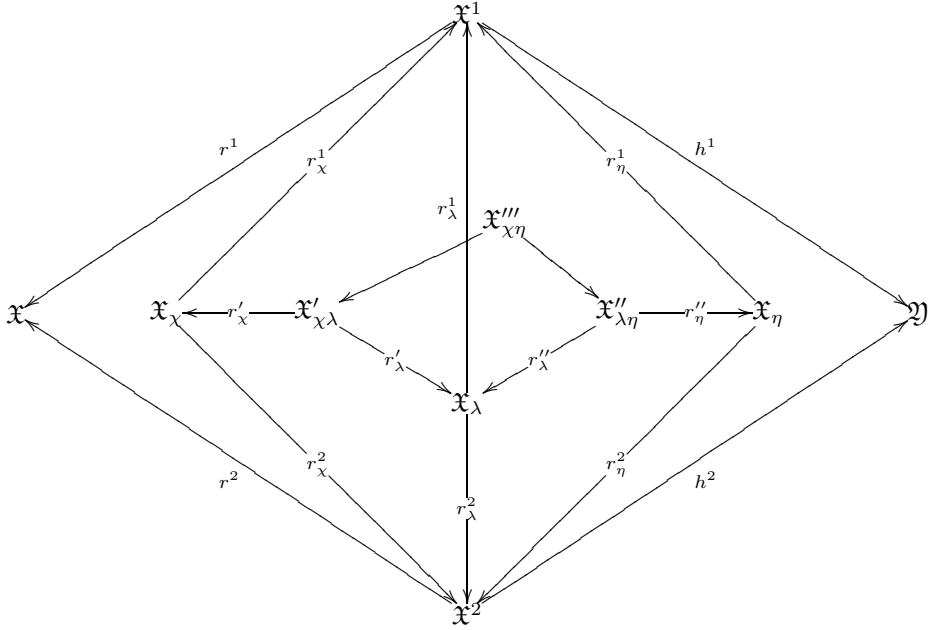
\[\begin{xy}
(0,0)*{\mathfrak{X}}="A";
(20,0)*{\mathfrak{X}_\chi}="B";
(40,0)*{\mathfrak{X}_{\chi\lambda}'}="F";
(60,-12)*{\mathfrak{X}_\lambda}="G";
(120,0)*{\mathfrak{Y}}="C";
(60,40)*{\mathfrak{X}^1}="D";
(60,-40)*{\mathfrak{X}^2}="E";
(80,0)*{\mathfrak{X}_{\lambda\eta}''}="H";
(100,0)*{\mathfrak{X}_\eta}="I";
(65,12)*{\mathfrak{X}_{\chi\eta}'''}="J";
{\ar@{->}_{r^1} "D"; "A"};
{\ar@{->}^{r^2} "E"; "A"};
{\ar@{->}^{h^1} "D"; "C"};
{\ar@{->}_{h^2} "E"; "C"};
{\ar@{->}|-{r_\chi^{1}} "B"; "D"};
{\ar@{->}|-{r^{2}_\chi} "B"; "E"};
{\ar@{->}^{r_\lambda^{1}} "G"; "D"};
{\ar@{->}|-{r_\lambda^{2}} "G"; "E"};
{\ar@{->}|-{r_\lambda'} "F"; "G"};
{\ar@{->}|-{r_\chi'} "F"; "B"};
{\ar@{->}|-{r_\lambda''} "H"; "G"};
{\ar@{->}|-{r_\eta''} "H"; "I"};
{\ar@{->}|-{r_\eta^{2}} "I"; "E"};
{\ar@{->}|-{r_\eta^{1}} "I"; "D"};
{\ar@{->}@[red] "J"; "F"};
{\ar@{->}@[red] "J"; "H"};
\end{xy}\]
The idea is to exhibit an equivalence between $\chi^{12}$ and $\eta^{12}$ by considering the common refinement
\[ \mathfrak{X}_{\chi\eta}''':= \mathfrak{X}_{\chi\lambda}'\times_{\mathfrak{X}_\lambda} \mathfrak{X}_{\lambda\eta}''.\]
By assumption, the restrictions of the homotopies $\chi^{12}$ and $\lambda^{12}$ onto $\mathfrak{X}_{\chi\eta}'''$ are equal, so are the restrictions of $\lambda^{12}$ and $\eta^{12}$. This implies that the restrictions of $\chi^{12}$ and $\eta^{12}$ onto $\mathfrak{X}_{\chi\eta}'''$ are also equal.\ed

\begin{defi}~\label{defi:2-mor}
Let $(h^1)(r^1)^{-1}, (h^2)(r^2)^{-1}: \mathfrak{X}\ra \mathfrak{Y}$  be morphisms in the category $\mathfrak{Kur}$.  A $2$-morphism from $(h^1)(r^1)^{-1}$ to $(h^2)(r^2)^{-1}$ is given by an equivalence class of homotopies from $(h^1)(r^1)^{-1}$ to $(h^2)(r^2)^{-1}$ (Definition~\ref{defi:2-hom}).
\end{defi}

Let $[\chi^{12}]: (h^1)(r^1)^{-1} \ra (h^2)(r^2)^{-1}$ and $[\chi^{23}]: (h^2)(r^2)^{-1} \ra (h^3)(r^3)^{-1}$ be composible $2$-morphisms. We define their vertical composition using the following diagram.

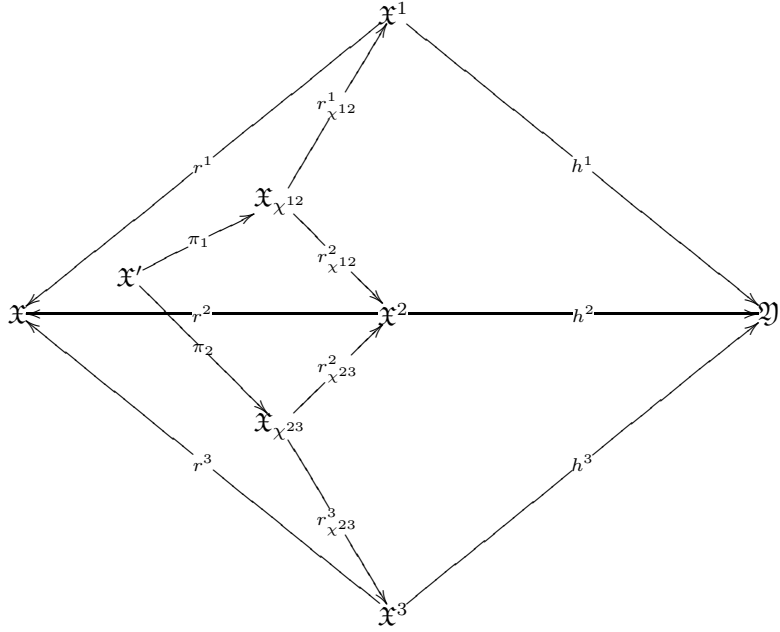
\[\begin{xy}
(0,0)*{\mathfrak{X}}="A";
(50,0)*{\mathfrak{X}^2}="B";
(100,0)*{\mathfrak{Y}}="C";
(50,40)*{\mathfrak{X}^1}="D";
(35,15)*{\mathfrak{X}_{\chi^{12}}}="F";
(15,5)*{\mathfrak{X}'}="H";
(35,-15)*{\mathfrak{X}_{\chi^{23}}}="G";
(50,-40)*{\mathfrak{X}^3}="E";
{\ar@{->}|-{r^2} "B"; "A"};
{\ar@{->}|-{r^1} "D"; "A"};
{\ar@{->}|-{r^3} "E"; "A"};
{\ar@{->}|-{h^1} "D"; "C"};
{\ar@{->}|-{h^2} "B"; "C"};
{\ar@{->}|-{h^3} "E"; "C"};
{\ar@{->}|-{r^2_{\chi^{12}}} "F"; "B"};
{\ar@{->}|-{r^1_{\chi^{12}}} "F"; "D"};
{\ar@{->}|-{r^2_{\chi^{23}}} "G"; "B"};
{\ar@{->}|-{r^3_{\chi^{23}}} "G"; "E"};
{\ar@{->}@[red]|-{\pi_1} "H"; "F"};
{\ar@{->}@[red]|-{\pi_2} "H"; "G"};
\end{xy}\]
In the above diagram, the notation 
\[\mathfrak{X}':=\mathfrak{X}_{\chi^{12}}\times_{\mathfrak{X}^2}\mathfrak{X}_{\chi^{23}}\]
is defined by the canonical pull-back construction. Since both $r^2_{\chi^{12}}$ and $r^2_{\chi^{23}}$ are refinements, both $\pi_1$ and $\pi_2$ are refinements. We then define the vertical composition of $[\chi^{12}]$ and $[\chi^{23}]$ by the equivalence class of the vertical composition of $\chi^{12}$ and $\chi^{23}$ restricted onto the common refinement $\mathfrak{X}'$:
\begin{equation}\label{eq-def-vert}
 [\chi^{23}]\circ_0 [\chi^{12}]:= [ \big( \chi^{23}|_{\mathfrak{X}'}\big)\circ_0\big( \chi^{12}|_{\mathfrak{X}'}\big)].
 \end{equation}
Here the compositions $\circ_1$ and $\circ_0$ are horizontal and vertical compositions in ${\sf pre-}\mathfrak{Kur}$.

\medskip
Next, we define the horizontal compositions of $2$-morphisms
\begin{align*}
& \chi: (h^1)(r^1)^{-1}, (h^2)(r^2)^{-1}: \mathfrak{X}\ra \mathfrak{Y}\\
& \lambda: (g^1)(t^1)^{-1}, (g^2)(t^2)^{-1}: \mathfrak{Y}\ra \mathfrak{Z}.
\end{align*}
Again, we illustrate the construction diagrammatically. In the following diagram, any square that has an red edge is a canonical pull-back square.

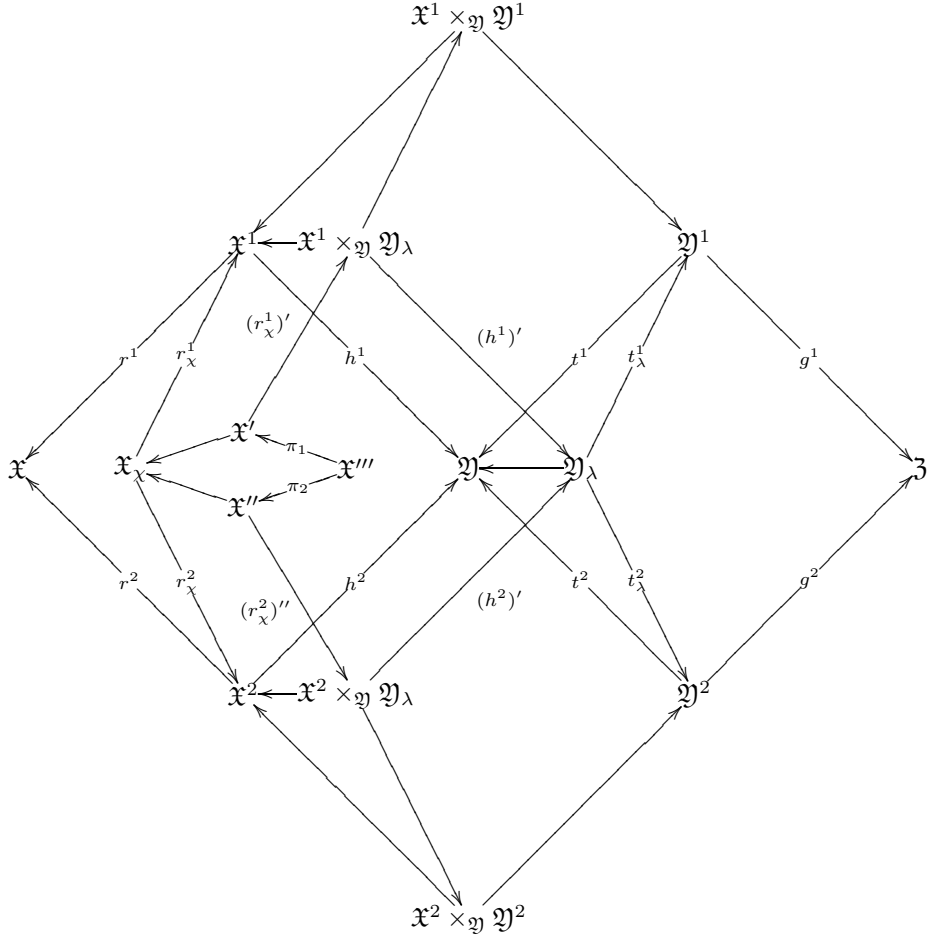
\[\begin{xy}
(0,0)*{\mathfrak{X}}="A";
(60,0)*{\mathfrak{Y}}="B";
(120,0)*{\mathfrak{Z}}="C";
(15,0)*{\mathfrak{X}_\chi}="D";
(75,0)*{\mathfrak{Y}_\lambda}="E";
(30,30)*{\mathfrak{X}^1}="F";
(30,5)*{\mathfrak{X}'}="L";
(30,-5)*{\mathfrak{X}''}="O";
(45,0)*{\mathfrak{X}'''}="P";
(45,30)*{\mathfrak{X}^1\times_{\mathfrak{Y}}\mathfrak{Y}_\lambda}="J";
(90,30)*{\mathfrak{Y}^1}="G";
(30,-30)*{\mathfrak{X}^2}="H";
(45,-30)*{\mathfrak{X}^2\times_{\mathfrak{Y}}\mathfrak{Y}_\lambda}="K";
(90,-30)*{\mathfrak{Y}^2}="I";
(60,60)*{\mathfrak{X}^1\times_{\mathfrak{Y}}\mathfrak{Y}^1}="M";
(60,-60)*{\mathfrak{X}^2\times_{\mathfrak{Y}}\mathfrak{Y}^2}="N";
{\ar@{->}|-{r^1} "F"; "A"};
{\ar@{->}|-{r^2} "H"; "A"};
{\ar@{->}|-{r^1_\chi} "D"; "F"};
{\ar@{->}|-{r^2_\chi} "D"; "H"};
{\ar@{->}|-{h^1} "F"; "B"};
{\ar@{->}|-{h^2} "H"; "B"};
{\ar@{->}|-{t^1} "G"; "B"};
{\ar@{->}|-{t^2} "I"; "B"};
{\ar@{->}|-{g^1} "G"; "C"};
{\ar@{->}|-{g^2} "I"; "C"};
{\ar@{->}|-{t^1_\lambda} "E"; "G"};
{\ar@{->}|-{t^2_\lambda} "E"; "I"};
{\ar@{->} "E"; "B"};
{\ar@{->} "J"; "F"};
{\ar@{->}@[red]^{(h^1)'} "J"; "E"};
{\ar@{->} "K"; "H"};
{\ar@{->}@[red]_{(h^2)'} "K"; "E"};
{\ar@{->} "L"; "D"};
{\ar@{->} "O"; "D"};
{\ar@{->}@[red]^{(r^1_\chi)'} "L"; "J"};
{\ar@{->}@[red]_{(r^2_\chi)''} "O"; "K"};
{\ar@{->}@[red]|-{\pi_1} "P"; "L"};
{\ar@{->}@[red]|-{\pi_2} "P"; "O"};
{\ar@{->} "M"; "F"};
{\ar@{->} "M"; "G"};
{\ar@{->} "J"; "M"};
{\ar@{->} "N"; "H"};
{\ar@{->} "N"; "I"};
{\ar@{->} "K"; "N"};
\end{xy}\]
Since $\mathfrak{Y}_\lambda\ra \mathfrak{Y}$ is a refinement, the restriction of $\chi$ onto $\mathfrak{X}'''$ may be viewed as a $2$-morphism (in ${\sf pre-}\mathfrak{Kur}$)
\[ \chi|_{\mathfrak{X}'''}:    (h^1)'(r^1_\chi)' \pi_1\ra (h^2)'(r^2_\chi)'\pi_2.\]
Then we define
\begin{equation}~\label{eq-def-hori}
 [\lambda]\circ_1 [\chi]:=[ \lambda\circ_1\big(\chi|_{\mathfrak{X}'''}\big)].
 \end{equation}
The following theorem is the main result of the paper. With the above constructions at hand, it is essentially a corollary of Theorem~\ref{thm:pre}. 

\medskip
\begin{thm}~\label{thm:main}
There exists a $2$-category structure on $\mathfrak{Kur}$ whose objects are Kuranishi manifolds, morphisms are roof diagrams~\ref{defi:hom}, and invertible $2$-morphisms defined as in~\ref{defi:2-mor}, with the vertical and horizontal compositions defined by Equations~\ref{eq-def-vert} and~\ref{eq-def-hori}. 
\end{thm}

\Pf. Like in the proof of Theorem~\ref{thm:pre}, we need to check associativity of $\circ_0$, $\circ_1$ and the interchange law identity. By choosing a common refinement, thee proof follows from the corresponding properties of the category ${\sf pre-}\mathfrak{Kur}$. We illustrate this argument by proving the vertical composition is associative. Indeed, according to Equation~\ref{eq-def-vert}, we have
\begin{align*}
[\chi^{23}]\circ_0 [\chi^{12}]& = [ \big( \chi^{23}|_{\mathfrak{X}'}\big)\circ_0\big( \chi^{12}|_{\mathfrak{X}'}\big)]\\
[\chi^{34}]\circ_0\big([\chi^{23}]\circ_0 [\chi^{12}]\big)&= [ \big(\chi^{34}|_{\mathfrak{X}''}\big) \circ_0 \big( \chi^{23}|_{\mathfrak{X}''} \circ_0 \chi^{12}|_{\mathfrak{X}''}\big)]
\end{align*}
Here $\mathfrak{X}''$  is a refinement of $\mathfrak{X}$. The other composition would look like
\[ \big([\chi^{34}]\circ_0[\chi^{23}]\big)\circ_0 [\chi^{12}]= [ \big(\chi^{34}|_{\mathfrak{X}'''} \circ_0 \chi^{23}|_{\mathfrak{X}'''} \big) \circ_0 \big(\chi^{12}|_{\mathfrak{X}'''}\big)].\]
Here $\mathfrak{X}'''$ is another refinement of $\mathfrak{X}$. Now to verify that the two compositions are equal, we restrict everything to the common refinement $\mathfrak{X}''\times_{\mathfrak{X}}\mathfrak{X}'''$ using the canonical pull-back construction. We then apply the associativity in ${\sf pre-}\mathfrak{Kur}$ to finish the proof. Other identities are verified in the same way.\ed

\section{Fiber products over manifolds}~\label{sec:fiber}

In this section, we prove certain $2$-fiber product property holds in $\mathfrak{Kur}$.

\subsection{The homotopy category of $\mathfrak{Kur}$} Let $\mathfrak{X}$ and $\mathfrak{Y}$ be two Kuranishi manifolds, the category
\[ \Hom_{\mathfrak{Kur}}(\mathfrak{X}, \mathfrak{Y})\]
is in fact a groupoid since all its morphisms are invertible. Taking $\pi_0$ of these Hom-groupoids, we obtain the homotopy category $\pi_0(\mathfrak{Kur})$.

\medskip
\begin{prop}
The category of smooth manifolds embeds fully faithfully into $\pi_0(\mathfrak{Kur})$.
\end{prop}

\Pf. By definition, a manifold is a topological space $X$, together with an equivalence class of smooth atlas on $X$. Choose a representative atlas $\AA$ in its equivalence class. We obtain a Kuranishi manifold $(X,\AA)$. Let $(X,\AA)$ and $(Y,\BB)$ be two such Kuranishi manifolds. A smooth map between the two manifolds may not be written in the chosen atlases $\AA$ and $\BB$. But there always exists a refinement atlas $\AA'$ of $\AA$ such that the smooth map is represented by a roof diagram of the form
\[ (X,\AA) \leftarrow (X,\AA') \ra (Y,\BB).\]
One easily checks that this gives a fully faithful embedding of manifolds into $\pi_0(\mathfrak{Kur})$, because the Kuranishi atlases defined from manifolds have no obstruction bundles which forces the $2$-morphisms among them to be the trivial.\ed

\subsection{Fiber products over manifolds} Recall the definition of $2$-fiber products in a $2$-category.

\medskip
\begin{defi}~\label{defi:2-fiber}
Let $\CC$ be a $2$-category. A diagram
\[\begin{xy}
(0,0)*{B}="B";
(15,0)*{A}="A";
(15,15)*{C}="C";
(0,15)*{D}="D";
(5,5)*{\;}="H";
(10,10)*{\;}="G";
{\ar@{=>}^{\theta} "H"; "G"};
{\ar@{->}_{h} "B"; "A"};
{\ar@{->}^{g} "C"; "A"};
{\ar@{->}_{\pi_1} "D"; "B"};
{\ar@{->}^{\pi_2} "D"; "C"};
\end{xy}\]
is called a $2$-fiber diagram if 
\begin{itemize}
\item[(A.)] For any other such diagram
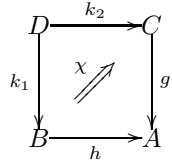
\[\begin{xy}
(0,0)*{B}="B";
(15,0)*{A}="A";
(15,15)*{C}="C";
(0,15)*{D}="E";
(5,5)*{\;}="H";
(10,10)*{\;}="G";
{\ar@{=>}^{\chi} "H"; "G"};
{\ar@{->}_{h} "B"; "A"};
{\ar@{->}^{g} "C"; "A"};
{\ar@{->}_{k_1} "E"; "B"};
{\ar@{->}^{k_2} "E"; "C"};
\end{xy}\]
there exists a morphism $u: E\ra D$ which fits into the following diagram
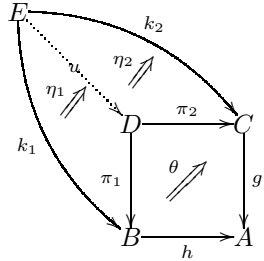
\[\begin{xy}
(-15,30)*{E}="E";
(0,0)*{B}="B";
(15,0)*{A}="A";
(15,15)*{C}="C";
(0,15)*{D}="D";
(5,5)*{\;}="H";
(10,10)*{\;}="G";
(-9,16)*{\;}="H_1";
(-6, 20)*{\;}="G_1";
(0,20)*{\;}="H_2";
(3,24)*{\;}="G_2";
{\ar@{=>}^{\theta} "H"; "G"};
{\ar@{=>}^{\eta_1} "H_1"; "G_1"};
{\ar@{=>}^{\eta_2} "H_2"; "G_2"};
{\ar@{->}_{h} "B"; "A"};
{\ar@{->}^{g} "C"; "A"};
{\ar@{->}_{\pi_1} "D"; "B"};
{\ar@{->}^{\pi_2} "D"; "C"};
{\ar@{.>}|-{u} "E"; "D"};
{\ar@{->}@/_1pc/_{k_1} "E"; "B"};
{\ar@{->}@/^1pc/^{k_2} "E"; "C"};
\end{xy}\]
such that 
\begin{equation}~\label{eq:2-fiber-1}
\chi=\big( g\circ_1\eta_2\big)\circ_0\big(\theta\circ_1 u \big)\circ_0\big(h\circ_1\eta_1\big).
\end{equation}
\item[(B.)] For any other $u': E\ra D$ satisfying part $(A.)$ above, there exists a unique $2$-morphism $\lambda: u\ra u'$ such that 
\begin{align}\label{eq:2-fiber-2}
\begin{split}
\eta_2 &=(\eta_2')\circ_0\big(\pi_2\circ_1\lambda\big)\\ \eta_1'&=\big(\pi_1\circ_1\lambda\big)\circ_0(\eta_1),
\end{split}
\end{align}
as illustrated in the following diagram.
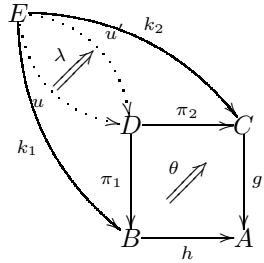
\[\begin{xy}
(-15,30)*{E}="E";
(0,0)*{B}="B";
(15,0)*{A}="A";
(15,15)*{C}="C";
(0,15)*{D}="D";
(5,5)*{\;}="H";
(10,10)*{\;}="G";
(-10,20)*{\;}="H_1";
(-5,25)*{\;}="G_2";
{\ar@{=>}^{\theta} "H"; "G"};
{\ar@{=>}^{\lambda} "H_1"; "G_2"};
{\ar@{.>}@/^1pc/^{u'} "E"; "D"};
{\ar@{->}_{h} "B"; "A"};
{\ar@{->}^{g} "C"; "A"};
{\ar@{->}_{\pi_1} "D"; "B"};
{\ar@{->}^{\pi_2} "D"; "C"};
{\ar@{.>}@/_1pc/_{u} "E"; "D"};
{\ar@{->}@/_1pc/_{k_1} "E"; "B"};
{\ar@{->}@/^1pc/^{k_2} "E"; "C"};
\end{xy}\]
\end{itemize}
\end{defi}

Now, suppose we are given the following diagram of morphisms in $\mathfrak{Kur}$ (where $\mathfrak{Z}$ and $\pi_1$, $\pi_2$ are to-be-constructed):
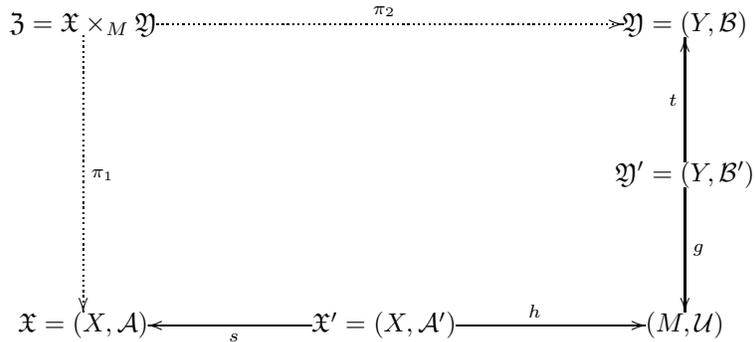
\[ \begin{xy}
(0,40)*{\mathfrak{Z}=\mathfrak{X}\times_M \mathfrak{Y}}="F";
(80,40)*{\mathfrak{Y}=(Y,\BB)}="A";
(80,20)*{\mathfrak{Y}'=(Y,\BB')}="B";
(80,0)*{(M,\UU)}="C";
(0,0)*{\mathfrak{X}=(X,\AA)}="D";
(40,0)*{\mathfrak{X}'=(X,\AA')}="E";
{\ar@{->}^t "B"; "A"};
{\ar@{->}^g "B"; "C"};
{\ar@{->}^s "E"; "D"};
{\ar@{->}^h "E"; "C"};
{\ar@{.>}^{\pi_1} "F"; "D"};
{\ar@{.>}^{\pi_2} "F"; "A"};
\end{xy}\] 
Assume that the Kuranishi manifold $(M,\UU)$ is an ordinary manifold. And other Kuranishi atlases are given by
\begin{align*}
\AA & :=\big\{ [V_i,\underline{\R^{m_i}}, s_i, \psi_i], [f_{ij},\hat{f}_{ij}], [\Lambda_{ijk}]\big\}_{i,j,k\in I}\\
\AA' & :=\big\{ [V_{i'},\underline{\R^{m_{i'}}}, s_{i'}, \psi_{i'}], [f_{i'j'},\hat{f}_{i'j'}], [\Lambda_{i'j'k'}]\big\}_{i',j',k'\in I'}\\
\BB & := \big\{ [V_p,\underline{\R^{m_p}}, s_p, \psi_p], [f_{pq},\hat{f}_{pq}], [\Lambda_{pqr}] \big\}_{p,q,r\in P}\\
\BB' & := \big\{ [V_{p'},\underline{\R^{m_{p'}}}, s_{p'}, \psi_{p'}], [f_{p'q'},\hat{f}_{p'q'}], [\Lambda_{p'q'r'}] \big\}_{p',q',r'\in P'}
\end{align*}
In this subsection, we prove the existence of a $2$-fiber product Kuranishi manifold $\mathfrak{X}\times_M \mathfrak{Y}$ in the $2$-category $\mathfrak{Kur}$. For the first step, we proceed to construct the following data:
\begin{itemize}
\item[(1.)] A Kuranishi manifold $\mathfrak{Z}$ (the to be fiber product $\mathfrak{X}\times_M \mathfrak{Y}$).
\item[(2.)] Strict Kuranishi morphisms $\pi_1: \mathfrak{Z} \ra \mathfrak{X}$ and $\pi_2:\mathfrak{Z}\ra \mathfrak{Y}$.
\item[(3.)] A $2$-morphism $[\theta]:      (hs^{-1})(\pi_1)\ra (gt^{-1})(\pi_2)$.
\end{itemize}
We describe the constructions in the following steps.
 
 \medskip
(1a.) The underlying topological space of $\mathfrak{X}\times_M \mathfrak{Y}$ is simply the fiber product $Z:=X\times_M Y$ in the category of topological spaces, i.e.
\[ Z:=\left\{ (x,y)| \underline{h}(x)=\underline{g}(y)\right\}.\]

 \medskip
(1b.) We define the index set $D$ of the to-be-constructed fiber product by
\[ D:=\left\{ (i',p') \in I'\times P'\mid \underline{h}(U_{i'})\cap \underline{g}(U_{p'})\neq \emptyset.\right\}\]
Let $\tau: I'\ra L$ and $\eta: P'\ra L$ be the maps of indices associated to the two strict Kuranishi morphisms $h: \mathfrak{X}'\ra M$ and $g: \mathfrak{Y}'\ra M$ respectively. Fixing an index $d=(i',p')\in D$, we have the following setup:
 \begin{align*}
& U_{i'} \stackrel{\psi_{i'}}{\ra} V_{i'} \stackrel{h_{i'}}{\ra} V_{\tau(i')}\subset \R^n,\\
& U_{p'} \stackrel{\psi_{p'}}{\ra} V_{p'} \stackrel{h_{p'}}{\ra} V_{\eta(p')}\subset \R^n.
 \end{align*}
We also have the transition map (on $M$)
\[ f_{\tau(i')\eta(p')}: V_{\underline{\tau(i')}\eta(p')} \ra V_{\tau(i')\underline{\eta(p')}},\]
defined on some open subsets $V_{\underline{\tau(i')}\eta(p')}\subset V_{\tau(i')}$ and $V_{\tau(i')\underline{\eta(p')}}\subset V_{\eta(p')}$ with footprint $U_{\tau(i')}\cap U_{\eta(p')}$. Define
\begin{align*}
 V_{d,X} &:=h_{i'}^{-1}\big( V_{\underline{\tau(i')}\eta(p')} \big) \\
U_{d,X} &:= \psi_{i'}^{-1} (V_{d,X})\\
 V_{d,Y} &:= g_{p'}^{-1}\big( V_{\tau(i')\underline{\eta(p')}}\big)\\
 U_{d,Y} &:= \psi_{p'}^{-1}( V_{d,Y})
 \end{align*}
 Define a smooth map $s_d: V_{d,X}\times V_{d,Y} \ra \R^{m_{i'}}\oplus \R^{m_{p'}} \oplus \R^n$ by sending a point $(x,y)\in V_{d,X}\times V_{d,Y}$ to 
\[ s_d(x,y):= \langle s_{i'}(x), s_{p'}(y), g_{p'}(y) - f_{\tau(i')\eta(p')} h_{i'}(x)\rangle.\]
It is clear that we have
\[Z\cap (U_{d,X}\times U_{d,Y})=(\psi_{i'}\times\psi_{p'})^{-1}\big( s_d^{-1}(0)\big).\]
Thus the quadruple $(V_d, \underline{\R^{m_d}}, s_d, \psi_d)$, with
\begin{align*}
V_d&:= V_{d,X}\times V_{d,Y} \\
m_d&:= m_{i'}+m_{p'}+n\\
\psi_d &:= (\psi_{i'}\times \psi_{p'})|_{V_d},
\end{align*}
defines a Kuranishi chart of the open subset $U_d:=Z\cap (U_{i'}\times U_{p'})$ in $Z$. We take its germ $[V_d,\underline{\R^{m_d}},s_d,\psi_d]$ to form a Kuranishi atlas on $Z$.

\medskip
(1c.) Let $d_1,d_2\in D$ be such that $U_{d_1}\cap U_{d_2}\neq \emptyset$. Define $$[f_{d_1d_2},\hat{f}_{d_1d_2}]: [ V_{\underline{d_1}d_2}, \underline{\R^{m_{d_1}}}, s_{d_1}, \psi_{d_1} ] \ra [ V_{d_1\underline{d_2}}, \underline{\R^{m_{d_2}}}, s_{d_2}, \psi_{d_2}\big]$$ by setting
\begin{align*}
f_{d_1d_2}&:= f_{i_1'i_2'}\times f_{p_1'p_2'},\\
\hat{f}_{d_1d_2}&:= \begin{bmatrix}
\hat{f}_{i_1'i_2'} & 0 &0 \\
0 & \hat{f}_{p_1'p_2'} &0\\
\delta f_{\tau(i_2')\eta(p_2')}\Delta_{i_1'i_2'} &-\Delta_{p_1'p_2'} & \delta f_{\eta(p_1')\eta(p_2')}
\end{bmatrix}: \underline{\R^{m_{d_1}}} \ra \underline{\R^{m_{d_2}}},
\end{align*}
where the notation $\delta$ is as in Lemma~\ref{lem:taylor}, and the matrix operator is written in the decomposition $\underline{\R^{m_{d_1}}}=\underline{\R^{m_{i_1'}}}\oplus \underline{\R^{m_{p_1'}}}\oplus \underline{\R^{{n}}}$ viewed as a column vector, and similarly for $\underline{\R^{m_{d_2}}}$. We need to verify $(f_{d_1d_2},\hat{f}_{d_1d_2})$ is a Kuranishi morphism, i.e. the identity $\hat{f}_{d_1d_2} s_{d_1} = s_{d_2} f_{d_1d_2}$ holds. Indeed, we compute
\begin{align*}
&\hat{f}_{d_1d_2} s_{d_1}(x,y)\\
&= \begin{bmatrix}
\hat{f}_{i_1'i_2'} & 0 &0 \\
0 & \hat{f}_{p_1'p_2'} &0\\
\delta f_{\tau(i_2')\eta(p_2')}\Delta_{i_1'i_2'} &-\Delta_{p_1'p_2'} & \delta f_{\eta(p_1')\eta(p_2')}
\end{bmatrix} \begin{bmatrix} s_{i_1'}(x)\\ s_{p_1'}(y)\\ g_{p_1'}(y)-f_{\tau(i_1')\eta(p_1')}h_{i_1'}(x)\end{bmatrix}\\
&=\begin{bmatrix} s_{i_2'}(f_{i_1'i_2'}(x))\\
s_{p_2'}(f_{p_1'p_2'}(y))\\
\delta f_{\tau(i_2')\eta(p_2')}\Delta_{i_1'i_2'} s_{i_1'}(x)-\Delta_{p_1'p_2'}s_{p_1'}(y)+\delta f_{\eta(p_1')\eta(p_2')}\big(g_{p_1'}(y)-f_{\tau(i_1')\eta(p_1')}h_{i_1'}(x)\big)\end{bmatrix}
\end{align*}
The last entry can be simplified as
\begin{align*}
& \delta f_{\tau(i_2')\eta(p_2')}\Delta_{i_1'i_2'} s_{i_1'}(x)-\Delta_{p_1'p_2'}s_{p_1'}(y)+\delta f_{\eta(p_1')\eta(p_2')}\big(g_{p_1'}(y)-f_{\tau(i_1')\eta(p_1')}h_{i_1'}(x)\big) \\
&= \delta f_{\tau(i_2')\eta(p_2')}[f_{\tau(i_1')\tau(i_2')}h_{i_1'}(x)-h_{i_2'}f_{i_1'i_2'}(x)]-[f_{\eta(p_1')\eta(p_2')}g_{p_1'}(y)-g_{p_2'}f_{p_1'p_2'}(y)]+\\ &+[f_{\eta(p_1')\eta(p_2')}g_{p_1'}(y)-f_{\eta(p_1')\eta(p_2')}f_{\tau(i_1')\eta(p_1')}h_{i_1'}(x)]\\
&= [f_{\tau(i_1')\eta(p_2')}h_{i_1'}(x)-f_{\tau(i_2')\eta(p_2')}h_{i_2'}f_{i_1'i_2'}(x)]-[f_{\eta(p_1')\eta(p_2')}g_{p_1'}(y)-g_{p_2'}f_{p_1'p_2'}(y)]+\\
&+ [f_{\eta(p_1')\eta(p_2')}g_{p_1'}(y)-f_{\tau(i_1')\eta(p_2')}h_{i_1'}(x)]\\
&= g_{p_2'}f_{p_1'p_2'}(y)-f_{\tau(i_2')\eta(p_2')}h_{i_2'}f_{i_1'i_2'}(x)
\end{align*}
This shows that $(f_{d_1d_2},\hat{f}_{d_1d_2})$ is a Kuranishi morphism.

\medskip
(1d.) Let $d_1,d_2,d_3\in D$ be such that $U_{d_1}\cap U_{d_2}\cap U_{d_3}\neq \emptyset$. We define the homotopy data $\Lambda_{d_1d_2d_3}: \underline{\R^{m_{d_1}}} \ra \underline{\R^{n_{i_3}}}\oplus \underline{\R^{n_{p_3}}}$ by the matrix
\[ \Lambda_{d_1d_2d_3}:= \begin{bmatrix}
\Lambda_{i_1'i_2'i_3'} &0 &0\\
0 & \Lambda_{p_1'p_2'p_3'} &0
\end{bmatrix},\]
We verify that this indeed defines a homotopy. Since $f_{d_1d_2}$ is defined to be $f_{i_1'i_2'}\times f_{p_1'p_2'}$, it is clear that 
\[ f_{d_2d_3}f_{d_1d_2}-f_{d_1d_3}=\Lambda_{d_1d_2d_3}\circ s_{d_1}.\]
Furthermore, after restriction to $U_{d_1}\cap U_{d_2}\cap U_{d_3}$, we have
\begin{align*}
& \Lambda_{d_1d_2d_3}ds_{d_1}= \begin{bmatrix}
\Lambda_{i_1'i_2'i_3'} &0 &0\\
0 & \Lambda_{p_1'p_2'p_3'} &0
\end{bmatrix} \begin{bmatrix}
ds_{i_1'} &0\\ 0 & ds_{p_1'}\\-df_{\tau(i_1')\eta(p_1')}dh_{i_1'} & dg_{p_1'}
\end{bmatrix}\\
&=\begin{bmatrix}
\Lambda_{i_1'i_2'i_3'} ds_{i_1'} & 0\\
0 & \Lambda_{p_1'p_2'p_3'} ds_{p_1'}
\end{bmatrix}\\
&=\begin{bmatrix}
df_{i_2'i_3'}df_{i_1'i_2'}-df_{i_1'i_3'} & 0\\
0& df_{p_2'p_3'}df_{p_1'p_2'}-df_{p_1'p_3'}
\end{bmatrix}\\
&=df_{d_2d_3}df_{d_1d_2}-df_{d_1d_3}.
\end{align*}
And the other composition is
\begin{align}~\label{composition1}
\begin{split}
ds_{d_3}\Lambda_{d_1d_2d_3}&=\begin{bmatrix}
ds_{i_1'} &0\\ 0 & ds_{p_1'}\\-df_{\tau(i_1')\eta(p_1')}dh_{i_1'} & dg_{p_1'}
\end{bmatrix}\begin{bmatrix}
\Lambda_{i_1'i_2'i_3'} &0 &0\\
0 & \Lambda_{p_1'p_2'p_3'} &0
\end{bmatrix}\\
&=\begin{bmatrix}
ds_{i_3'}\Lambda_{i_1'i_2'i_3'} &0 &0\\
0& ds_{p_3'}\Lambda_{p_1'p_2'p_3'} &0\\
-df_{\tau(i_3')\eta(p_3')}dh_{i_3'}\Lambda_{i_1'i_2'i_3'} & dg_{p_3'}\Lambda_{p_1'p_2'p_3'} &0
\end{bmatrix}\end{split}
\end{align}
We need to show this agrees with
\begin{align}~\label{composition2}\begin{split}
&\hat{f}_{d_2d_3}\hat{f}_{d_1d_2}-\hat{f}_{d_1d_3}\\
=&\begin{bmatrix}
\hat{f}_{i_2'i_3'} & 0 &0 \\
0 & \hat{f}_{p_2'p_3'} &0\\
d f_{\tau(i_3')\eta(p_3')}\Delta_{i_2'i_3'} &-\Delta_{p_2'p_3'} & d f_{\eta(p_2')\eta(p_3')}
\end{bmatrix}\begin{bmatrix}
\hat{f}_{i_1'i_2'} & 0 &0 \\
0 & \hat{f}_{p_1'p_2'} &0\\
d f_{\tau(i_2')\eta(p_2')}\Delta_{i_1'i_2'} &-\Delta_{p_1'p_2'} & d f_{\eta(p_1')\eta(p_2')}
\end{bmatrix}\\
&-\begin{bmatrix}
\hat{f}_{i_1'i_3'} & 0 &0 \\
0 & \hat{f}_{p_1'p_3'} &0\\
d f_{\tau(i_3')\eta(p_3')}\Delta_{i_1'i_3'} &-\Delta_{p_1'p_3'} & d f_{\eta(p_1')\eta(p_3')}
\end{bmatrix}\\
=& \begin{bmatrix}
ds_{i_3'}\Lambda_{i_1'i_2'i_3'}&0&0\\
0&ds_{p_3'}\Lambda_{p_1'p_2'p_3'} &0\\
A & B & 0
\end{bmatrix}
\end{split}\end{align}
where the entries $A$ and $B$ are given by
\begin{align*}
A&= df_{\tau(i_3')\eta(p_3')}\Delta_{i_2'i_3'}\hat{f}_{i_1'i_2'}+df_{\tau(i_2')\eta(p_3')}\Delta_{i_1'i_2'}-df_{\tau(i_3')\eta(p_3')}\Delta_{i_1'i_3'}\\
&= df_{\tau(i_3')\eta(p_3')}[\Delta_{i_2'i_3'}\hat{f}_{i_1'i_2'}+df_{\tau(i_2')\tau(i_3')}\Delta_{i_1'i_2'}-\Delta_{i_1'i_3'}]\\
B&= -\Delta_{p_2'p_3'}\hat{f}_{p_1'p_2'}-df_{\eta(p_2')\eta(p_3')}\Delta_{p_1'p_2'}+\Delta_{p_1'p_3'}.
\end{align*}
Comparing the matrices in the last equation of~\ref{composition1} and~\ref{composition2}, we need to prove that $A=-df_{\tau(i_3')\eta(p_3')}dh_{i_3'}\Lambda_{i_1'i_2'i_3'}$ and $B=dg_{p_3'}\Lambda_{p_1'p_2'p_3'}$. Now, both identities follows from Equation~\ref{eq:2-hom}.

\medskip
(1e.) Finally we need to check the homotopies $\Lambda_{d_1d_2d_3}$ defined in $(4.)$ are compatible on quadruple intersections, i.e. to prove that
\[ \Lambda_{d_1d_3d_4}-\Lambda_{d_2d_3d_4}\circ \hat{f}_{d_1d_2} -\Lambda_{d_1d_2d_4}+df_{d_3d_4}\circ\Lambda_{d_1d_2d_3}=0\]
over a quadruple intersection $U_{d_1d_2d_3d_4}$. This follows directly from the corresponding equations satisfied by $\Lambda_{ijk}$'s and $\Lambda_{pqr}$'s. As a conclusion, we have shown the data $$\big\{ [V_d, \underline{\R^{m_d}}, s_d, \psi_d], [f_{d_1d_2}, \hat{f}_{d_1d_2}], [\Lambda_{d_1d_2d_3}]\big\}_{z,d_1,d_2,d_3\in D}$$ as defined above forms a Kuranishi atlas on the topological space $Z$. We denote this Kuranishi manifold by $\mathfrak{Z}$.

\medskip
(2.) Recall that $V_d= V_{d,X}\times V_{d,Y}\subset V_{i'}\times V_{p'}\subset V_i\times V_p$, and $\underline{\R^{m_d}}=\R^{m_{i}}\oplus \R^{m_{p}} \oplus \R^n$. Define the strict Kuranishi morphism  $\pi_1:\mathfrak{Z}\ra \mathfrak{X}$ by the natural projection onto the first component
\[ [(\pi_1)_d,(\widehat{\pi_1})_d]: [V_d, \underline{\R^{m_d}}, s_d, \psi_d] \ra [V_i, \underline{\R^{m_i}}, s_i, \psi_i],\]
with trivial $\Delta$'s. And similarly $\pi_2: \mathfrak{Z}\ra \mathfrak{Y}$ is given by the natural projection onto the second component, again with trivial $\Delta$'s.

\medskip
(3.) We define the a homotopy $\theta:(hs^{-1})(\pi_1)\ra (gt^{-1})(\pi_2)$ by
\[ \theta_d: \underline{\R^{m_d}}=\R^{m_{i'}}\oplus \R^{m_{p'}} \oplus \R^n \ra \R^n, \;\; (a,b,c)\mapsto -c.\]
One immediately checks that
\[ f_{\tau(i')\eta(p')}h_{i'}(x)-g_{p'}(y)=\theta_d\big( s_d(x,y)\big).\]
We also need to verify the Equation~\ref{eq:2-equiv}. Thus let $d_1, d_2$ be two indices. By construction in $(2.)$, we have $\Delta^{\pi_1}$ and $\Delta^{\pi_2}$ are both zero. Furthermore, since $M$ is a manifold, its $\Lambda$'s are zero. Thus in Equation~\ref{eq:2-equiv}, we are left to verify that
\[ \theta_{d_2}*[f_{d_1d_2},\hat{f}_{d_1d_2}]= [f_{\eta(p_1')\eta(p_2')},0]*\theta_{d_1}.\]
Evaluating the two operators on a vector $(a,b,c)\in\underline{\R^{m_{d_1}}}=\R^{m_{i_1'}}\oplus \R^{m_{p_1'}} \oplus \R^n$  by their definition, both yield $-df_{\eta(p_1')\eta(p_2')}(c)$.

\medskip
\begin{thm}~\label{thm:fiber}
With the constructions of $(1.), (2.)$ and $(3.)$ above, the following diagram is a $2$-fiber product diagram in the $2$-category $\mathfrak{Kur}$.
\[ \begin{xy}
(0,40)*{\mathfrak{Z}=\mathfrak{X}\times_M \mathfrak{Y}}="F";
(40,40)*{\mathfrak{Y}=(Y,\BB)}="A";
(40,0)*{(M,\UU)}="C";
(0,0)*{\mathfrak{X}=(X,\AA)}="D";
(15,15)*{\;}="E";
(25,25)*{\;}="G";
{\ar@{=>}_{[\theta]} "E"; "G"};
{\ar@{->}^{gt^{-1}} "A"; "C"};
{\ar@{->}^{hr^{-1}} "D"; "C"};
{\ar@{->}^{\pi_1} "F"; "D"};
{\ar@{->}^{\pi_2} "F"; "A"};
\end{xy}\] 
\end{thm}

\Pf. Assume that we are given another diagram
\[ \begin{xy}
(0,40)*{\mathfrak{W}}="F";
(40,40)*{\mathfrak{Y}=(Y,\BB)}="A";
(40,0)*{(M,\UU)}="C";
(0,0)*{\mathfrak{X}=(X,\AA)}="D";
(15,15)*{\;}="E";
(25,25)*{\;}="G";
{\ar@{=>}_{[\chi]} "E"; "G"};
{\ar@{->}^{gt^{-1}} "A"; "C"};
{\ar@{->}^{hr^{-1}} "D"; "C"};
{\ar@{->}^{k_1a^{-1}} "F"; "D"};
{\ar@{->}^{k_2b^{-1}} "F"; "A"};
\end{xy}\] 
in the category $\mathfrak{Kur}$. We use the notations
\[ \mathfrak{W}=\big\{ [V_\alpha,\underline{\R^{m_\alpha}}, s_\alpha, \psi_\alpha], [f_{\alpha\beta},\hat{f}_{\alpha\beta}], [\Lambda_{\alpha\beta\gamma}]\big\}_{\alpha,\beta,\gamma\in Q}.\]
Recall by definition, the $2$-morphism $[\chi]$ is represented by a homotopy $\chi$ in the following diagram
\[ \begin{xy}
(0,60)*{\mathfrak{W}}="F";
(60,60)*{\mathfrak{Y}=(Y,\BB)}="A";
(60,0)*{(M,\UU)}="C";
(0,0)*{\mathfrak{X}=(X,\AA)}="D";
(20,0)*{\mathfrak{X}'}="B";
(60,40)*{\mathfrak{Y}'}="E";
(40,60)*{\mathfrak{W}_2}="G";
(0,20)*{\mathfrak{W}_1}="H";
(20,40)*{\mathfrak{W}'}="I";
(40,40)*{\mathfrak{W}_4}="J";
(20,20)*{\mathfrak{W}_3}="K";
(35,15)*{\;}="L";
(45,25)*{\;}="M";
{\ar@{=>}_{[\chi]} "L"; "M"};
{\ar@{->}^{t} "E"; "A"};
{\ar@{->}^{g} "E"; "C"};
{\ar@{->}^{h} "B"; "C"};
{\ar@{->}^{r} "B"; "D"};
{\ar@{->}^{k_1} "H"; "D"};
{\ar@{->}^{k_2} "G"; "A"};
{\ar@{->}^{a} "H"; "F"};
{\ar@{->}^{b} "G"; "F"};
{\ar@{->}@[red] "K"; "B"};
{\ar@{->}@[red] "K"; "H"};
{\ar@{->}@[red] "J"; "E"};
{\ar@{->}@[red] "J"; "G"};
{\ar@{->}@[red] "I"; "J"};
{\ar@{->}@[red] "I"; "K"};
{\ar@{->}@[red]@/^1.5pc/^{k_1'} "I"; "B"};
{\ar@{->}@[red]@/_1.5pc/_{k_2'} "I"; "E"};
\end{xy}\] 
Denote charts of the Kuranishi manifold $\mathfrak{W}'$ by
\[ \mathfrak{W}'= \big\{ [V_{\alpha'},\underline{\R^{m_{\alpha'}}}, s_{\alpha'}, \psi_{\alpha'}]\big\}_{\alpha'\in Q'}.\]
As shown in the above diagram, we also have strict Kuranishi morphisms 
\[ k_1': \mathfrak{W}' \ra \mathfrak{X}', \mbox{\;\; and \;\;} k_2': \mathfrak{W}' \ra \mathfrak{Y}'.\]
Let $\alpha'\in Q'$ be a index. Assume that $\alpha'\mapsto i'$ and $\alpha'\mapsto p'$ under the maps $k_1'$ and $k_2$'.
Define a morphism $u: \mathfrak{W}' \ra \mathfrak{Z}$ by
\begin{itemize}
\item[--] On indices, we set $\alpha' \mapsto d:=(i', p')\in D$.
\item[--] The morphism $[u_{\alpha'},\hat{u}_{\alpha'}]: [V_{\alpha'},\underline{\R^{m_{\alpha'}}}, s_{\alpha'}, \psi_{\alpha'}]\ra [V_d, \underline{\R^{m_d}}, s_d, \psi_d]$ is defined by
\begin{align*}
u_{\alpha'}: & w\in V_{\alpha'} \mapsto \big((k_1')_{\alpha'}(w), (k_2')_{\alpha'} (w)\big) \in V_d;\\
\hat{u}_{\alpha'}: & s \in \underline{\R^{m_{\alpha'}}} \mapsto \Big(\widehat{(k_1')}_{\alpha'}(s),\widehat{(k_2')}_{\alpha'}(s), -\chi_{\alpha'}(s)\Big)\in \underline{\R^{m_d}}
\end{align*}
\item[--] The homotopies on transition is defined by
\[ \Delta^u_{\alpha'}:  s\in \underline{\R^{m_{\alpha'}}} \mapsto \Big(\Delta_{\alpha'}^{k_1'}(s), \Delta_{\alpha'}^{k_2'}(s)\Big)\in \underline{\R^{i'}}\oplus\underline{\R^{p'}}.\]
\end{itemize}
The morphism $u$ fits into the following diagram.
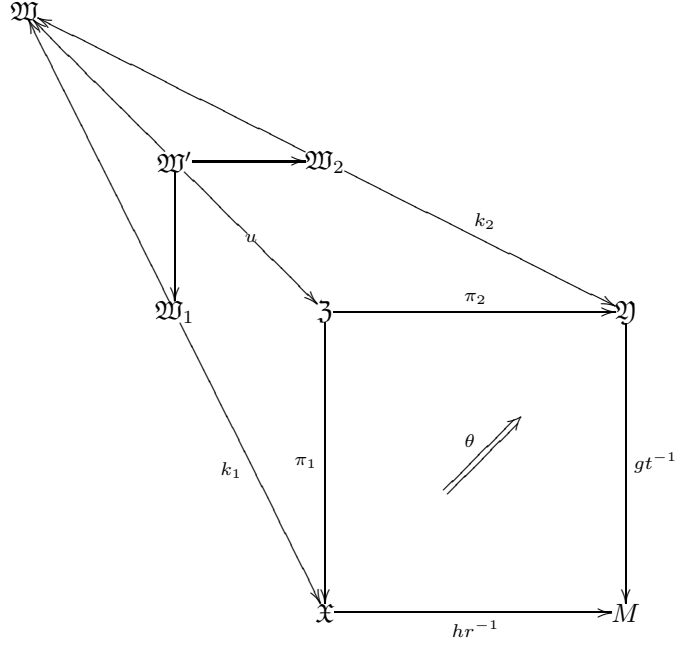
\[\begin{xy}
(-40,80)*{\mathfrak{W}}="E";
(0,0)*{\mathfrak{X}}="B";
(40,0)*{M}="A";
(40,40)*{\mathfrak{Y}}="C";
(0,40)*{\mathfrak{Z}}="D";
(16,16)*{\;}="H";
(26,26)*{\;}="G";
(-20,60)*{\mathfrak{W}'}="F";
(-20,40)*{\mathfrak{W}_1}="I";
(0,60)*{\mathfrak{W}_2}="J";
{\ar@{=>}^{\theta} "H"; "G"};
{\ar@{->}_{hr^{-1}} "B"; "A"};
{\ar@{->}^{gt^{-1}} "C"; "A"};
{\ar@{->}_{\pi_1} "D"; "B"};
{\ar@{->}^{\pi_2} "D"; "C"};
{\ar@{->}@[red]|-{u} "F"; "D"};
{\ar@{->} "F"; "E"};
{\ar@{->}_{k_1} "I"; "B"};
{\ar@{->}^{k_2} "J"; "C"};
{\ar@{->} "I"; "E"};
{\ar@{->} "J"; "E"};
{\ar@{->} "F"; "I"};
{\ar@{->} "F"; "J"};
\end{xy}\]
The two trapezoidal regions are strictly commutative diagrams. Hence to prove the identity~\ref{eq:2-fiber-1}, it suffices to show that
\[ \theta\circ_1 u = \chi.\]
Indeed, for an index $\alpha'$, and $s\in \underline{\R^{m_{\alpha'}}}$, we have
\begin{align*}
&[\theta\circ_1 u]_\alpha (s)\\
=& \theta_d \hat{u}_{\alpha'} (s)\\
=& \theta_d \Big(\widehat{(k_1')}_{\alpha'}(s),\widehat{(k_2')}_{\alpha'}(s), -\chi_{\alpha'}(s)\Big)\\
=& \chi_{\alpha'}(s).
\end{align*}
Since by definition, $\theta_d$ is the negative of the third projection map. 

Next, we prove part $(B.)$ of Definition~\ref{defi:2-fiber}. By taking possibly further refinement of $\mathfrak{W}'$, we can start with the following diagram, 
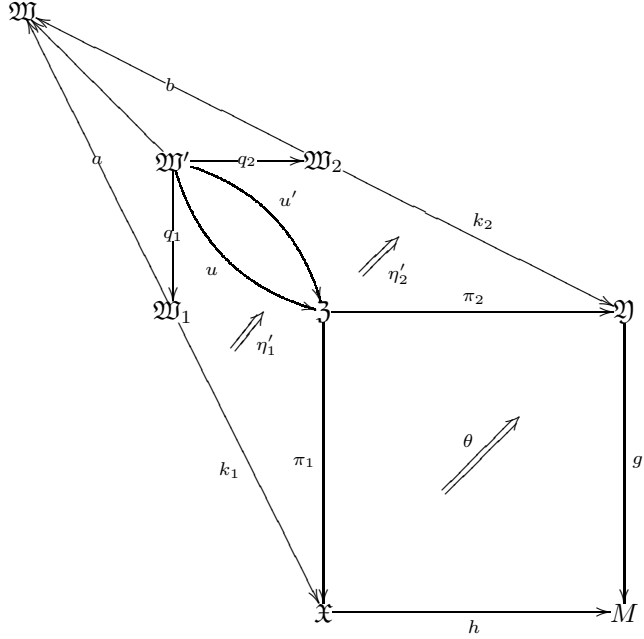
\[\begin{xy}
(-40,80)*{\mathfrak{W}}="E";
(0,0)*{\mathfrak{X}}="B";
(40,0)*{M}="A";
(40,40)*{\mathfrak{Y}}="C";
(0,40)*{\mathfrak{Z}}="D";
(16,16)*{\;}="H";
(26,26)*{\;}="G";
(-20,60)*{\mathfrak{W}'}="F";
(-20,40)*{\mathfrak{W}_1}="I";
(0,60)*{\mathfrak{W}_2}="J";
(-12,35)*{\;}="K";
(-8,40)*{\;}="L";
(5,45)*{\;}="M";
(10,50)*{\;}="N";
{\ar@{=>}_{\eta_2'} "M"; "N"};
{\ar@{=>}_{\eta_1'} "K"; "L"};
{\ar@{=>}^{\theta} "H"; "G"};
{\ar@{->}_{h} "B"; "A"};
{\ar@{->}^{g} "C"; "A"};
{\ar@{->}_{\pi_1} "D"; "B"};
{\ar@{->}^{\pi_2} "D"; "C"};
{\ar@{->}@[red]@/_1pc/_{u} "F"; "D"};
{\ar@{->}@[red]@/^1pc/^{u'} "F"; "D"};
{\ar@{->} "F"; "E"};
{\ar@{->}_{k_1} "I"; "B"};
{\ar@{->}^{k_2} "J"; "C"};
{\ar@{->}|-{a} "I"; "E"};
{\ar@{->}|-{b} "J"; "E"};
{\ar@{->}|-{q_1} "F"; "I"};
{\ar@{->}|-{q_2} "F"; "J"};
\end{xy}\]
Here we are given $2$-morphisms
\[ \eta_1': k_1 q_1 \ra \pi_1 u', \;\; \mbox{ and }\;\; \eta_2': \pi_2 u' \ra k_2 q_2.\]
We define a $2$-morphism $\lambda: u \ra u'$ by formula
\begin{align*}
\lambda_{\alpha'}: & \underline{\R^{m_{\alpha'}}} \ra \underline{\R^{n_{i'}}}\oplus\underline{\R^{n_{p'}}}\\
& s \mapsto \Big( (\eta'_1)_{\alpha'} (s), -(\eta'_2)_{\alpha'} (s)\Big).
\end{align*}
One easily checks that
\[ \eta_1'= \pi_1\circ_1 \lambda, \;\; \mbox{ and }\;\; 0=\eta_2'\circ_0 (\pi_2\circ_1\lambda).\]
Furthermore, the fact that $\pi_1$ and $\pi_2$ are projections forces the uniqueness of $\lambda$.\ed

\appendix

\section{Proof of Theorem~\ref{thm:pre}}~\label{app}

In this appendix, for simplicity, we shall drop the germ notation $[-]$, and write a Kuranishi map $[h,\hat{h}]$ between charts simply by $h$. 

\medskip
\begin{prop}
The vertical composition $\circ_0$ defined by Equation~\ref{eq:2-hom} is associative.
\end{prop}
\Pf. Let $\Upsilon^{12}: h^1\ra h^2,\Upsilon^{23}:h^2\ra h^3,\Upsilon^{34}: h^3\ra h^4$ be $2$-morphisms from $\mathfrak{X}$ to $\mathfrak{Y}$. We calculate the two ways to compose them, using Equation~\ref{eq:2-hom}.
\begin{align*}
&\big(\Upsilon^{34}\circ_0  (\Upsilon^{23}\circ_0\Upsilon^{12})\big)_i= \Upsilon^{34}_i-\Lambda_{\tau_1(i)\tau_3(i)\tau_4(i)}*h^1_i+f_{\tau_3(i)\tau_4(i)}*(\Upsilon^{23}\circ_0\Upsilon^{12})_i\\
&=\Upsilon^{34}_i-\Lambda_{\tau_1(i)\tau_3(i)\tau_4(i)}*h^1_i+\\
&f_{\tau_3(i)\tau_4(i)}*\Upsilon^{23}_i-f_{\tau_3(i)\tau_4(i)}*\Lambda_{\tau_1(i)\tau_2(i)\tau_3(i)}*h^1_i+f_{\tau_3(i)\tau_4(i)}*f_{\tau_2(i)\tau(3)}*\Upsilon^{12}_i\\
& \\
&\big((\Upsilon^{34}\circ_0 \Upsilon^{23})\circ_0\Upsilon^{12}\big)_i=
(\Upsilon^{34}\circ_0  \Upsilon^{23})_i -\Lambda_{\tau_1(i)\tau_2(i)\tau_4(i)}*h^1_i+f_{\tau_2(i)\tau_4(i)}*\Upsilon^{12}_i\\
&= \Upsilon^{34}_i-\Lambda_{\tau_2(i)\tau_3(i)\tau_4(i)}*h^2_i+f_{\tau_3(i)\tau_4(i)}*\Upsilon^{23}_i-\Lambda_{\tau_1(i)\tau_2(i)\tau_4(i)}*h^1_i+f_{\tau_2(i)\tau_4(i)}*\Upsilon^{12}_i\\
\end{align*}
We use the identity
\[ f_{\tau_3(i)\tau_4(i)}*f_{\tau_2(i)\tau(3)}*\Upsilon^{12}_i-f_{\tau_2(i)\tau_4(i)}*\Upsilon^{12}_i
= \Lambda_{\tau_2(i)\tau_3(i)\tau_4(i)}*\big( f_{\tau_1(i)\tau_2(i)}*h^1_i-h^2_i\big)
\]
to compute the difference of the above two compositions:
\begin{align*}
&\big(\Upsilon^{34}\circ_0(\Upsilon^{23}\circ_0\Upsilon^{12})\big)_i-\big((\Upsilon^{34}\circ_0 \Upsilon^{23})\circ_0\Upsilon^{12}\big)_i\\
&=\big(-\Lambda_{\tau_1(i)\tau_3(i)\tau_4(i)}-f_{\tau_3(i)\tau_4(i)}*\Lambda_{\tau_1(i)\tau_2(i)\tau_3(i)}+\\
&\Lambda_{\tau_2(i)\tau_3(i)\tau_4(i)}
*f_{\tau_1(i)\tau_2(i)}+\Lambda_{\tau_1(i)\tau_2(i)\tau_4(i)}\big)*h^1_i
\end{align*}
The right hand side vanishes by Equation~\ref{eq:cocycle2}.\ed

\medskip
\begin{prop}
The horizontal composition $\circ_1$ defined by Equation~\ref{eq:horizontal} is associative.
\end{prop}

\Pf. Let 
\begin{align*}
\Upsilon & : h^1\ra h^2: \mathfrak{X}\ra \mathfrak{Y},\\
\Gamma &: g^1\ra g^2: \mathfrak{Y}\ra \mathfrak{Z},\\
\Sigma &: k^1\ra k^2: \mathfrak{Z}\ra \mathfrak{W}
\end{align*}
be three horizontally composable $2$-morphisms. We first compute
\begin{align*}
&\big(\Sigma\circ_1(\Gamma\circ_1\Upsilon)\big)_i\\
=&k^2_{\sigma_2\tau_2(i)}*(\Gamma\circ_1\Upsilon)_i+f_{\eta_2\sigma_1\tau_1(i),\eta_2\sigma_2\tau_2(i)}*\Sigma_{\sigma_1\tau_1(i)}*g^1_{\tau_1(i)}*h^1_i+\\
&\Delta^{k^2}_{\sigma_1\tau_1(i)\sigma_2\tau_2(i)}**g^1_{\tau_1(i)}*h^1_i-\Lambda_{\eta_1\sigma_1\tau_1(i)\eta_2
\sigma_1\tau_1(i)\eta_2\sigma_2\tau_2(i)}*k^1_{\sigma_1\tau_1(i)}*g^1_{\tau_1(i)}*h^1_i\\
=&k^2_{\sigma_2\tau_2(i)}*g^2_{\tau_2(i)}*\Upsilon_i+k^2_{\sigma_2\tau_2(i)}*f_{\sigma_2\tau_1(i)\sigma_2\tau_2(i)}
*\Gamma_{\tau_1(i)}*h_i^1\\
&+k^2_{\sigma_2\tau_2(i)}*\Delta^{g^2}_{\tau_1(i)\tau_2(i)}*h^1_i-k^2_{\sigma_2\tau_2(i)}*\Lambda_{\sigma_1\tau_1(i)\sigma_2\tau_1(i)\sigma_2\tau_2(i)}*g^1_{\tau_1(i)}*h^1_i\\
&+f_{\eta_2\sigma_1\tau_1(i),\eta_2\sigma_2\tau_2(i)}*\Sigma_{\sigma_1\tau_1(i)}*g^1_{\tau_1(i)}*h^1_i+
\Delta^{k^2}_{\sigma_1\tau_1(i)\sigma_2\tau_2(i)}*g^1_{\tau_1(i)}*h^1_i\\
&-\Lambda_{\eta_1\sigma_1\tau_1(i)\eta_2
\sigma_1\tau_1(i)\eta_2\sigma_2\tau_2(i)}*k^1_{\sigma_1\tau_1(i)}*g^1_{\tau_1(i)}*h^1_i\\
\end{align*}
The other way of composition is
\begin{align*}
&\big((\Sigma\circ_1\Gamma)\circ_1\Upsilon\big)_i\\
=& k^2_{\sigma_2\tau_2(i)}*g^2_{\tau_2(i)}*\Upsilon_i +f_{\eta_2\sigma_2\tau_1(i)\eta_2\sigma_2\tau_2(i)}
*(\Sigma\circ_1\Gamma)_{\tau_1(i)}*h^1_i\\
&+\Delta^{k^2g^2}_{\tau_1(i)\tau_2(i)}*h^1_i-\Lambda_{\eta_1\sigma_1\tau_1(i)\eta_2\sigma_2\tau_1(i)\eta_2\sigma_2\tau_2(i)}
*k^1_{\sigma_1\tau_1(i)}*g^1_{\tau_1(i)}*h^1_i\\
=& k^2_{\sigma_2\tau_2(i)}*g^2_{\tau_2(i)}*\Upsilon_i +f_{\eta_2\sigma_2\tau_1(i)\eta_2\sigma_2\tau_2(i)}
* k^2_{\sigma_2\tau_1(i)}*\Gamma_{\tau_1(i)}*h^1_i\\
& +f_{\eta_2\sigma_2\tau_1(i)\eta_2\sigma_2\tau_2(i)}*f_{\eta_2\sigma_1\tau_1(i)\eta_2\sigma_2\tau_1(i)}*\Sigma_{\sigma_1\tau_1(i)}*g^1_{\tau_1(i)}*h^1_i
\\ &+f_{\eta_2\sigma_2\tau_1(i)\eta_2\sigma_2\tau_2(i)}*\Delta^{k^2}_{\sigma_1\tau_1(i)\sigma_2\tau_1(i)}*g^1_{\tau_1(i)}*h^1_i\\
&-f_{\eta_2\sigma_2\tau_1(i)\eta_2\sigma_2\tau_2(i)}\Lambda_{\eta_1\sigma_1\tau_1(i)\eta_2\sigma_1\tau_1(i)\eta_2\sigma_2\tau_1(i)}*k^1_{\sigma_1\tau_1(i)}*g^1_{\tau_1(i)}*h^1_i\\
&+\Delta^{k^2g^2}_{\tau_1(i)\tau_2(i)}*h^1_i-\Lambda_{\eta_1\sigma_1\tau_1(i)\eta_2\sigma_2\tau_1(i)\eta_2\sigma_2\tau_2(i)}
*k^1_{\sigma_1\tau_1(i)}*g^1_{\tau_1(i)}*h^1_i
\end{align*}
We group the difference $\big((\Sigma\circ_1\Gamma)\circ_1\Upsilon\big)_i-\big(\Sigma\circ_1(\Gamma\circ_1\Upsilon)\big)_i$ by
\begin{align*}
&\big((\Sigma\circ_1\Gamma)\circ_1\Upsilon\big)_i-\big(\Sigma\circ_1(\Gamma\circ_1\Upsilon)\big)_i=A+B+C\\
A=&\big(\Delta^{k^2g^2}_{\tau_1(i)\tau_2(i)}*h^1_i-k^2_{\sigma_2\tau_2(i)}*\Delta^{g^2}_{\tau_1(i)\tau_2(i)}*h^1_i\big)\\
&-\Delta^{k^2}_{\sigma_2\tau_1(i)\sigma_2\tau_2(i)}*f_{\sigma_1\tau_1(i),\sigma_2\tau_1(i)}*g^1_{\tau_1(i)}*h^1_i\\
&-\big( k^2_{\sigma_2\tau_2(i)}*f_{\sigma_2\tau_1(i)\sigma_2\tau_2(i)}-f_{\eta_2\sigma_2\tau_1(i)\eta_2\sigma_2\tau_2(i)}
* k^2_{\sigma_2\tau_1(i)}\big)
*\Gamma_{\tau_1(i)}*h_i^1\\
B=&f_{\eta_2\sigma_2\tau_1(i)\eta_2\sigma_2\tau_2(i)}*\Delta^{k^2}_{\sigma_1\tau_1(i)\sigma_2\tau_1(i)}*g^1_{\tau_1(i)}*h^1_i
-\Delta^{k^2}_{\sigma_1\tau_1(i)\sigma_2\tau_2(i)}*g^1_{\tau_1(i)}*h^1_i\\
&+k^2_{\sigma_2\tau_2(i)}*\Lambda_{\sigma_1\tau_1(i)\sigma_2\tau_1(i)\sigma_2\tau_2(i)}*g^1_{\tau_1(i)}*h^1_i
\\&+\Delta^{k^2}_{\sigma_2\tau_1(i)\sigma_2\tau_2(i)}*f_{\sigma_1\tau_1(i),\sigma_2\tau_1(i)}*g^1_{\tau_1(i)}*h^1_i\\
C=& f_{\eta_2\sigma_2\tau_1(i)\eta_2\sigma_2\tau_2(i)}*f_{\eta_2\sigma_1\tau_1(i)\eta_2\sigma_2\tau_1(i)}*\Sigma_{\sigma_1\tau_1(i)}*g^1_{\tau_1(i)}*h^1_i\\
&-f_{\eta_2\sigma_2\tau_1(i)\eta_2\sigma_2\tau_2(i)}\Lambda_{\eta_1\sigma_1\tau_1(i)\eta_2\sigma_1\tau_1(i)\eta_2\sigma_2\tau_1(i)}*k^1_{\sigma_1\tau_1(i)}*g^1_{\tau_1(i)}*h^1_i\\
&-\Lambda_{\eta_1\sigma_1\tau_1(i)\eta_2\sigma_2\tau_1(i)\eta_2\sigma_2\tau_2(i)}
*k^1_{\sigma_1\tau_1(i)}*g^1_{\tau_1(i)}*h^1_i\\
&-f_{\eta_2\sigma_1\tau_1(i),\eta_2\sigma_2\tau_2(i)}*\Sigma_{\sigma_1\tau_1(i)}*g^1_{\tau_1(i)}*h^1_i\\&+\Lambda_{\eta_1\sigma_1\tau_1(i)\eta_2
\sigma_1\tau_1(i)\eta_2\sigma_2\tau_2(i)}*k^1_{\sigma_1\tau_1(i)}*g^1_{\tau_1(i)}*h^1_i
\end{align*}
where observe that in the first term $A$ we subtracted an extra term $\Delta^{k^2}_{\sigma_2\tau_1(i)\sigma_2\tau_2(i)}*f_{\sigma_1\tau_1(i),\sigma_2\tau_1(i)}*g^1_{\tau_1(i)}*h^1_i$, and added it back in the term $B$. We first prove that $A=0$. For this, observe that
\[ \Delta^{k^2g^2}_{\tau_1(i)\tau_2(i)}*h^1_i-k^2_{\sigma_2\tau_2(i)}*\Delta^{g^2}_{\tau_1(i)\tau_2(i)}*h^1_i=\Delta^{k^2}_{\sigma_2\tau_1(i)\sigma_2\tau_2(i)}*g^2_{\tau_1(i)}*h^1_i.\]
Using $\Gamma$ is a $2$-morphisms, we obtain
\begin{align*}
&\Delta^{k^2}_{\sigma_2\tau_1(i)\sigma_2\tau_2(i)}*g^2_{\tau_1(i)}*h^1_i-\Delta^{k^2}_{\sigma_2\tau_1(i)\sigma_2\tau_2(i)}*f_{\sigma_1\tau_1(i),\sigma_2\tau_1(i)}*g^1_{\tau_1(i)}*h^1_i\\
=& \big( k^2_{\sigma_2\tau_2(i)}*f_{\sigma_2\tau_1(i)\sigma_2\tau_2(i)}-f_{\eta_2\sigma_2\tau_1(i)\eta_2\sigma_2\tau_2(i)}
* k^2_{\sigma_2\tau_1(i)}\big)
*\Gamma_{\tau_1(i)}*h_i^1
\end{align*}
This proves that $A=0$. We apply Equation~\ref{eq:2-hom} to the morphism $k^2$ to get
\[ B=\Lambda_{\eta_2\sigma_1\tau_1(i)\eta_2\sigma_2\tau_1(i)\eta_2\sigma_2\tau_2(i)}*k^2_{\sigma_1\tau_1(i)}*g^1_{\tau_1(i)}*h^1_i.\]
In part $C$, we have
\begin{align*}
&\big(f_{\eta_2\sigma_2\tau_1(i)\eta_2\sigma_2\tau_2(i)}*f_{\eta_2\sigma_1\tau_1(i)\eta_2\sigma_2\tau_1(i)}-
f_{\eta_2\sigma_1\tau_1(i),\eta_2\sigma_2\tau_2(i)}\big)*\Sigma_{\sigma_1\tau_1(i)}*g^1_{\tau_1(i)}*h^1_i\\
=& \Lambda_{\eta_2\sigma_1\tau_1(i)\eta_2\sigma_2\tau_1(i)\a_2\sigma_2\tau_2(i)}
*\big(f_{\eta_1\sigma_1\tau_1(i)\eta_2\sigma_1\tau_1(i)}*k^1_{\sigma_1\tau_1(i)}-k^2_{\sigma_1\tau_1(i)}\big)*g^1_{\tau_1(i)}*h^1_i
\end{align*}
Thus, combining with part $B$, we have
\begin{align*}
&B+C\\
=&\Lambda_{\eta_2\sigma_1\tau_1(i)\eta_2\sigma_2\tau_1(i)\a_2\sigma_2\tau_2(i)}
*f_{\eta_1\sigma_1\tau_1(i)\eta_2\sigma_1\tau_1(i)}*k^1_{\sigma_1\tau_1(i)}*g^1_{\tau_1(i)}*h^1_i\\
&-f_{\eta_2\sigma_2\tau_1(i)\eta_2\sigma_2\tau_2(i)}\Lambda_{\eta_1\sigma_1\tau_1(i)\eta_2\sigma_1\tau_1(i)\eta_2\sigma_2\tau_1(i)}*k^1_{\sigma_1\tau_1(i)}*g^1_{\tau_1(i)}*h^1_i\\
&-\Lambda_{\eta_1\sigma_1\tau_1(i)\eta_2\sigma_2\tau_1(i)\eta_2\sigma_2\tau_2(i)}
*k^1_{\sigma_1\tau_1(i)}*g^1_{\tau_1(i)}*h^1_i\\
&+\Lambda_{\eta_1\sigma_1\tau_1(i)\eta_2
\sigma_1\tau_1(i)\eta_2\sigma_2\tau_2(i)}*k^1_{\sigma_1\tau_1(i)}*g^1_{\tau_1(i)}*h^1_i 
\\=&0
\end{align*}
The last identity is by Equation~\ref{eq:cocycle2}.\ed

\medskip
\begin{prop}
We have the interchange law identity
\[ \big(\Gamma^{23}\circ_0\Gamma^{12}\big)\circ_1\big(\Upsilon^{23}\circ_0
\Upsilon^{12}\big)=\big(\Gamma^{23}\circ_1\Upsilon^{23}\big)\circ_0
\big(\Gamma^{12}\circ_1\Upsilon^{12}\big).\]
\end{prop}

\Pf. Computing both sides yields
\begin{align*}
&[\big(\Gamma^{23}\circ_0\Gamma^{12}\big)\circ_1\big(\Upsilon^{23}\circ_0
\Upsilon^{12}\big)]_i\\
=&g^3_{\tau_3(i)}*\big(\Upsilon^{23}\circ_0
\Upsilon^{12}\big)_i+f_{\sigma_3\tau_1(i)\sigma_3\tau_3(i)}*\big(\Gamma^{23}\circ_0\Gamma^{12}\big)_{\tau_1(i)}*h^1_i\\
&+\Delta^{g^3}_{\tau_1(i)\tau_3(i)}*h^1_i-\Lambda_{\sigma_1\tau_1(i)\sigma_3\tau_1(i)\sigma_3\tau_3(i)}*g^1_{\tau_1(i)}*h^1_i\\
=&g^3_{\tau_3(i)}*\Upsilon^{23}_i-g^3_{\tau_3(i)}*\Lambda_{\tau_1(i)\tau_2(i)\tau_3(i)}*h^1_i+g^3_{\tau_3(i)}*f_{\tau_2(i)\tau_3(i)}*\Upsilon^{12}_i\\
&+f_{\sigma_3\tau_1(i)\sigma_3\tau_3(i)}*\Gamma^{23}_{\tau_1(i)}*h^1_i-f_{\sigma_3\tau_1(i)\sigma_3\tau_3(i)}*\Lambda_{\sigma_1\tau_1(i)\sigma_2\tau_1(i)\sigma_3\tau_1(i)}*g^1_{\tau_1(i)}*h^1_i
\\&+f_{\sigma_3\tau_1(i)\sigma_3\tau_3(i)}*f_{\sigma_2\tau_1(i)\sigma_3\tau_1(i)}*\Gamma^{12}_{\tau_1(i)}*h^1_i
\\&+\Delta^{g^3}_{\tau_1(i)\tau_3(i)}*h^1_i-\Lambda_{\sigma_1\tau_1(i)\sigma_3\tau_1(i)\sigma_3\tau_3(i)}*g^1_{\tau_1(i)}*h^1_i
\end{align*}
\begin{align*}
&[\big(\Gamma^{23}\circ_1\Upsilon^{23}\big)\circ_0\big(\Gamma^{12}\circ_1\Upsilon^{12}\big)]_i\\
=&\big(\Gamma^{23}\circ_1\Upsilon^{23}\big)_i-\Lambda_{\sigma_1\tau_1(i)\sigma_2\tau_2(i)\sigma_3\tau_3(i)}*g^1_{\tau_1(i)}*h^1_i
+f_{\sigma_2\tau_2(i)\sigma_3\tau_3(i)}*\big(\Gamma^{12}\circ_1\Upsilon^{12}\big)_i\\
=&g^3_{\tau_3(i)}*\Upsilon^{23}_i + f_{\sigma_3\tau_2(i)\sigma_3\tau_3(i)}*\Gamma^{23}_{\tau_2(i)}*h^2_i+\Delta^{g^3}_{\tau_2(i)\tau_3(i)}*h^2_i
\\&-\Lambda_{\sigma_2\tau_2(i)\sigma_3\tau_2(i)\sigma_3\tau_3(i)}*g^2_{\tau_2(i)}*h^2_i\\
&-\Lambda_{\sigma_1\tau_1(i)\sigma_2\tau_2(i)\sigma_3\tau_3(i)}*g^1_{\tau_1(i)}*h^1_i+f_{\sigma_2\tau_2(i)\sigma_3\tau_3(i)}*
g^2_{\tau_2(i)}*\Upsilon^{12}_i\\
&+f_{\sigma_2\tau_2(i)\sigma_3\tau_3(i)}*f_{\sigma_2\tau_1(i)\sigma_2\tau_2(i)}*\Gamma^{12}_{\tau_1(i)}*h^1_i+f_{\sigma_2\tau_2(i)\sigma_3\tau_3(i)}*
\Delta^{g^2}_{\tau_1(i)\tau_2(i)}*h^1_i\\
&-f_{\sigma_2\tau_2(i)\sigma_3\tau_3(i)}*\Lambda_{\sigma_1\tau_1(i)\sigma_2\tau_1(i)\sigma_2\tau_2(i)}*g^1_{\tau_1(i)}*h^1_i\\
\end{align*}
We regroup terms in the difference of the two by
\begin{align*}
&[\big(\Gamma^{23}\circ_1\Upsilon^{23}\big)\circ_0\big(\Gamma^{12}\circ_1\Upsilon^{12}\big)]_i
-[\big(\Gamma^{23}\circ_0\Gamma^{12}\big)\circ_1\big(\Upsilon^{23}\circ_0
\Upsilon^{12}\big)]_i=A+B+C\\
A=&f_{\sigma_3\tau_2(i)\sigma_3\tau_3(i)}*\Gamma^{23}_{\tau_2(i)}*h^2_i+\Delta^{g^3}_{\tau_2(i)\tau_3(i)}*h^2_i\\
&-\Lambda_{\sigma_2\tau_2(i)\sigma_3\tau_2(i)\sigma_3\tau_3(i)}*g^2_{\tau_2(i)}*h^2_i
+f_{\sigma_2\tau_2(i)\sigma_3\tau_3(i)}*g^2_{\tau_2(i)}*\Upsilon^{12}_i\\
&+g^3_{\tau_3(i)}*\Lambda_{\tau_1(i)\tau_2(i)\tau_3(i)}*h^1_i-g^3_{\tau_3(i)}*f_{\tau_2(i)\tau_3(i)}*\Upsilon^{12}_i-\Delta^{g^3}_{\tau_1(i)\tau_3(i)}*h^1_i\\
B=&B_1+B_2\\
B_1=&-f_{\sigma_2\tau_2(i)\sigma_3\tau_3(i)}*\Lambda_{\sigma_1\tau_1(i)\sigma_2\tau_1(i)\sigma_2\tau_2(i)}*g^1_{\tau_1(i)}*h^1_i
-\\&\Lambda_{\sigma_1\tau_1(i)\sigma_2\tau_2(i)\sigma_3\tau_3(i)}*g^1_{\tau_1(i)}*h^1_i\\
&+\Lambda_{\sigma_1\tau_1(i)\sigma_3\tau_1(i)\sigma_3\tau_3(i)}*g^1_{\tau_1(i)}*h^1_i
\\&+f_{\sigma_3\tau_1(i)\sigma_3\tau_3(i)}*\Lambda_{\sigma_1\tau_1(i)\sigma_2\tau_1(i)\sigma_3\tau_1(i)}*g^1_{\tau_1(i)}*h^1_i\\
B_2=&f_{\sigma_2\tau_2(i)\sigma_3\tau_3(i)}*f_{\sigma_2\tau_1(i)\sigma_2\tau_2(i)}*\Gamma^{12}_{\tau_1(i)}*h^1_i\\&-f_{\sigma_3\tau_1(i)\sigma_3\tau_3(i)}*f_{\sigma_2\tau_1(i)\sigma_3\tau_1(i)}*\Gamma^{12}_{\tau_1(i)}*h^1_i\\
C=& f_{\sigma_2\tau_2(i)\sigma_3\tau_3(i)}*
\Delta^{g^2}_{\tau_1(i)\tau_2(i)}*h^1_i-f_{\sigma_3\tau_1(i)\sigma_3\tau_3(i)}*\Gamma^{23}_{\tau_1(i)}*h^1_i
\end{align*}
We first simplify part $A$. Observe that by Equation~\ref{eq:2-hom} applied to the morphism $g^3$, we get
\begin{align*}
 &g^3_{\tau_3(i)}*\Lambda_{\tau_1(i)\tau_2(i)\tau_3(i)}*h^1_i  -\Delta^{g^3}_{\tau_1(i)\tau_3(i)}*h^1_i=\\
&\big( \Lambda_{\sigma_3\tau_1(i)\sigma_3\tau_2(i)\sigma_3\tau_3(i)}*g^3_{\tau_1(i)}-f_{\sigma_3\tau_2(i)\sigma_3\tau_3(i)}*\Delta^{g^3}_{\tau_1(i)\tau_2(i)}-\Delta^{g^3}_{\tau_2(i)\tau_3(i)}*f_{\tau_1(i)\tau_2(i)}\big)*h^1_i
\end{align*}
The last term of the right side coupled with the term $\Delta^{g^3}_{\tau_2(i)\tau_3(i)}*h^2_i$ in $A$ gives
\begin{align*}
& \Delta^{g^3}_{\tau_2(i)\tau_3(i)}*h^2_i-\Delta^{g^3}_{\tau_2(i)\tau_3(i)}*f_{\tau_1(i)\tau_2(i)}*h^1_i\\
=& \big( -f_{\sigma_3\tau_2(i)\sigma_3\tau_3(i)}*g^3_{\tau_2(i)}+g^3_{\tau_3(i)}*f_{\tau_2(i)\tau_3(i)}\big)*\Upsilon^{12}_i.
\end{align*}
Furthermore, we also have
\begin{align*}
&f_{\sigma_2\tau_2(i)\sigma_3\tau_3(i)}*g^2_{\tau_2(i)}*\Upsilon^{12}_i-f_{\sigma_3\tau_2(i)\sigma_3\tau_3(i)}*g^3_{\tau_2(i)}*\Upsilon^{12}_i\\
=&\big( -f_{\sigma_3\tau_2(i)\sigma_3\tau_3(i)}*f_{\sigma_2\tau_2(i)\sigma_3\tau_2(i)}+f_{\sigma_2\tau_2(i)\sigma_3\tau_3(i)}\big)*g^2_{\tau_2(i)}*\Upsilon^{12}_i
\\&+f_{\sigma_3\tau_2(i)\sigma_3\tau_3(i)}*\big(f_{\sigma_2\tau_2(i)\sigma_3\tau_2(i)}*g^2_{\tau_2(i)}-g^3_{\tau_2(i)}\big)*\Upsilon^{12}_i\\
=&\Lambda_{\sigma_2\tau_2(i)\sigma_3\tau_2(i)\sigma_3\tau_3(i)}*g^2_{\tau_2(i)}*(h^2_i-f_{\tau_1(i)\tau_2(i)}*h^1_i)
\\&+f_{\sigma_3\tau_2(i)\sigma_3\tau_3(i)}*\Gamma^{23}_{\tau_2(i)}*(f_{\tau_1(i)\tau_2(i)}*h^1_i-h^2_i)
\end{align*}
Putting these identities together we arrived at
\begin{align}~\label{eq:A}\begin{split}
A&=\Lambda_{\sigma_3\tau_1(i)\sigma_3\tau_2(i)\sigma_3\tau_3(i)}*g^3_{\tau_1(i)}*h^1_i-f_{\sigma_3\tau_2(i)\sigma_3\tau_3(i)}*\Delta^{g^3}_{\tau_1(i)\tau_2(i)}*h^1_i\\
&-\Lambda_{\sigma_2\tau_2(i)\sigma_3\tau_2(i)\sigma_3\tau_3(i)}*g^2_{\tau_2(i)}*f_{\tau_1(i)\tau_2(i)}*h^1_i
\\&+f_{\sigma_3\tau_2(i)\sigma_3\tau_3(i)}*\Gamma^{23}_{\tau_2(i)}*f_{\tau_1(i)\tau_2(i)}*h^1_i
\end{split}
\end{align}
Secondly, we simplify the terms in $B_1$ using Equation~\ref{eq:cocycle2}. We get
\[ B_1=\big(\Lambda_{\sigma_2\tau_1(i)\sigma_3\tau_1(i)\sigma_3\tau_3(i)}-\Lambda_{\sigma_2\tau_1(i)\sigma_2\tau_2(i)
\sigma_3\tau_3(i)}\big)*f_{\sigma_1\tau_1(i)\sigma_2\tau_1(i)}*g^1_{\tau_1(i)}*h^1_i.\]
Adding this to $B_2$, we have
\begin{align*}
B=&B_1+B_2\\
=&B_1+f_{\sigma_2\tau_2(i)\sigma_3\tau_3(i)}*f_{\sigma_2\tau_1(i)\sigma_2\tau_2(i)}*\Gamma^{12}_{\tau_1(i)}*h^1_i\\&-f_{\sigma_3\tau_1(i)\sigma_3\tau_3(i)}*f_{\sigma_2\tau_1(i)\sigma_3\tau_1(i)}*\Gamma^{12}_{\tau_1(i)}*h^1_i\\
=&B_1+\big(f_{\sigma_2\tau_2(i)\sigma_3\tau_3(i)}*f_{\sigma_2\tau_1(i)\sigma_2\tau_2(i)}-f_{\sigma_2\tau_1(i)\sigma_3\tau_3(i)}\big)*\Gamma^{12}_{\tau_1(i)}*h^1_i-\\
&\big(f_{\sigma_3\tau_1(i)\sigma_3\tau_3(i)}*f_{\sigma_2\tau_1(i)\sigma_3\tau_1(i)}-f_{\sigma_2\tau_1(i)\sigma_3\tau_3(i)}\big)*\Gamma^{12}_{\tau_1(i)}*h^1_i\\
=&B_1+\Lambda_{\sigma_2\tau_1(i)\sigma_2\tau_2(i)\sigma_3\tau_3(i)}*(f_{\sigma_1\tau_1(i)\sigma_2\tau_1(i)}*
g^1_{\tau_1(i)}-g^2_{\tau_1(i)})*h^1_i-\\
&-\Lambda_{\sigma_2\tau_1(i)\sigma_3\tau_1(i)\sigma_3\tau_3(i)}*(f_{\sigma_1\tau_1(i)\sigma_2\tau_1(i)}*
g^1_{\tau_1(i)}-g^2_{\tau_1(i)})*h^1_i\\
=&\big(\Lambda_{\sigma_2\tau_1(i)\sigma_3\tau_1(i)\sigma_3\tau_3(i)}-\Lambda_{\sigma_2\tau_1(i)\sigma_2\tau_2(i)\sigma_3\tau_3(i)}\big)
*g^2_{\tau_1(i)}*h^1_i
\end{align*}
The last term in $A$~\ref{eq:A} and the last term of $C$ gives
\begin{align*}
&f_{\sigma_3\tau_2(i)\sigma_3\tau_3(i)}*\Gamma^{23}_{\tau_2(i)}*f_{\tau_1(i)\tau_2(i)}*h^1_i-f_{\sigma_3\tau_1(i)\sigma_3\tau_3(i)}*\Gamma^{23}_{\tau_1(i)}*h^1_i\\
=& f_{\sigma_3\tau_2(i)\sigma_3\tau_3(i)}*\big(\Gamma^{23}_{\tau_2(i)}*f_{\tau_1(i)\tau_2(i)}-f_{\sigma_3\tau_1(i)\sigma_3\tau_2(i)}*\Gamma^{23}_{\tau_1(i)}\big)*h^1_i\\
&+\big(f_{\sigma_3\tau_2(i)\sigma_3\tau_3(i)}*f_{\sigma_3\tau_1(i)\sigma_3\tau_2(i)}-f_{\sigma_3\tau_1(i)\sigma_3\tau_3(i)}\big)*\Gamma^{23}_{\tau_1(i)}*h^1_i\\
=&f_{\sigma_3\tau_2(i)\sigma_3\tau_3(i)}*\big(-f_{\sigma_2\tau_2(i)\sigma_3\tau_2(i)}*\Delta^{g^2}_{\tau_1(i)\tau_2(i)}
-\Lambda_{\sigma_2\tau_1(i)\sigma_3\tau_1(i)\sigma_3\tau_2(i)}*g^2_{\tau_1(i)}+\\
&+\Lambda_{\sigma_2\tau_1(i)\sigma_2\tau_2(i)\sigma_3\tau_2(i)}*g^2_{\tau_1(i)}+\Delta^{g^3}_{\tau_1(i)\tau_2(i)}\big)*h^1_i\\
&+\Lambda_{\sigma_3\tau_1(i)\sigma_3\tau_2(i)\sigma_3\tau_3(i)}*\big(f_{\sigma_2\tau_1(i)\sigma_3\tau_1(i)}*g^2_{\tau_1(i)}-g^3_{\tau_1(i)}\big)*h^1_i
\end{align*}
In the second equality we used Equation~\ref{eq:2-equiv} applied to $\Gamma^{23}$. Putting everything together, we obtain
\begin{align*} 
&A+B+C=D+E\\
D=&\Big(-f_{\sigma_3\tau_2(i)\sigma_3\tau_3(i)}*\Lambda_{\sigma_2\tau_1(i)\sigma_3\tau_1(i)\sigma_3\tau_2(i)}+
f_{\sigma_3\tau_2(i)\sigma_3\tau_3(i)}*\Lambda_{\sigma_2\tau_1(i)\sigma_2\tau_2(i)\sigma_3\tau_2(i)}+\\
&\Lambda_{\sigma_3\tau_1(i)\sigma_3\tau_2(i)\sigma_3\tau_3(i)}*f_{\sigma_2\tau_1(i)\sigma_3\tau_1(i)}+
\Lambda_{\sigma_2\tau_1(i)\sigma_3\tau_1(i)\sigma_3\tau_3(i)}-\Lambda_{\sigma_2\tau_1(i)\sigma_2\tau_2(i)\sigma_3\tau_3(i)}\Big)
\\&*g^2_{\tau_1(i)}*h^1_i\\
E=&f_{\sigma_2\tau_2(i)\sigma_3\tau_3(i)}*
\Delta^{g^2}_{\tau_1(i)\tau_2(i)}*h^1_i-f_{\sigma_3\tau_2(i)\sigma_3\tau_3(i)}*
f_{\sigma_2\tau_2(i)\sigma_3\tau_2(i)}*\Delta^{g^2}_{\tau_1(i)\tau_2(i)}*h^1_i\\
&-\Lambda_{\sigma_2\tau_2(i)\sigma_3\tau_2(i)\sigma_3\tau_3(i)}*g^2_{\tau_2(i)}*f_{\tau_1(i)\tau_2(i)}*h^1_i
\end{align*}
Using Equation~\ref{eq:cocycle2} to simplify the terms in $D$ we get
\[ D=\Lambda_{\sigma_2\tau_2(i)\sigma_3\tau_2(i)\sigma_3\tau_3(i)}*f_{\sigma_2\tau_1(i)\sigma_2\tau_2(i)}*g^2_{\tau_1(i)}*h^1_i.\]
Adding the last term of $E$ into $D$, we have
\begin{align*}
&\Lambda_{\sigma_2\tau_2(i)\sigma_3\tau_2(i)\sigma_3\tau_3(i)}*f_{\sigma_2\tau_1(i)\sigma_2\tau_2(i)}*g^2_{\tau_1(i)}*h^1_i
\\&-\Lambda_{\sigma_2\tau_2(i)\sigma_3\tau_2(i)\sigma_3\tau_3(i)}*g^2_{\tau_2(i)}*f_{\tau_1(i)\tau_2(i)}*h^1_i\\
=&\Lambda_{\sigma_2\tau_2(i)\sigma_3\tau_2(i)\sigma_3\tau_3(i)}*\big(f_{\sigma_2\tau_1(i)\sigma_2\tau_2(i)}*g^2_{\tau_1(i)}
-g^2_{\tau_2(i)}*f_{\tau_1(i)\tau_2(i)}\big)*h^1_i\\
=&\big(f_{\sigma_3\tau_2(i)\sigma_3\tau_3(i)}*f_{\sigma_2\tau_2(i)\sigma_3\tau_2(i)}-f_{\sigma_2\tau_2(i)\sigma_3\tau_3(i)}\big)*\Delta^{g^2}_{\tau_1(i)\tau_2(i)}*h^1_i
\end{align*}
This cancels precisely the first two terms in $E$. Thus $A+B+C=0$, as desired.\ed

\end{document}